\documentclass[sn-mathphys-num]{sn-jnl}

\usepackage{graphicx}%
\usepackage{multirow}%
\usepackage{amsmath,amssymb,amsfonts}%
\usepackage{amsthm}%
\usepackage{mathrsfs}%
\usepackage[title]{appendix}%
\usepackage{xcolor}%
\usepackage{textcomp}%
\usepackage{manyfoot}%
\usepackage{booktabs}%
\usepackage{algorithm}%
\usepackage{algorithmicx}%
\usepackage{algpseudocode}%
\usepackage{listings}%
\usepackage{tikz}
\usepackage{subfigure}

\theoremstyle{thmstyleone}%
\newtheorem{theorem}{Theorem}
\theoremstyle{thmstyletwo}%
\newtheorem{example}{Example}%
\newtheorem{remark}{Remark}%

\theoremstyle{thmstylethree}%
\newtheorem{definition}{Definition}%

\newtheorem{lemma}{Lemma}
\newtheorem{corollary}{Corollary}
\newtheorem{result}{Result}
\begin{document}

\title[Article Title]{Discrete the solving model of time-variant standard Sylvester-conjugate matrix equations using Euler-forward formula: An analysis of the differences between sampling discretion errors and space compressive approximation errors in optimizing neural dynamics}


\author[1,2]{\fnm{Jiakuang} \sur{He}}\email{hjkq101860089@163.com}

\author*[2]{\fnm{Dongqing} \sur{Wu}}\email{rickwu@zhku.edu.cn}

\affil[1]{School of Clinical Medical College of Acupuncture moxibustion and Rehabilitation, Guangzhou University of Chinese Medicine, Guangzhou 510006, P.
	R. China}
\affil*[2]{School of Mathematics and Data Science, Zhongkai University of Agriculture and Engineering, Guangzhou 510225, P.
	R. China}



\abstract{Time-variant standard Sylvester-conjugate matrix equations are presented as early time-variant versions of the complex conjugate matrix equations. Current solving methods include Con-CZND1 and Con-CZND2 models, both of which use ode45 for continuous model. Given practical computational considerations, discrete these models is also important. Based on Euler-forward formula discretion, Con-DZND1-2i model and Con-DZND2-2i model are proposed. Numerical experiments using step sizes of 0.1 and 0.001. The above experiments show that Con-DZND1-2i model and Con-DZND2-2i model exhibit different neural dynamics compared to their continuous counterparts, such as trajectory correction in Con-DZND2-2i model and the swallowing phenomenon in Con-DZND1-2i model, with convergence affected by step size. These experiments highlight the differences between optimizing sampling discretion errors and space compressive approximation errors in neural dynamics.}




\maketitle

\section{Introduction}
Standard Sylvester-conjugate matrix equations (SSCME)
\cite{bevisConsimilarityMatrixEquation1987}
are the earliest version of complex conjugate matrix equations (CCME)
\cite{Wu2017}.
And SSCME is time-invariant. In recent years of studies, Wu et al. provided matrix algebraic formula methods
\cite{wuSolutionsMatrixEquations2006}
and iterative solving methods
\cite{Wu2017}
based on approximation theory. The essence of the iterative methods is to gradually approach the theoretical solution using multi-step computations. And time-variant standard Sylvester-conjugate matrix equations (TVSSCME)
\cite{He2024ZeroingND} 
is the time-variant extension of SSCME. TVSSCME is supplemented by the differences between differential algebra and linear algebra operations
\cite{ongDeferredCorrectionMethods2020,CHEN2022161}.
The difference between SSCME and TVSSCME solutions is shown in Fig. \ref{fig.dif.sscme.tvsscme}.
Unless otherwise specified, let 
$O$ represent ``null matrix", and only consider the unique theoretical solution $X^{\ast}(\tau)$, same as below.
\begin{figure}[!h]
	\centering 
	\begin{tikzpicture}[node distance=2cm]
		\node (A) {(a)$X(\tau)F-AX(\tau)-C=O$, where $\tau\to +\infty$, $X(\tau)\to X^{\ast}(\tau).$};
		\node (B) [below of=A] {(b)$X(\tau)F(\tau)-A(\tau)X(\tau)-C(\tau)=O$, where $\tau\to +\infty$, $X(\tau)\to X^{\ast}(\tau).$};
		
		\draw[<->] (A) -- (B); 
	\end{tikzpicture}
	\caption{Differences between SSCME(a) and TVSSCME(b).}\label{fig.dif.sscme.tvsscme}
\end{figure}

TVSSCME is currently primarily solved using zeroing neural dynamics (ZND) models Con-CZND1
\cite{He2024ZeroingND}
and Con-CZND2
\cite{He2024ZeroingND}.
Above two models structure can be seen in Fig. \ref{fig.1.2.structure}.
\begin{figure*}[!h]
	\centering  
	\subfigure[]{%
		
		\resizebox{0.4875\textwidth}{!}{%
			\def\layersep{1.5cm}
			
			\begin{tikzpicture}[shorten >=1pt,->,draw=black!50, node distance=\layersep]
			\tikzstyle{every pin edge}=[<-,shorten <=1pt]
			\tikzstyle{neuron}=[circle,fill=black!25,minimum size=17pt,inner sep=0pt]
			\tikzstyle{input neuron}=[neuron, fill=green!50];
			\tikzstyle{input neuron_omit}=[neuron, fill=white!50];
			
			\tikzstyle{hidden neuron1}=[neuron, fill=blue!50];
			\tikzstyle{hidden neuron1_omit}=[neuron, fill=white!50];
			\tikzstyle{hidden neuron2}=[neuron, fill=blue!50];
			\tikzstyle{hidden neuron2_omit}=[neuron, fill=white!50];
			\tikzstyle{hidden neuron23}=[neuron, fill=white!50];
			\tikzstyle{hidden neuron3}=[neuron, fill=green!50];
			\tikzstyle{hidden neuron3_omit}=[neuron, fill=white!50];
			\tikzstyle{hidden neuron4}=[neuron, fill=red!50];
			\tikzstyle{hidden neuron4_omit}=[neuron, fill=white!50];
			\tikzstyle{annot} = [text width=4em, text centered]
			
			\foreach \y in {1, 3,4, 6} {
				\node[input neuron, pin=left:Random Input] (I-\y) at (0,-\y) {};
			}
			\foreach \y in {2, 5} {
				\node[input neuron_omit,label={left: {\Huge \textbf{n\{}}}] (I-\y) at (0,-\y) {$\vdots$};
			}
			
			\foreach \name / \y in {1,3}
			\path[yshift=-6cm]
			node[hidden neuron1] (H1-\name) at (\layersep,-\y cm) {};
			\foreach \name / \y in {2}
			\path[yshift=-6cm]
			node[hidden neuron1_omit,label={left: {\Huge \textbf{n\{}}}] (H1-\name) at (\layersep,-\y cm) {$\vdots$};
			\foreach \name / \y in {1,3}
			\path[yshift=-6cm]
			node[hidden neuron2,right of=H1] (H2-\name) at (\layersep,-\y cm){};
			\foreach \name / \y in {2}
			\path[yshift=-6cm]
			node[hidden neuron2_omit,right of=H1] (H2-\name) at (\layersep,-\y cm){$\vdots$};
			
			\foreach \name / \y in {1}
			\path[yshift=-7cm]
			node[hidden neuron23,right of=H2] (H23-\name) at (0.5*\layersep,-\y cm){...};
			
			\foreach \name / \y in {1,3,4,6}
			\path[yshift=0cm]
			node[hidden neuron3,right of=H2] (H3-\name) at (2*\layersep,-\y cm){};
			
			\foreach \name / \y in {2,5}
			\path[yshift=0cm]
			node[hidden neuron3_omit,right of=H2] (H3-\name) at (2*\layersep,-\y cm){$\vdots$};
			
			\foreach \y in {1, 3, 4, 6} {
				\path[yshift=0cm]
				node[hidden neuron4, right of=H3, pin={[pin edge={->}]right:Output}] (H4-\y) at (3*\layersep,-\y cm) {};
			}
			
			\foreach \y in {2, 5} {
				\path[yshift=0cm]
				node[hidden neuron4_omit, right of=H3,label={left: {\Huge \textbf{n\{}}}] (H4-\y) at (3*\layersep,-\y cm) {$\vdots$};
			}
			
			\foreach \source in {1,3}
			\path (I-\source) edge (H1-\source);
			\path (I-4) edge (H1-1);
			\path (I-6) edge (H1-3);
			
			\foreach \source in {1,3}
			\path (H1-\source) edge (H23-1);
			
			\foreach \dest in {1,3}
			\path (H23-1) edge (H2-\dest);
			
			\foreach \source in {1,3}
			\path (H2-\source) edge (H3-\source);
			\path (H2-1) edge (H3-4);
			\path (H2-3) edge (H3-6);

			\foreach \source in {1,3,4,6}
			\path (H3-\source) edge (H4-\source);

			\foreach \source in {1,3,4,6} {
				\path (H3-\source) edge (I-1);
				\path (H3-\source) edge (I-3);
				\path (H3-\source) edge (I-4);
				\path (H3-\source) edge (I-6);
			}

			
			
			\node[annot,below of=H1-1, node distance=3cm] (hl) {Hidden Complex Layer1};
			\node[annot,above of=I-1, node distance=1cm](IRltext) {Input Real Layer};
			\node[annot,right of=hl] (h2){Hidden Complex Layer2};
			\node[annot,right of=IRltext, node distance=4.5cm] (h3){Hidden Real Layer1};
			\node[annot,right of=h3] {Output Real Layer};
		\end{tikzpicture}
		}%
		
		\label{fig.1.structure}
	} 
	\hfill 
	\subfigure[]{%
		\resizebox{0.4875\textwidth}{!}{%
			\def\layersep{1.5cm}
			\raisebox{2.45cm}{
	\begin{tikzpicture}[shorten >=1pt,->,draw=black!50, node distance=\layersep]
	\tikzstyle{every pin edge}=[<-,shorten <=1pt]
	\tikzstyle{neuron}=[circle,fill=black!25,minimum size=17pt,inner sep=0pt]
	\tikzstyle{input neuron}=[neuron, fill=green!50];
	\tikzstyle{input neuron_omit}=[neuron, fill=white!50];
	
	\tikzstyle{hidden neuron12}=[neuron, fill=white!50];
	\tikzstyle{hidden neuron2}=[neuron, fill=green!50];
	\tikzstyle{hidden neuron2_omit}=[neuron, fill=white!50];
	\tikzstyle{hidden neuron4}=[neuron, fill=red!50];
	\tikzstyle{hidden neuron4_omit}=[neuron, fill=white!50];
	\tikzstyle{annot} = [text width=4em, text centered]
	
	\foreach \y in {1, 3,4, 6} {
		\node[input neuron, pin=left:Random Input] (I-\y) at (0,-\y) {};
	}
	\foreach \y in {2, 5} {
		\node[input neuron_omit,label={left: {\Huge \textbf{n\{}}}] (I-\y) at (0,-\y) {$\vdots$};
	}

	\foreach \name / \y in {1}
	\path[yshift=-7cm]
	node[hidden neuron12,right of=I] (H12-\name) at (0.5*\layersep,-\y cm){...};
	
	\foreach \y in {1, 3, 4, 6} {
		\path[yshift=0cm]
		node[hidden neuron2, right of=H2] (H2-\y) at (2*\layersep,-\y cm) {};
	}
	
	\foreach \y in {2, 5} {
		\path[yshift=0cm]
		node[hidden neuron2_omit, right of=H2] (H2-\y) at (2*\layersep,-\y cm) {$\vdots$};
	}
	
	\foreach \y in {1, 3, 4, 6} {
		\path[yshift=0cm]
		node[hidden neuron4, right of=H3, pin={[pin edge={->}]right:Output}] (H3-\y) at (3*\layersep,-\y cm) {};
	}
	
	\foreach \y in {2, 5} {
		\path[yshift=0cm]
		node[hidden neuron4_omit, right of=H3,label={left: {\Huge \textbf{n\{}}}] (H3-\y) at (3*\layersep,-\y cm) {$\vdots$};
	}
	
	\foreach \source in {1, 3, 4, 6} {
		\path (I-\source) edge (H12-1);
	}
	
	\foreach \dest in {1, 3, 4, 6} {
		\path (H12-1) edge (H2-\dest);
	}
	
	\foreach \source in {1,3,4,6} {
		\path (H2-\source) edge (I-1);
		\path (H2-\source) edge (I-3);
		\path (H2-\source) edge (I-4);
		\path (H2-\source) edge (I-6);
	}

	\foreach \source in {1,3,4,6}		
	\path (H2-\source) edge (H3-\source);
	
	\node[annot,above of=I-1, node distance=1cm] (hl) {Input Real Layer};
	\node[annot,right of=hl] (h2){};
	\node[annot,right of=h2] (h3){};
	\node[annot,right of=h3] (h4){Hidden Real Layer1};
	\node[annot,right of=h4] {Output Real Layer};
\end{tikzpicture}}
		}%
		\label{fig.2.structure}
	}  
	\caption{Different between Con-CZND1 \cite{He2024ZeroingND} model and Con-CZND2 \cite{He2024ZeroingND} model. \subref{fig.1.structure} Con-CZND1 model. \subref{fig.2.structure} Con-CZND2 model.} 
	\label{fig.1.2.structure} 
\end{figure*}

However, Con-CZND1 model essentially approximates using the complex field error, while Con-CZND2 model approximates using the real field error. In ode45
\cite{dabounouAdaptiveFeedforwardGradient2024}
solver, Con-CZND2 model does not perform as well as Con-CZND1 model.

Discrete neural dynamics is validated in previous studies to reduce the error between theoretical and numerical solutions
\cite{qiuAccuracyArchitectureStudies2023}.
Zhang et al. continued to develop discretion in the real field, progressing from Euler-forward formula
\cite{chenContinuousDiscreteZeroing2020,rodgersAdaptiveIntegrationNonlinear2022}
to an 11-point sampling discretion
\cite{10355921,uhligZhangNeuralNetworks2024,wuGaussNewtonMethod2024}.
However, there is no exploration of neural models for solving TVSSCME using sampling discretion in the existing literature. According to the known studies, the two continuous solution models, Con-CZND1 and Con-CZND2, show significant differences due to the approximation effects of the internal ode45
\cite{flynn2024exploringoriginsswitchingdynamics}
solver. Additionally, Con-CZND1 model exhibits space compressive approximation phenomenon. Therefore, it is essential to rigorously establish a discrete neural dynamics model for TVSSCME.

The rest of this article is organized as follows:
Section 2 provides the definition of TVSSCME and supplementary knowledge. Section 3 defines the Con-DZND1-2i discrete solving model over the complex field and the Con-DZND2-2i discrete solving model over the real field.
Section 4 presents simulations that validate the effectiveness of each model and compares their strengths and weaknesses.
Sections 5 and 6 summarizes this article and suggests future directions. Before proceeding to the next section, the main contributions of this article are as follows:

\begin{itemize}
	\item[(1)]Con-DZND1-2i model, which directly defines complex field error, and Con-DZND2-2i model, which maps to real field error, are proposed for solving TVSSCME. 
	
	\item[(2)]Based on Euler-forward formula, both discrete models, Con-DZND1-2i and Con-DZND2-2i, which use different step sizes, can ultimately approximate the theoretical solution.
	
	\item[(3)]Con-DZND1-2i model defines complex field error, while Con-DZND2-2i model maps to real field error. These models highlight a significant difference between optimizing space compressive approximation errors and optimizing sampling discretion errors in neural network optimization. Both aspects should be considered from different perspectives.
\end{itemize}

\section{Problem formulation}

Consider a TVSSCME\eqref{eq.sccsme.variant} 
\cite{He2024ZeroingND,He2024RevisitingTC}
as
\begin{equation} \label{eq.sccsme.variant}
	X(\tau)F(\tau)-A(\tau)\overline{X}(\tau)-C(\tau)=O,
\end{equation}
where $F(\tau)\in
\mathbb{C}^{n\times n}$, $A(\tau)\in
\mathbb{C}^{m\times m}$, $C(\tau)\in
\mathbb{C}^{m\times n}$ are known as time-variant matrices,
$X(\tau)\in
\mathbb{C}^{m\times n}$ is a time-variant matrix to be computed,
and $\tau \ge 0$ represents the real-time.

To facilitate subsequent derivations, supplementary definitions and relevant formulas are provided.
\begin{definition}
	\cite{needhamVisualComplexAnalysis2023,Krmer2024ATO,4203075,He2024RevisitingTC}
The time-variant matrix elements $\tilde{m}_{st}(\tau)$ are defined as follows:
	\begin{equation}\label{eq.define.complexmatrix.elements}
		\tilde{m}_{st}(\tau)=m_{\mathrm{r},st}(\tau)+\mathrm{i}m_{\mathrm{i},st}(\tau),
	\end{equation}
where $s\in
\mathbb{I}[1,p]$, $t\in
\mathbb{I}[1,q]$, $\mathbb{I}[m,n]$ means the set of integers from m to n, $m_{\mathrm{r},st}(\tau)\in
\mathbb{R}$ is a real coefficient, $m_{\mathrm{i},st}(\tau)\in
\mathbb{R}$ is a imaginary coefficient, $\mathrm{i}$ is an imaginary unit, same as below. Also, if $m_{\mathrm{i},st}(\tau)=0$, the elements $\tilde{m}_{st}(\tau)=\left (m_{\mathrm{r},st}(\tau)+\mathrm{i}\times0\right )\in
\mathbb{R}$. For the sake of uniformity, the elements are represented by $m_{st}(\tau)=\tilde{m}_{st}(\tau)$ under the real matrices.
	
Correspondingly, the time-variant matrix $M(\tau)$ and time-variant conjugate matrix $\overline{M}(\tau)$ are defined as follows:
	\begin{equation}\label{eq.define.complexmatrix}
		M(\tau)=M_{\mathrm{r}}(\tau)+\mathrm{i}M_{\mathrm{i}}(\tau),
	\end{equation}
\begin{equation}\label{eq.define.complexmatrix.conj}
	\overline{M}(\tau)=M_{\mathrm{r}}(\tau)-\mathrm{i}M_{\mathrm{i}}(\tau).
\end{equation}	
\end{definition}
\begin{lemma}\cite{He2024ZeroingND,He2024RevisitingTC}
	Where $\tau \ge 0$ represents the real-time, $A(\tau)\in
	\mathbb{C}^{m\times n}$, $B(\tau)\in
	\mathbb{C}^{s\times t}$, $X(\tau)\in
	\mathbb{C}^{n\times s}$, are time-variant matrices,
	the following equation can be obtained:
	\begin{equation}\label{eq.infer.complexkroneckerproductvectorization} 
		\mathrm{vec}(A(\tau)X(\tau)B(\tau))=(\overline{B^{\mathrm{H}}}(\tau)\otimes A(\tau))\mathrm{vec}(X(\tau)).	
	\end{equation}
	\begin{corollary}
		Where $A(\tau)\in
		\mathbb{R}^{m\times n}$, $B(\tau)\in
		\mathbb{R}^{s\times t}$, $X(\tau)\in
		\mathbb{R}^{n\times s}$, and $\tau \ge 0$ represents the real-time, \eqref{eq.infer.complexkroneckerproductvectorization} converts to: \begin{equation}\label{eq.infer.realkroneckerproductvectorization} 
			\mathrm{vec}(A(\tau)X(\tau)B(\tau))=(B^{\mathrm{T}}(\tau)\otimes A(\tau))\mathrm{vec}(X(\tau)).	
		\end{equation}
	\end{corollary}
\end{lemma}
\subsection{Con-CZND1 model}
The following error of \eqref{eq.sccsme.variant} is then defined:
\begin{equation} \label{eq.define.errconcznd1}
		E_{\mathrm{C}}(\tau)=X(\tau)F(\tau)-A(\tau)\overline{X}(\tau)-C(\tau),
\end{equation}
where $E_{\mathrm{C}}(\tau)\in\mathbb{C}^{m\times n}$.
Using the complex field ZND
\cite{He2024ZeroingND,He2024RevisitingTC}, that is:
\begin{subequations} \label{eq.infer.linearerrconcznd1}
	\begin{align}
\frac{\partial E_{\mathrm{C}}(\tau)}{\partial \tau} +\gamma E_{\mathrm{C}}(\tau)=O,\\
		\dot{\tilde{m}}_{st}(\tau)+\gamma  \tilde{m}_{st}(\tau)=0,
	\end{align}
\end{subequations}
where $\gamma\in\mathbb{R^+}$ denotes the regulation parameter controlling the convergence rate, $\dot{\tilde{m}}_{st}(\tau)\in\mathbb{C}$ is $E_{\mathrm{C}}(\tau)$ elements differentiated from $\tau$.
Therefore, substituting \eqref{eq.define.errconcznd1} into \eqref{eq.infer.linearerrconcznd1} results in \eqref{eq.join.linearerrconcznd1}.
\begin{equation}\label{eq.join.linearerrconcznd1}
	\frac{\left ( X(\tau)F(\tau)-A(\tau)\overline{X}(\tau)-C(\tau)\right )}{\partial \tau} +\gamma \left ( X(\tau)F(\tau)-A(\tau)\overline{X}(\tau)-C(\tau) \right )=O.
\end{equation}
Using \eqref{eq.infer.complexkroneckerproductvectorization}, \eqref{eq.join.linearerrconcznd1} is converted to \eqref{eq.use.complexkroneckerproductvectorization}:
\begin{equation} \label{eq.use.complexkroneckerproductvectorization}
	\begin{split}
	(\overline{F^{\mathrm{H}}}(\tau)\otimes I_{m})\mathrm{vec}(\dot{X}(\tau))-(\overline{I_{n}^{\mathrm{H}}}\otimes A(\tau))\mathrm{vec}(\dot{\overline{X}}(\tau))
	\\
	-\mathrm{vec}(\dot{C}(\tau)+\dot{A}(\tau)\overline{X}(\tau)-X(\tau)\dot{F}(\tau))
	\\
	+\gamma \mathrm{vec} (X(\tau)F(\tau)-A(\tau)\overline{X}(\tau)-C(\tau))=O.
\end{split}
\end{equation}
Then, \eqref{eq.use.complexkroneckerproductvectorization} is further reformulated as
\begin{equation}
	\label{eq.simplify.complexkroneckerproductvectorization}
	U(\tau)\mathrm{vec}(\dot{X}(\tau))-V(\tau)\mathrm{vec}(\dot{\overline{X}}(\tau))
	-G(\tau)=O,
\end{equation}
where $U(\tau)=(\overline{F^{\mathrm{H}}}(\tau) \otimes I_{m})\in
\mathbb{C}^{nm\times mn}$, $V(\tau)=(\overline{I_{n}^{\mathrm{H}}} \otimes A(\tau))=(I_{n}\otimes A(\tau))\in
\mathbb{C}^{mn\times nm}$,
$G(\tau)=\mathrm{vec}(\dot{C}(\tau)+\dot{A}(\tau)\overline{X}(\tau)-X(\tau)\dot{F}(\tau))
-\gamma \mathrm{vec} (X(\tau)F(\tau)-A(\tau)\overline{X}(\tau)-C(\tau))\in
\mathbb{C}^{mn\times 1}$. Based on the linearity of the derivative as well as \eqref{eq.define.complexmatrix} and \eqref{eq.define.complexmatrix.conj}, \eqref{eq.simplify.complexkroneckerproductvectorization} can be written in the form of the following real-only matrix operation:
\begin{align}\label{eq.divide.complexkroneckerproductvectorization}
	\begin{bmatrix}
		U_{\mathrm{r}}(\tau)-V_{\mathrm{r}}(\tau)	&-(U_{\mathrm{i}}(\tau)+V_{\mathrm{i}}(\tau)) \\
		U_{\mathrm{i}}(\tau)-V_{\mathrm{i}}(\tau)	&U_{\mathrm{r}}(\tau)+V_{\mathrm{r}}(\tau)
	\end{bmatrix}
	\begin{bmatrix}
		\dot{Z}_{\mathrm{r}}(\tau)\\
		\dot{Z}_{\mathrm{i}}(\tau)	
	\end{bmatrix}
	-\begin{bmatrix}
		G_{\mathrm{r}}(\tau)\\
		G_{\mathrm{i}}(\tau)	
	\end{bmatrix}=O,
\end{align}
where $Z(\tau)=\mathrm{vec}(X(\tau))\in
\mathbb{C}^{mn\times 1}$, $\dot{Z}(\tau)=\mathrm{vec}(\dot{X}(\tau))\in
\mathbb{C}^{mn\times 1}$. To simplify, let $W_{\mathrm{C}}(\tau)=\left [U_{\mathrm{r}}(\tau)-V_{\mathrm{r}}(\tau), -(U_{\mathrm{i}}(\tau)+V_{\mathrm{i}}(\tau));
U_{\mathrm{i}}(\tau)-V_{\mathrm{i}}(\tau),	U_{\mathrm{r}}(\tau)+V_{\mathrm{r}}(\tau)\right ]
\in\mathbb{R}^{2mn\times 2mn}$, $\dot{X}_{\mathrm{C}}(\tau)= \left [\dot{Z}_{\mathrm{r}}(\tau);
\dot{Z}_{\mathrm{i}}(\tau)\right ]\in
\mathbb{R}^{2mn\times 1}$, $B_{\mathrm{C}}(\tau)= \left [G_{\mathrm{r}}(\tau);
G_{\mathrm{i}}(\tau)\right ]\in
\mathbb{R}^{2mn\times 1}$. The final solution model Con-CZND1 is obtained:
\begin{align} \label{eq.solve.linearerrconcznd1}
	\dot{X}_{\mathrm{C}}(\tau) -W^{+}_{\mathrm{C}}(\tau)B_{\mathrm{C}}(\tau)=O,
\end{align}
where $W^{+}_{\mathrm{C}}(\tau)$ is the pseudo-inverse matrix of $W_{\mathrm{C}}(\tau)$.
\subsection{Con-CZND2 model}
Using \eqref{eq.define.complexmatrix}, \eqref{eq.sccsme.variant} is first transformed into \eqref{eq.transccsme.variant}. 
\begin{equation}\label{eq.transccsme.variant}
	\begin{split}
			(X_{\mathrm{r}}(\tau)+\mathrm{i}X_{\mathrm{i}}(\tau))(F_{\mathrm{r}}(\tau)+\mathrm{i}F_{\mathrm{i}}(\tau))-(A_{\mathrm{r}}(\tau)+\mathrm{i}A_{\mathrm{i}}(\tau))(X_{\mathrm{r}}(\tau)-\mathrm{i}X_{\mathrm{i}}(\tau))
		\\
		-(C_{\mathrm{r}}(\tau)+\mathrm{i}C_{\mathrm{i}}(\tau))=O.
	\end{split}
\end{equation}
Then, \eqref{eq.transccsme.variant} is performed separation to obtain \eqref{eq.dividetransccsme.variant} with real-only matrix operation:
\begin{equation}\label{eq.dividetransccsme.variant}
	\left\{\begin{matrix}
		X_{\mathrm{r}}(\tau)F_{\mathrm{r}}(\tau)-X_{\mathrm{i}}(\tau)F_{\mathrm{i}}(\tau)-A_{\mathrm{r}}(\tau)X_{\mathrm{r}}(\tau)-A_{\mathrm{i}}(\tau)X_{\mathrm{i}}(\tau)-C_{\mathrm{r}}(\tau)=O,	\\
		X_{\mathrm{i}}(\tau)F_{\mathrm{r}}(\tau)+X_{\mathrm{r}}(\tau)F_{\mathrm{i}}(\tau)-A_{\mathrm{i}}(\tau)X_{\mathrm{r}}(\tau)+A_{\mathrm{r}}(\tau)X_{\mathrm{i}}(\tau)-C_{\mathrm{i}}(\tau)=O.	
	\end{matrix}\right.
\end{equation}
According to \eqref{eq.infer.realkroneckerproductvectorization}, \eqref{eq.dividetransccsme.variant} is formulated into the following:
\begin{equation}\label{eq.veckrondividetransccsme.variant}
	\begin{bmatrix}
		K_{11}(\tau)	&K_{12}(\tau) \\
		K_{21}(\tau)	&K_{22}(\tau)
	\end{bmatrix}
	\begin{bmatrix}
		\mathrm{vec}(X_{\mathrm{r}}(\tau))\\
		\mathrm{vec}(X_{\mathrm{i}}(\tau))	
	\end{bmatrix}
	-\begin{bmatrix}
		\mathrm{vec}(C_{\mathrm{r}}(\tau))\\
		\mathrm{vec}(C_{\mathrm{i}}(\tau))	
	\end{bmatrix}=O,
\end{equation}
where $K_{11}(\tau)=(F_{\mathrm{r}}^{\mathrm{T}}(\tau) \otimes I_{m})-(I_{n} \otimes A_{\mathrm{r}}(\tau))\in
\mathbb{R}^{nm\times mn}$, $K_{12}(\tau)=-(F_{\mathrm{i}}^{\mathrm{T}}(\tau) \otimes I_{m}+I_{n} \otimes A_{\mathrm{i}}(\tau))\in
\mathbb{R}^{nm\times mn}$,
$K_{21}(\tau)=(F_{\mathrm{i}}^{\mathrm{T}}(\tau) \otimes I_{m})-(I_{n} \otimes A_{\mathrm{i}}(\tau))\in
\mathbb{R}^{nm\times mn}$,
$K_{22}(\tau)=(F_{\mathrm{r}}^{\mathrm{T}}(\tau) \otimes I_{m}+I_{n}\otimes A_{\mathrm{r}}(\tau))\in
\mathbb{R}^{nm\times mn}$. Then let $W_{\mathrm{R}}(\tau)=[K_{11}(\tau),K_{12}(\tau);
K_{21}(\tau),K_{22}(\tau)]\in
\mathbb{R}^{2mn\times 2mn}$, $X_{\mathrm{R}}(\tau)=[\mathrm{vec}(X_{\mathrm{r}}(\tau));
\mathrm{vec}(X_{\mathrm{i}}(\tau))]\in\mathbb{R}^{2mn\times 1}$,
$B_{\mathrm{R}}(\tau)=[\mathrm{vec}(C_{\mathrm{r}}(\tau)); \mathrm{vec}(C_{\mathrm{i}}(\tau))]\in\mathbb{R}^{2mn\times 1}$,
\eqref{eq.time.linearerrconcznd2} is obtained:
\begin{equation} \label{eq.time.linearerrconcznd2}
	W_{\mathrm{R}}(\tau)X_{\mathrm{R}}(\tau) -B_{\mathrm{R}}(\tau)=O.
\end{equation}
Using real field ZND
\cite{He2024ZeroingND,He2024RevisitingTC},
the following error of \eqref{eq.time.linearerrconcznd2} is first defined as follows.
\begin{equation} \label{eq.define.errconcznd2}
	E_{\mathrm{R}}(\tau)=W_{\mathrm{R}}(\tau)X_{\mathrm{R}}(\tau)-B_{\mathrm{R}}(\tau),
\end{equation}
where $E_{\mathrm{R}}(\tau)\in\mathbb{R}^{2mn\times 1}$.
Next, the formula under the real field of ZND is proposed to make all elements of \eqref{eq.define.errconcznd2} 
converge to zero, which is obtained as
\begin{subequations} \label{eq.infer.linearerrconcznd2}
	\begin{align}
		\frac{\partial E_{\mathrm{R}}(\tau)}{\partial \tau} +\gamma E_{\mathrm{R}}(\tau)=O,\\
		\dot{m}_{st}(\tau)+\gamma  m_{st}(\tau)=0,
	\end{align}
\end{subequations}
where $\gamma\in\mathbb{R^+}$ denotes the regulation parameter controlling the convergence rate, $\dot{m}_{st}(\tau)\in\mathbb{R}$ is $E_{\mathrm{R}}(\tau)$ elements differentiated from $\tau$.
Therefore, substituting \eqref{eq.define.errconcznd2} into \eqref{eq.infer.linearerrconcznd2} results in \eqref{eq.join.linearerrconcznd2}.
\begin{equation}\label{eq.join.linearerrconcznd2}
	\frac{\left ( W_{\mathrm{R}}(\tau)X_{\mathrm{R}}(\tau)-B_{\mathrm{R}}(\tau)\right )}{\partial \tau} +\gamma \left ( W_{\mathrm{R}}(\tau)X_{\mathrm{R}}(\tau)-B_{\mathrm{R}}(\tau) \right )=O.
\end{equation}
Finally, the final solution model Con-CZND2 \eqref{eq.solve.linearerrconcznd2} is obtained:
\begin{equation} \label{eq.solve.linearerrconcznd2}
	\begin{split}
		\dot{X}_{\mathrm{R}}(\tau)
		-W^{+}_{\mathrm{R}}(\tau)(\dot{B}_{\mathrm{R}}(\tau)-\dot{W}_{\mathrm{R}}(\tau)X_{\mathrm{R}}(\tau)
		-\gamma(W_{\mathrm{R}}(\tau)X_{\mathrm{R}}(\tau) -B_{\mathrm{R}}(\tau)))=O,
	\end{split}
\end{equation}
where $W^{+}_{\mathrm{R}}(\tau)$ is the pseudo-inverse matrix of $W_{\mathrm{R}}(\tau)$.
\begin{theorem}[Convergence theorem]\label{thm.1}
	Given differentiable time-variant matrices $F(\tau)\in
	\mathbb{C}^{n\times n}$, $A(\tau)\in
	\mathbb{C}^{m\times m}$, and $C(\tau)\in
	\mathbb{C}^{m\times n}$, if TVSSCME \eqref{eq.sccsme.variant} only has one theoretical time-variant solution $X^{\ast}(\tau)\in
	\mathbb{C}^{m\times n}$, then each solving element of \eqref{eq.solve.linearerrconcznd1} and \eqref{eq.solve.linearerrconcznd2} converge to the corresponding theoretical time-variant solving elements.
	\begin{proof}
		Refer to Theorems 2 and 3 in 
		\cite{He2024ZeroingND}. 
	\end{proof} 
\end{theorem}
\section{Euler-forward formula, Con-DZND1-2i and Con-DZND2-2i}
In this section, to facilitate the implementation of mathematical hardware, Euler-forward formula
\cite{VILLATORO1999111,youFastConstructionForward2020}
is applied to discrete \eqref{eq.solve.linearerrconcznd1} and \eqref{eq.solve.linearerrconcznd2}.

Generally, Euler-forward formula
\cite{VILLATORO1999111,youFastConstructionForward2020}
 can be expressed as:
\begin{equation}\label{eq.euler.forward}
	\dot{x}(\tau )=\frac{x(\tau_{k+1})-x(\tau_{k})}{\varepsilon} +\mathcal{O}(\varepsilon), 
\end{equation}
where $\varepsilon$ represents different step sizes, $\mathcal{O}(\varepsilon)$ represents the corresponding truncation residual, same as below.

By substituting \eqref{eq.euler.forward} into \eqref{eq.solve.linearerrconcznd1} and \eqref{eq.solve.linearerrconcznd2}, two equivalent \eqref{eq.euler.forward.linearerrconcznd1} and \eqref{eq.euler.forward.linearerrconcznd2} are obtained.
\begin{equation}\label{eq.euler.forward.linearerrconcznd1}
	X_{\mathrm{C}}(\tau_{k+1})=X_{\mathrm{C}}(\tau_{k})+{\varepsilon}\times W^{+}_{\mathrm{C}}(\tau_{k})B_{\mathrm{C}}(\tau_{k}) +\mathcal{O}(\varepsilon^{2}), 
\end{equation}
\begin{equation}\label{eq.euler.forward.linearerrconcznd2}
	\begin{split}
			X_{\mathrm{R}}(\tau_{k+1})=X_{\mathrm{R}}(\tau_{k})+
			{\varepsilon} \times W^{+}_{\mathrm{R}}(\tau_{k})(\dot{B}_{\mathrm{R}}(\tau_{k})-\dot{W}_{\mathrm{R}}(\tau_{k})X_{\mathrm{R}}(\tau_{k})
			\\
		-\gamma(W_{\mathrm{R}}(\tau_{k})X_{\mathrm{R}}(\tau_{k}) -B_{\mathrm{R}}(\tau_{k})))+\mathcal{O}(\varepsilon^{2}).
	\end{split}
\end{equation}
In \eqref{eq.euler.forward.linearerrconcznd1} and \eqref{eq.euler.forward.linearerrconcznd2}, 
$\varepsilon$ represents the sampling interval with different step sizes, and 
$\mathcal{O}(\varepsilon^{2})$ represents that each component of the vector is $\mathcal{O}(\varepsilon^{2})$. Then \eqref{eq.euler.forward.linearerrconcznd1} and \eqref{eq.euler.forward.linearerrconcznd2} are similar to ``Gradient descent"
\cite{10.5555/2987994},
as seen in Appendix \ref{appendix.A}. Finally, the two major discretion models, Con-DZND1-2i \eqref{eq.euler.forward.solve.linearerrconcznd1} and Con-DZND2-2i \eqref{eq.euler.forward.solve.linearerrconcznd2}, are obtained with a truncation residual of $\mathcal{O}(\varepsilon^{2})$.
\begin{equation}\label{eq.euler.forward.solve.linearerrconcznd1}
	X_{\mathrm{C}}(\tau_{k+1})=X_{\mathrm{C}}(\tau_{k})+{\varepsilon}\times W^{+}_{\mathrm{C}}(\tau_{k})B_{\mathrm{C}}(\tau_{k}), 
\end{equation}
\begin{equation}\label{eq.euler.forward.solve.linearerrconcznd2}
	\begin{split}
		X_{\mathrm{R}}(\tau_{k+1})=X_{\mathrm{R}}(\tau_{k})+{\varepsilon}\times W^{+}_{\mathrm{R}}(\tau_{k})(\dot{B}_{\mathrm{R}}(\tau_{k})-\dot{W}_{\mathrm{R}}(\tau_{k})X_{\mathrm{R}}(\tau_{k})
		\\-\gamma(W_{\mathrm{R}}(\tau_{k})X_{\mathrm{R}}(\tau_{k}) -B_{\mathrm{R}}(\tau_{k}))).
	\end{split}
\end{equation}
\begin{theorem}[Stability theorem]\label{thm.2}
	Let $\mathcal{O}(\varepsilon^{2})$ represent that each component of the vector is $\mathcal{O}(\varepsilon^{2})$, where 
	$\varepsilon\in \left ( 0,1 \right )$ represents the sampling interval. Then, Con-DZND1-2i \eqref{eq.euler.forward.solve.linearerrconcznd1} model and Con-DZND2-2i \eqref{eq.euler.forward.solve.linearerrconcznd2} model are 0-stable, consistent, and convergent, with a truncation residual of 
	$\mathcal{O}(\varepsilon^{2})$.
	\begin{proof}
		The characteristic polynomials corresponding to Con-DZND1-2i \eqref{eq.euler.forward.solve.linearerrconcznd1} model and Con-DZND2-2i \eqref{eq.euler.forward.solve.linearerrconcznd2} are
		\begin{equation}\label{eq.euler.forward.root}
			P_{1}(\delta)=\delta-1.
		\end{equation}
		Then, \eqref{eq.euler.forward.root} has only one root $\delta=1$, which is located on the unit circle. Therefore, according to Result \ref{result.1} in Appendix \ref{appendix.B}, Con-DZND1-2i \eqref{eq.euler.forward.solve.linearerrconcznd1} model and Con-DZND2-2i \eqref{eq.euler.forward.solve.linearerrconcznd2} model are 0-stable.
		
		Furthermore, from \eqref{eq.euler.forward.linearerrconcznd1} and \eqref{eq.euler.forward.linearerrconcznd2}, it can be seen that Con-DZND1-2i \eqref{eq.euler.forward.solve.linearerrconcznd1} model and Con-DZND2-2i \eqref{eq.euler.forward.solve.linearerrconcznd2} model have a truncation residual of 
		$\mathcal{O}(\varepsilon^{2})$.
		
		Thus, based on Results \ref{result.2}, \ref{result.3}, and \ref{result.4} in Appendix \ref{appendix.B}, Con-DZND1-2i \eqref{eq.euler.forward.solve.linearerrconcznd1} model and Con-DZND2-2i \eqref{eq.euler.forward.solve.linearerrconcznd2} model are consistent and convergent, with a truncation residual of 
		$\mathcal{O}(\varepsilon^{2})$.
	\end{proof} 
\end{theorem}
\begin{theorem}[Residual theorem]\label{thm.3}
	Let $\left \|\cdot   \right \|_{\mathrm{F}}$ represents Frobenius norm, while 
	$\left \|\cdot   \right \|_{\mathrm{F}} \ge 0$. And
	$\varepsilon\in \left ( 0,1 \right )$ represents the sampling interval. The maximum steady-state residuals for the two models, Con-DZND1-2i \eqref{eq.euler.forward.solve.linearerrconcznd1} and Con-DZND2-2i \eqref{eq.euler.forward.solve.linearerrconcznd2}, are both 
	$\mathcal{O}(\varepsilon^{2})$.
	\begin{proof}
	As the real field is a subset of the complex field, this article only proves Con-DZND1-2i \eqref{eq.euler.forward.solve.linearerrconcznd1} model based on complex field. The proof for the Con-DZND2-2i \eqref{eq.euler.forward.solve.linearerrconcznd2} model based on the real field is similar, so it will not be discussed further.
	
	Before the proof, the following two norm formulas are provided:
	\begin{lemma}\cite{hornMatrixAnalysis2010}
		Norm equality.
		\begin{equation}\label{eq.frobenius.conj} 
			\left \|\overline{X}(\tau)\right \|_{\mathrm{F}}=\left \|X(\tau)\right \|_{\mathrm{F}}.	
		\end{equation}
	\end{lemma}	
	\begin{lemma}\cite{hornMatrixAnalysis2010}
		Norm triangle inequality.
		\begin{equation}\label{eq.frobenius.triangle} 
				\left \|A(\tau)\right \|_{\mathrm{F}}-\left \|B(\tau)\right \|_{\mathrm{F}}\le\left \|A(\tau)-B(\tau)\right \|_{\mathrm{F}} \le \left \|A(\tau)\right \|_{\mathrm{F}}+\left \|B(\tau)\right \|_{\mathrm{F}}.	
		\end{equation}
	\end{lemma}
	Let 
	$X^{\ast}(\tau_{k+1})$ be the theoretical solution of $X(\tau_{k+1})F(\tau_{k+1})-A(\tau_{k+1})\overline{X}(\tau_{k+1})-C(\tau_{k+1})=O$. When $k$ is sufficiently large, 
	$X(\tau_{k+1})=X^{\ast}(\tau_{k+1})+\mathcal{O}(\varepsilon^{2})$. Furthermore, it is further obtained that:
	\begin{equation}
		\begin{split}
			&\lim_{k \to +\infty } \mathrm{sup} \left \|X(\tau_{k+1})F(\tau_{k+1})-A(\tau_{k+1})\overline{X}(\tau_{k+1})-C(\tau_{k+1})  \right \|_{\mathrm{F}}
			\\
			=&\lim_{k \to +\infty } \mathrm{sup} \left \|\left ( X^{\ast}(\tau_{k+1})+\mathcal{O}(\varepsilon^{2})\right )F(\tau_{k+1})-A(\tau_{k+1})\overline{\left (X^{\ast}(\tau_{k+1})+\mathcal{O}(\varepsilon^{2})\right )}-C(\tau_{k+1})  \right \|_{\mathrm{F}}
			\\
			=&\lim_{k \to +\infty } \mathrm{sup} \left \|\left ( X^{\ast}(\tau_{k+1})+\mathcal{O}(\varepsilon^{2})\right )F(\tau_{k+1})-A(\tau_{k+1})\left (\overline{X^{\ast}}(\tau_{k+1})+\overline{\mathcal{O}}(\varepsilon^{2})\right )-C(\tau_{k+1})  \right \|_{\mathrm{F}}
			\\
			=&\lim_{k \to +\infty } \mathrm{sup} \left \| \mathcal{O}(\varepsilon^{2})F(\tau_{k+1})-A(\tau_{k+1})\overline{\mathcal{O}}(\varepsilon^{2})  \right \|_{\mathrm{F}}.  
		\end{split}
	\end{equation}
	Next, according to \eqref{eq.frobenius.conj} and \eqref{eq.frobenius.triangle}, along with squeeze theorem, \eqref{eq.frobenius.triangle.linearerrconcznd1} can be obtained:
	\begin{equation}\label{eq.frobenius.triangle.linearerrconcznd1} 
		\begin{split}
			&\left \|\mathcal{O}(\varepsilon^{2})F(\tau_{k+1})\right \|_{\mathrm{F}}-\left \|A(\tau_{k+1})\overline{\mathcal{O}}(\varepsilon^{2})\right \|_{\mathrm{F}}
			\\
			\le&\left \|\mathcal{O}(\varepsilon^{2})F(\tau_{k+1})-A(\tau_{k+1})\overline{\mathcal{O}}(\varepsilon^{2})\right \|_{\mathrm{F}}
			\\
			\le&\left \|\mathcal{O}(\varepsilon^{2})F(\tau_{k+1})\right \|_{\mathrm{F}}+\left \|A(\tau_{k+1})\overline{\mathcal{O}}(\varepsilon^{2})\right \|_{\mathrm{F}}.
		\end{split}	
	\end{equation}
	For $\left \|\mathcal{O}(\varepsilon^{2})F(\tau_{k+1})\right \|_{\mathrm{F}}-\left \|A(\tau_{k+1})\overline{\mathcal{O}}(\varepsilon^{2})\right \|_{\mathrm{F}}$,
	\begin{equation}\label{eq.frobenius.conj1.linearerrconcznd1} 
	\begin{split}
		&\lim_{k \to +\infty } \mathrm{sup} \left \|\mathcal{O}(\varepsilon^{2})F(\tau_{k+1})\right \|_{\mathrm{F}}-\left \|A(\tau_{k+1})\overline{\mathcal{O}}(\varepsilon^{2})\right \|_{\mathrm{F}}
		\\
		=&\lim_{k \to +\infty } \mathrm{sup} \left \|\mathcal{O}(\varepsilon^{2})F(\tau_{k+1})\right \|_{\mathrm{F}}-\left \|\overline{A}(\tau_{k+1})\mathcal{O}(\varepsilon^{2})\right \|_{\mathrm{F}}
		\\
		=&\mathcal{O}(\varepsilon^{2}).		
	\end{split}	
\end{equation}
	For $\left \|\mathcal{O}(\varepsilon^{2})F(\tau_{k+1})\right \|_{\mathrm{F}}+\left \|A(\tau_{k+1})\overline{\mathcal{O}}(\varepsilon^{2})\right \|_{\mathrm{F}}$,
	\begin{equation}\label{eq.frobenius.conj2.linearerrconcznd1} 
		\begin{split}
			&\lim_{k \to +\infty } \mathrm{sup} \left \|\mathcal{O}(\varepsilon^{2})F(\tau_{k+1})\right \|_{\mathrm{F}}+\left \|A(\tau_{k+1})\overline{\mathcal{O}}(\varepsilon^{2})\right \|_{\mathrm{F}}
			\\
			=&\lim_{k \to +\infty } \mathrm{sup} \left \|\mathcal{O}(\varepsilon^{2})F(\tau_{k+1})\right \|_{\mathrm{F}}+\left \|\overline{A}(\tau_{k+1})\mathcal{O}(\varepsilon^{2})\right \|_{\mathrm{F}}
			\\
			=&\mathcal{O}(\varepsilon^{2}).		
		\end{split}	
	\end{equation}
	Then, combined with squeeze theorem, \eqref{eq.frobenius.conj.triangle.linearerrconcznd1} can be obtained:
	\begin{equation}\label{eq.frobenius.conj.triangle.linearerrconcznd1}
		\lim_{k \to +\infty } \mathrm{sup} \left \| \mathcal{O}(\varepsilon^{2})F(\tau_{k+1})-A(\tau_{k+1})\overline{\mathcal{O}}(\varepsilon^{2})  \right \|_{\mathrm{F}}=\mathcal{O}(\varepsilon^{2}).
	\end{equation}
	The proof is thus completed.	
	\end{proof} 
\end{theorem}

\section{Numerical experiments and validation}
Subsequent experiments are made in MATLAB but throw out ode45
\cite{He2024ZeroingND,He2024RevisitingTC,wuGaussNewtonMethod2024,Kong2024ReservoircomputingBA} solver, as ode45 is a variable-step size solver. Inspired by the mathematical concept of ``Measure theory", to establish a concept of step size, $\varepsilon$ of 0.1 and 0.001 are adopted, with a total simulation time of $\tau$ equals 10. It means
\begin{equation}\label{eq.stepsize}
	k=\frac{10}{\varepsilon},
\end{equation}
where $\varepsilon$ represents the step size between points, which is constant-step sizes. $k$ is the number of points generated for each step size. For $\varepsilon$ of 0.1, $k$ equals 100, and for $\varepsilon$ of 0.001, $k$ equals 10000. At the same time, the data precision is set using the ``format long" command in Matlab.

In this article, to highlight the differences between sampling discretion errors and space compressive approximation errors in neural dynamics, two typical examples are selected.
\begin{example}\label{example1}
Consider the following SSCME
\cite{wuSolutionsMatrixEquations2006}.
\begin{equation}\label{eq.example1}
		\begin{split} 
			F
			=\begin{bmatrix}
				0	&0 \\
				1	&-1
			\end{bmatrix}
			&+\mathrm{i}\begin{bmatrix}
				2	& 1 \\
				0 & 1
			\end{bmatrix}\in
			\mathbb{C}^{2\times 2},		
			\\
			A
			=\begin{bmatrix}
				1& -2& -1 \\
				0& 0& 0 \\
				0& -1& 1 \\
			\end{bmatrix}
			&+\mathrm{i}\begin{bmatrix}
		0& -1& 1 \\
			0& 1& 0 \\
			0& 0& -1 \\
			\end{bmatrix}\in
			\mathbb{C}^{3\times 3},		
			\\
			C
			=\begin{bmatrix}
				-1& 1\\
				0& 0\\
				0& 1\\
			\end{bmatrix}
			&+\mathrm{i}\begin{bmatrix}
				1&0 \\
				0&1 \\
				-1&-2 \\
			\end{bmatrix}\in
			\mathbb{C}^{3\times 2},			
		\end{split}
\end{equation}
The only theoretical solution to this example is
 \begin{equation}\label{eq.example1.solve}
 	\renewcommand{\arraystretch}{2}
 			X^*(\tau)=\begin{bmatrix}
 					\dfrac{21}{40} & -\dfrac{9}{8}\\
 					\dfrac{1}{2} & \dfrac{3}{4}\\
 					-\dfrac{6}{5} & -\dfrac{13}{20}\\
 				\end{bmatrix}+\mathrm{i}\begin{bmatrix}
 					-\dfrac{33}{40} & \dfrac{5}{8} \\
 					\dfrac{1}{4} & -\dfrac{1}{2} \\
 					\dfrac{21}{20} & \dfrac{9}{5} \\
 				\end{bmatrix}\in
 			\mathbb{C}^{3\times 2}.
 \end{equation}
\end{example}
\begin{example}\label{example2}
	Consider the following TVSSCME
	\cite{He2024ZeroingND,He2024RevisitingTC},
	where 
	$s(\tau)$ and 
	$c(\tau)$ represent 
	$\mathrm{sin}(\tau)$ and 
	$\mathrm{cos}(\tau)$, respectively.
	\begin{equation}\label{eq.example2}
	\begin{split} 
		F(\tau)
		=\begin{bmatrix}
			6+s(\tau)	&c(\tau) \\
			c(\tau)	&4+s(\tau)
		\end{bmatrix}
		&+\mathrm{i}\begin{bmatrix}
			c(\tau)	& s(\tau) \\
			s(\tau) & c(\tau)
		\end{bmatrix}\in
		\mathbb{C}^{2\times 2},		
	\\
		A(\tau)
		=\begin{bmatrix}
			c(\tau)	&s(\tau) \\
			-s(\tau)	&c(\tau)
		\end{bmatrix}
		&+\mathrm{i}\begin{bmatrix}
			s(\tau)	& c(\tau) \\
			c(\tau) & -s(\tau)
		\end{bmatrix}\in
		\mathbb{C}^{2\times 2},		
\\
		C(\tau)
		=\begin{bmatrix}
			c_{\mathrm{r},11}(\tau) & c_{\mathrm{r},12}(\tau) \\
			c_{\mathrm{r},21}(\tau) & c_{\mathrm{r},22}(\tau) \\
		\end{bmatrix}
		&+\mathrm{i}\begin{bmatrix}
			c_{\mathrm{i},11}(\tau) & 	c_{\mathrm{i},12}(\tau) \\
			c_{\mathrm{i},21}(\tau) & 	c_{\mathrm{i},22}(\tau) \\
		\end{bmatrix}\in
		\mathbb{C}^{2\times 2},			
	\end{split}
	\end{equation}
		where
		$c_{\mathrm{r},11}(\tau)=2c^{2}(\tau)-2c(\tau)s(\tau)+6s(\tau)$,
		$c_{\mathrm{r},12}(\tau)=4c(\tau)+2c(\tau)s(\tau)-2c^{2}(\tau)$, $c_{\mathrm{r},21}(\tau)=-2s(2\tau)-6c(\tau)+2$, 
		$c_{\mathrm{r},22}(\tau)=2s(2\tau)-4s(\tau)-2$ 
		and 
		$c_{\mathrm{i},11}(\tau)=2c^{2}(\tau)+2c(\tau)s(\tau)+6s(\tau)$, $c_{\mathrm{i},12}(\tau)=4c(\tau)+2c(\tau)s(\tau)+2c^{2}(\tau)$, 
		$c_{\mathrm{i},21}(\tau)=-2s(2\tau)-6c(\tau)-2$, $c_{\mathrm{i},22}(\tau)=-2s(2\tau)-4s(\tau)-2$.
		
		The only theoretical solution to this example is
		\begin{equation}\label{eq.example2.solve}
				X^*(\tau)=\begin{bmatrix}
						s(\tau) & c(\tau)\\
						-c(\tau) & -s(\tau)\\
					\end{bmatrix}+\mathrm{i}\begin{bmatrix}
						s(\tau) & c(\tau) \\
						-c(\tau) & -s(\tau) \\
					\end{bmatrix}\in
				\mathbb{C}^{2\times 2}.
		\end{equation}	
\end{example}
\begin{remark}
Unless otherwise specified, the following numerical experiments use random initial values in the interval $\left [-5,5  \right ] $. Like y$\left ( \cdot \right )$ represents the corresponding theoretical solution elements of x$\left ( \cdot \right )$, $\left ( \cdot \right )$ means element position and associated matrix. Then $x_{\mathrm{r},st}(\tau)$'s theoretical solution is $y_{\mathrm{r},st}(\tau)$. The red line represents Con-DZND1-2i, while the green line represents Con-DZND2-2i. Because of Fig. \ref{fig.dif.sscme.tvsscme}, the residuals are uniformly defined $\left \|X(\tau)-X^*(\tau)   \right \|_{\mathrm{F}}$, where $\left \|\cdot   \right \|_{\mathrm{F}}$ stands for Frobenius norm.
\end{remark}
\subsection{For Con-DZND1-2i model and Con-DZND2-2i model using step sizes $\varepsilon$ of 0.1 and 0.001 with $\gamma$ equals 10}
\begin{figure}[!h]\centering
	\subfigure[]{\includegraphics[width=0.32\columnwidth]{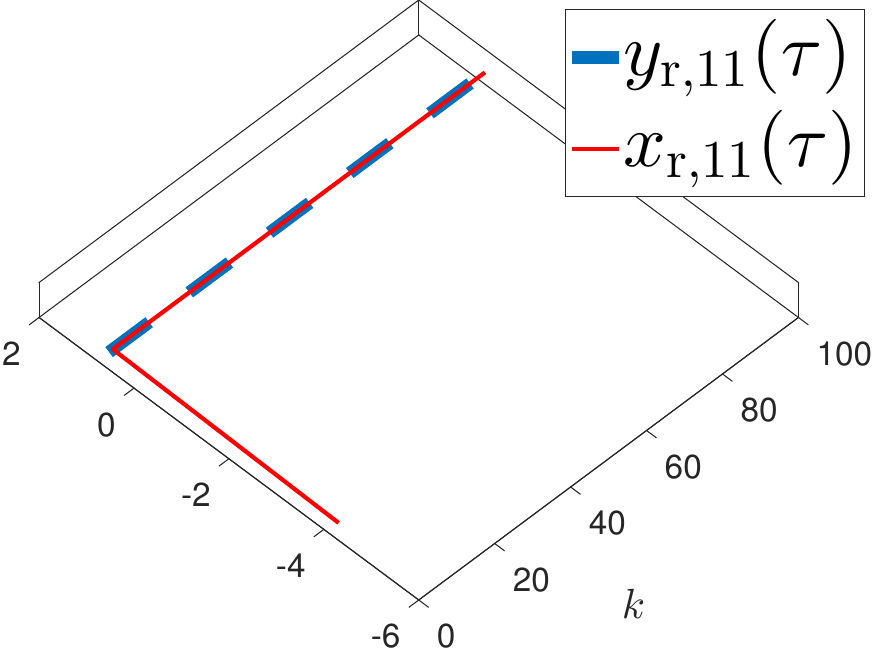}}
	\subfigure[]{\includegraphics[width=0.32\columnwidth]{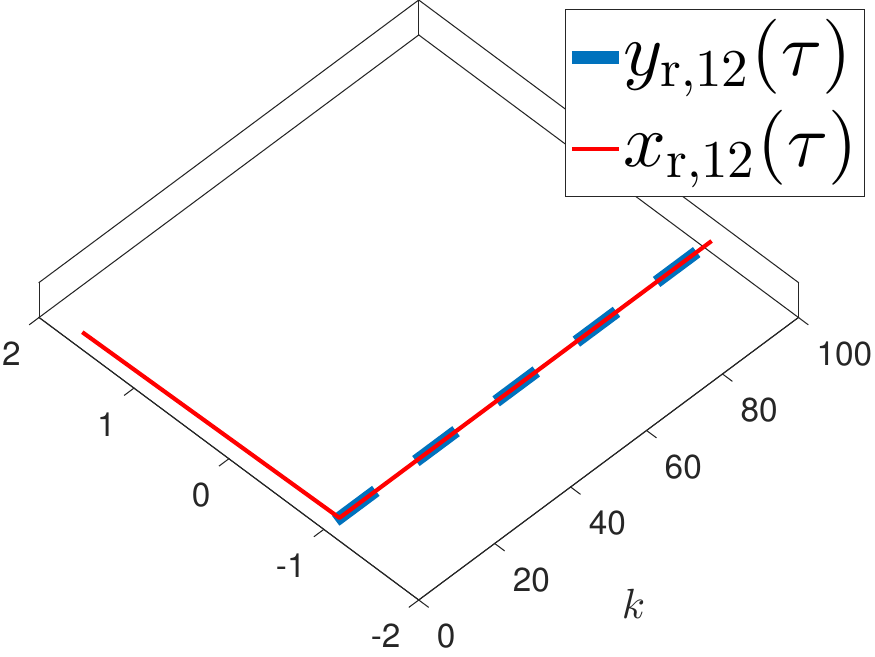}}
	\subfigure[]{\includegraphics[width=0.32\columnwidth]{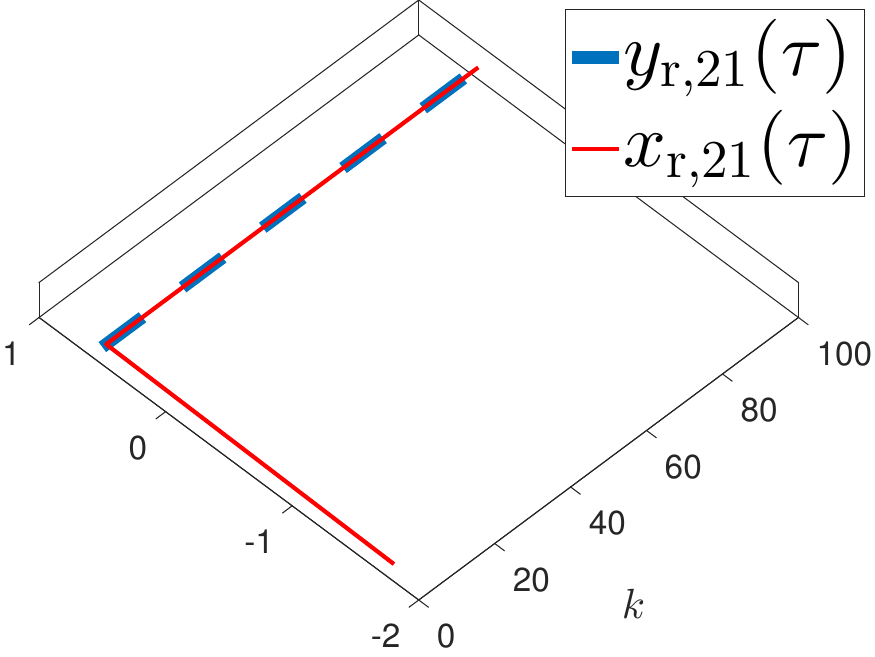}}
	\subfigure[]{\includegraphics[width=0.32\columnwidth]{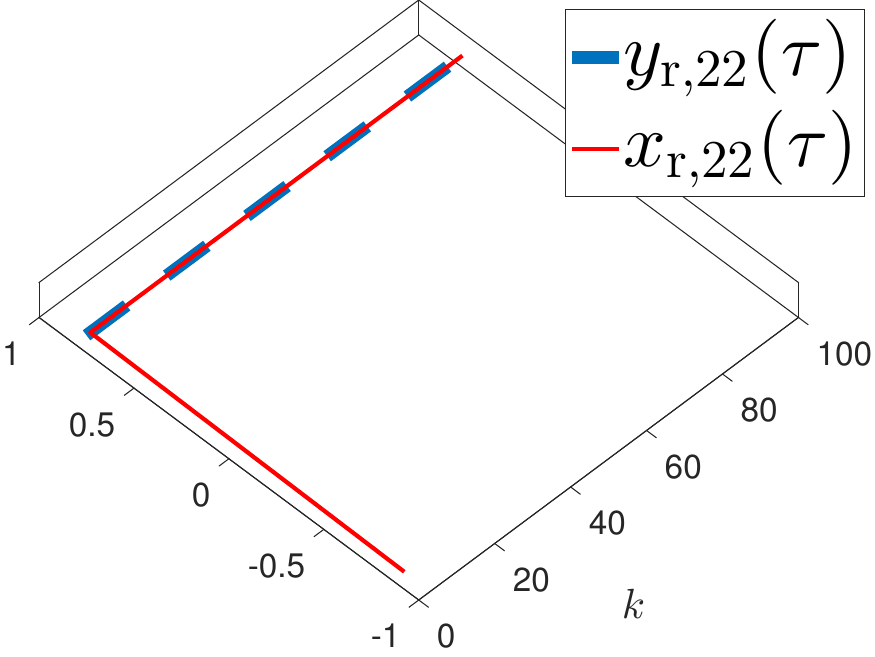}}
	\subfigure[]{\includegraphics[width=0.32\columnwidth]{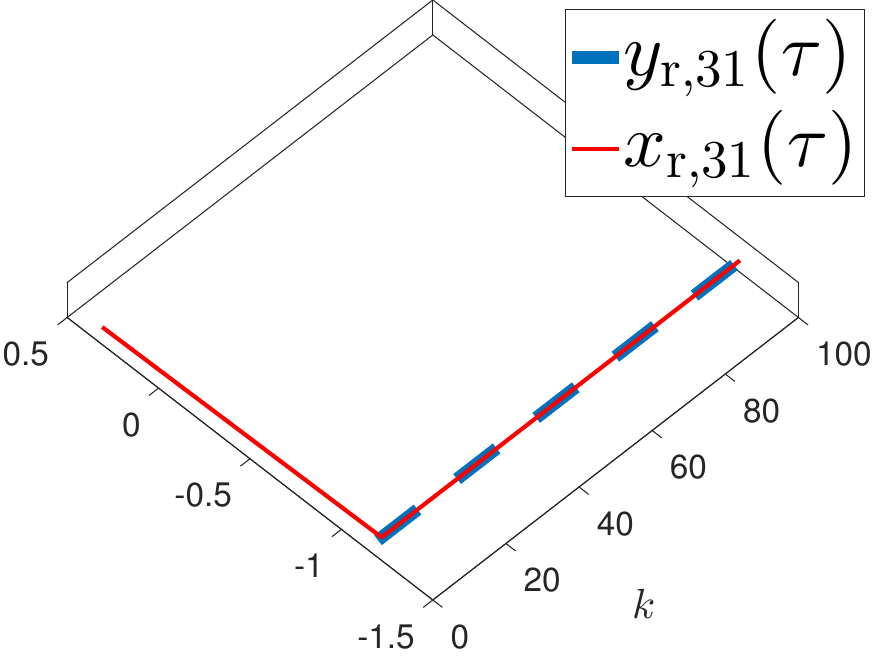}}
	\subfigure[]{\includegraphics[width=0.32\columnwidth]{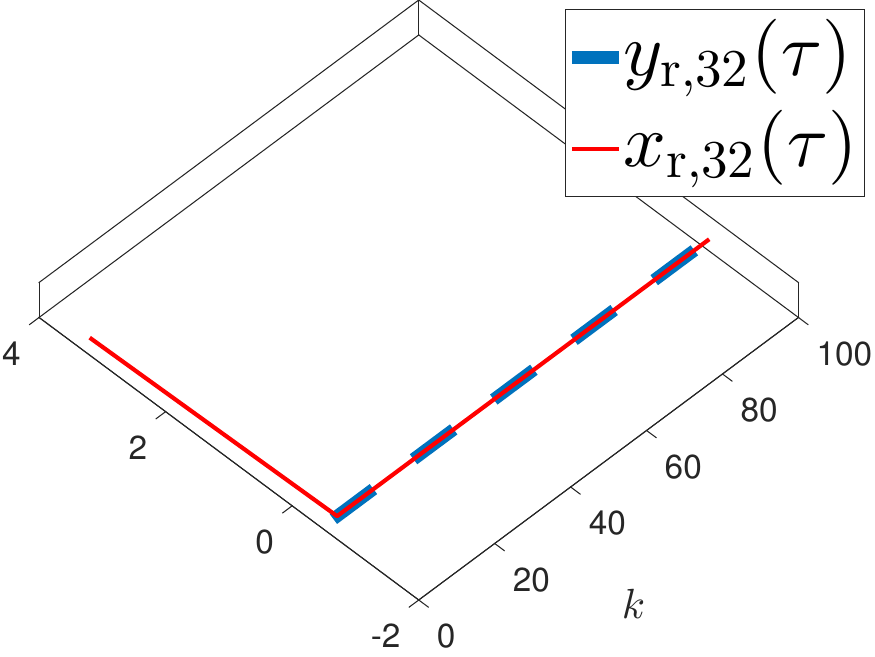}}
	\subfigure[]{\includegraphics[width=0.32\columnwidth]{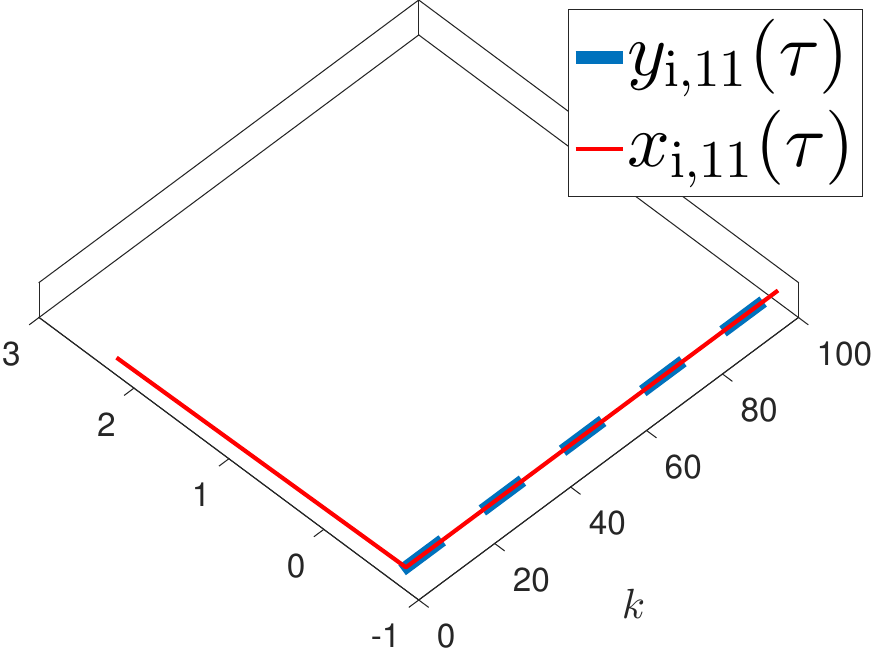}}
	\subfigure[]{\includegraphics[width=0.32\columnwidth]{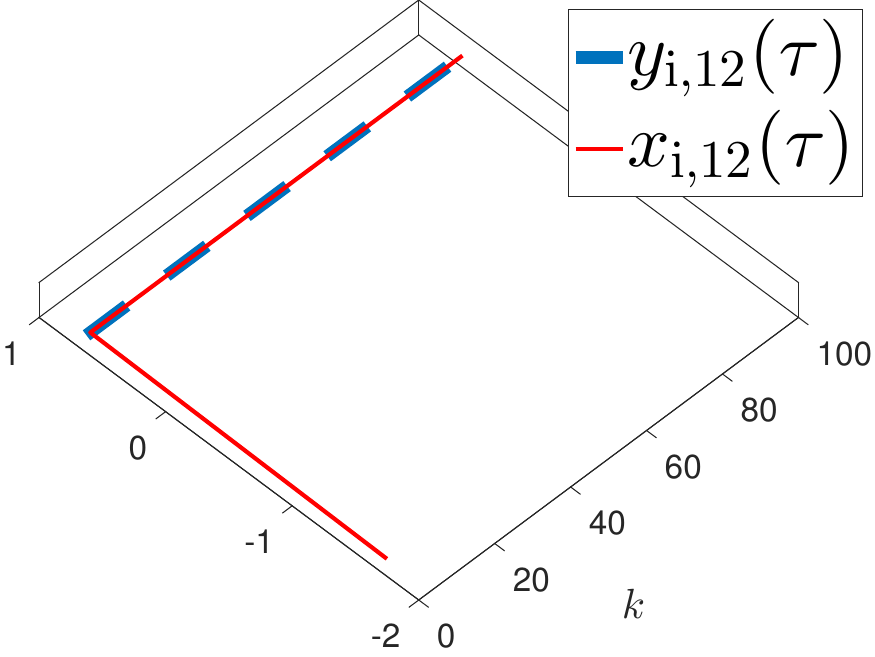}}
	\subfigure[]{\includegraphics[width=0.32\columnwidth]{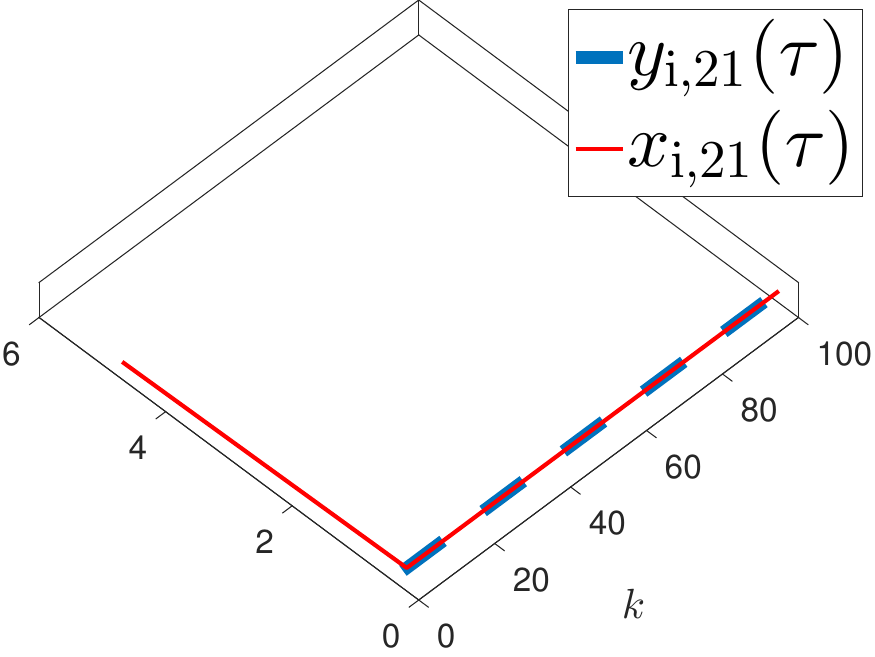}}
	\subfigure[]{\includegraphics[width=0.32\columnwidth]{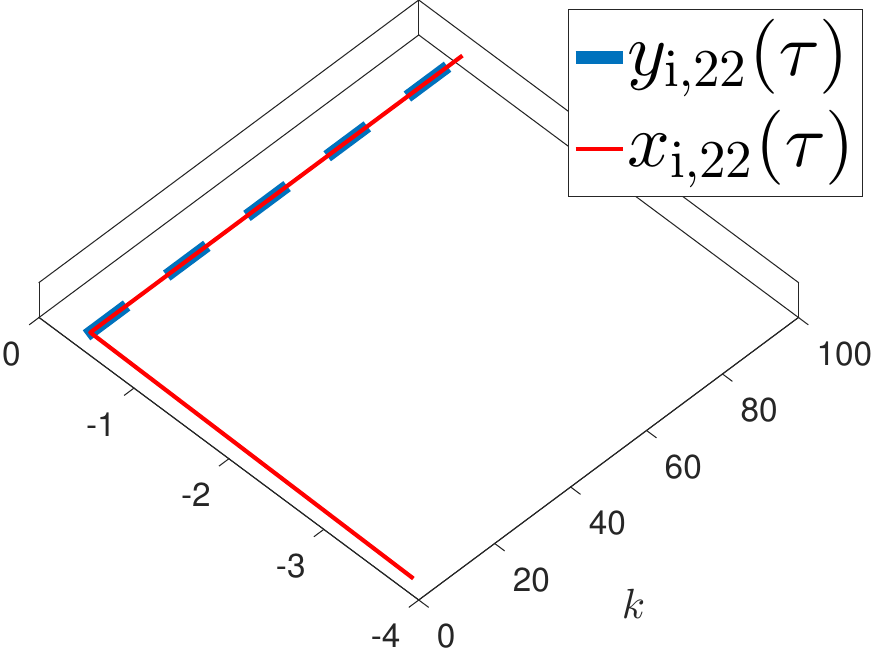}}
	\subfigure[]{\includegraphics[width=0.32\columnwidth]{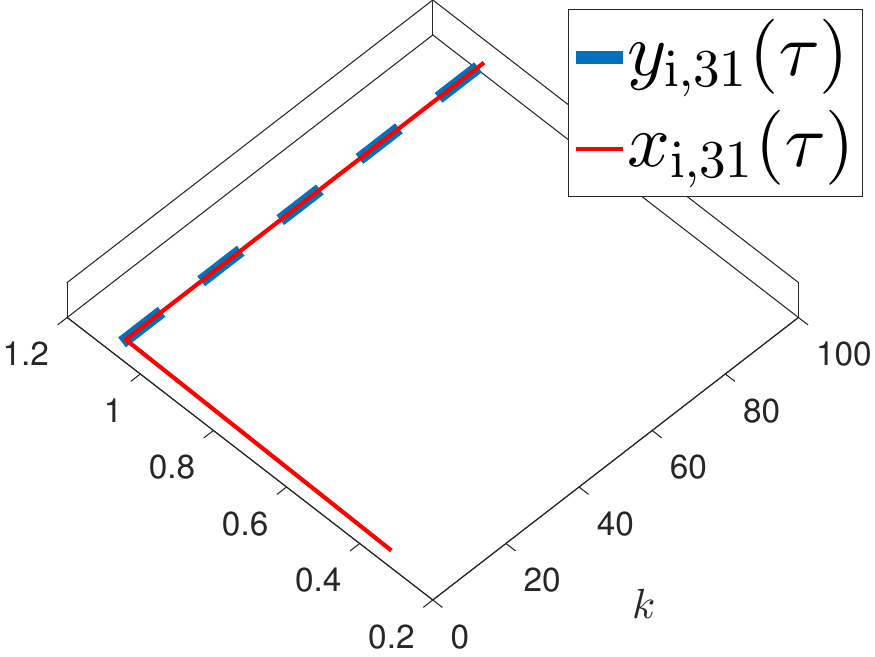}}
	\subfigure[]{\includegraphics[width=0.32\columnwidth]{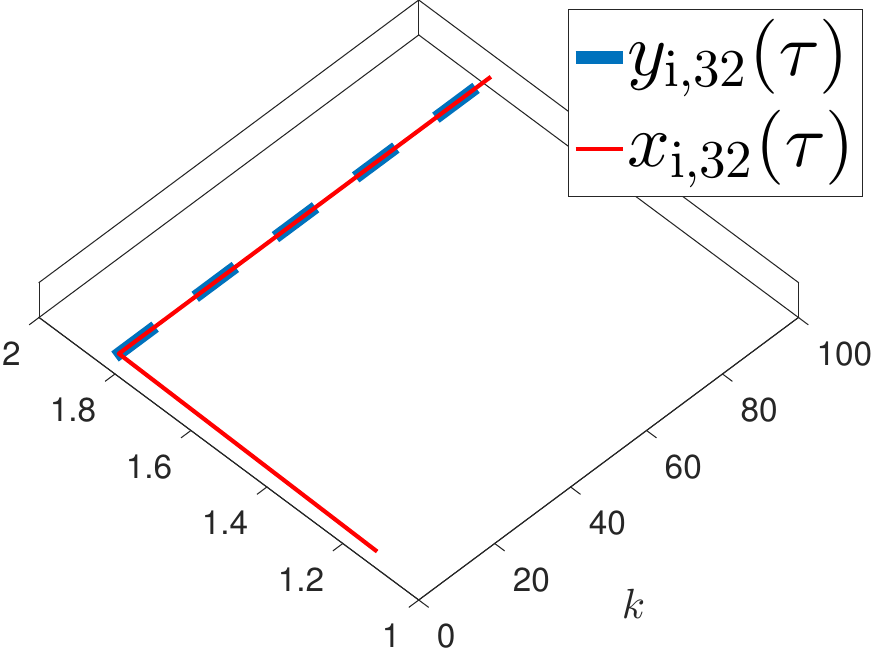}}
	\caption{Solution $X(\tau)$ computed by Con-DZND1-2i \eqref{eq.euler.forward.solve.linearerrconcznd1} model in Example \ref{example1} where $\gamma$ equals 10 and $\varepsilon$ equals 0.1.}
	\label{fig.e1.Con-DZND1-2i.solve.10.0.1}
\end{figure}

\begin{figure}[!h]\centering
	\subfigure[]{\includegraphics[width=0.32\columnwidth]{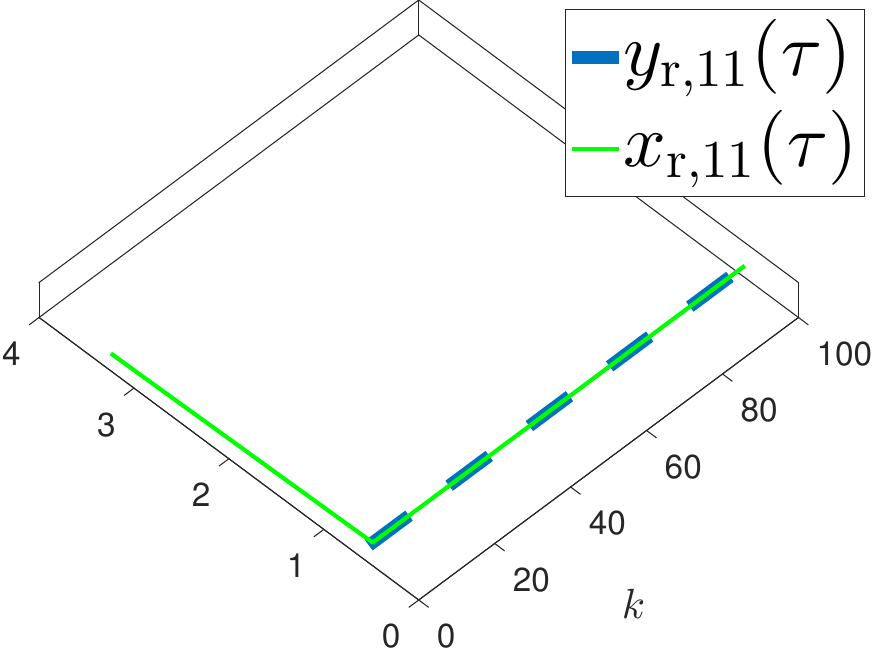}}
	\subfigure[]{\includegraphics[width=0.32\columnwidth]{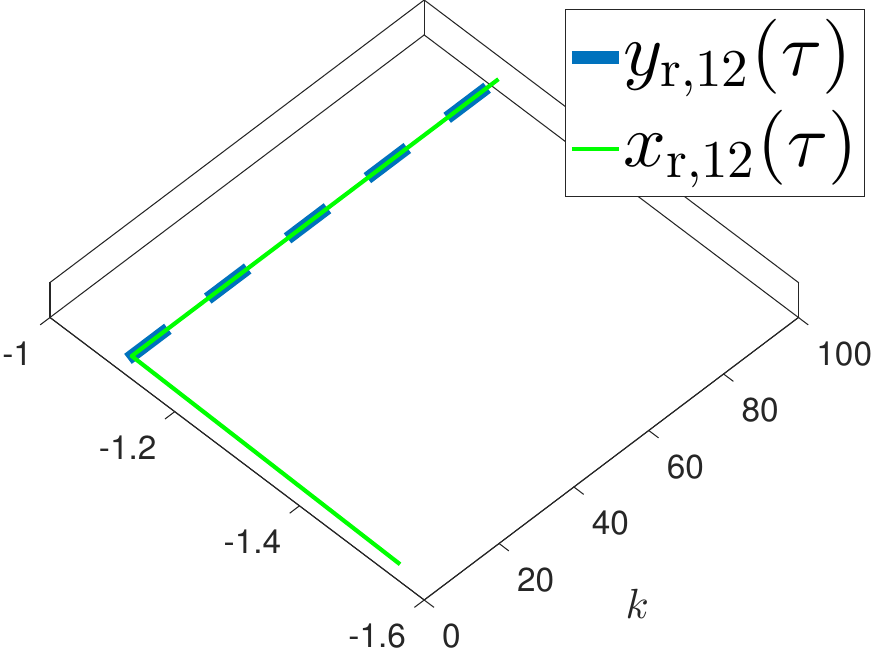}}
	\subfigure[]{\includegraphics[width=0.32\columnwidth]{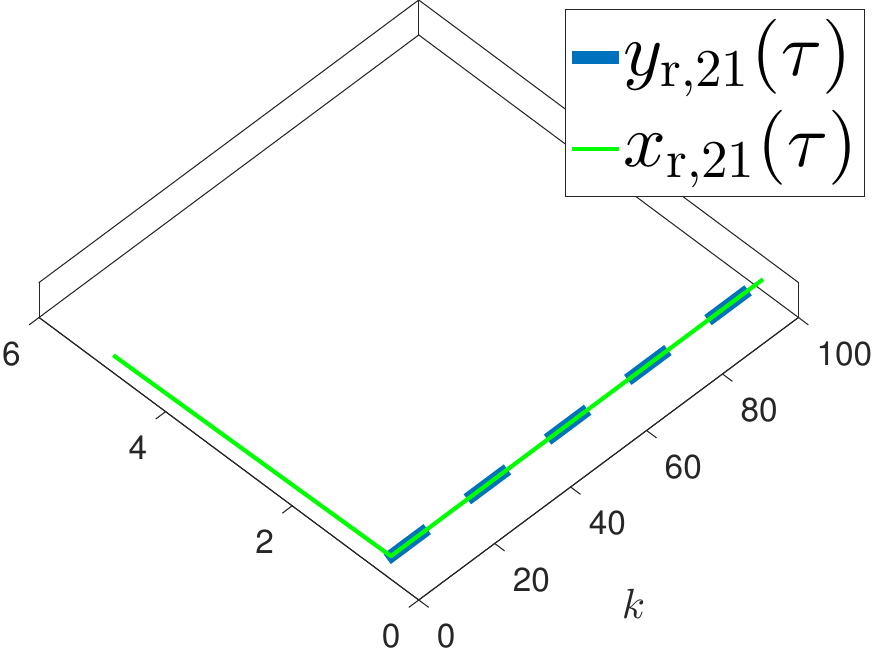}}
	\subfigure[]{\includegraphics[width=0.32\columnwidth]{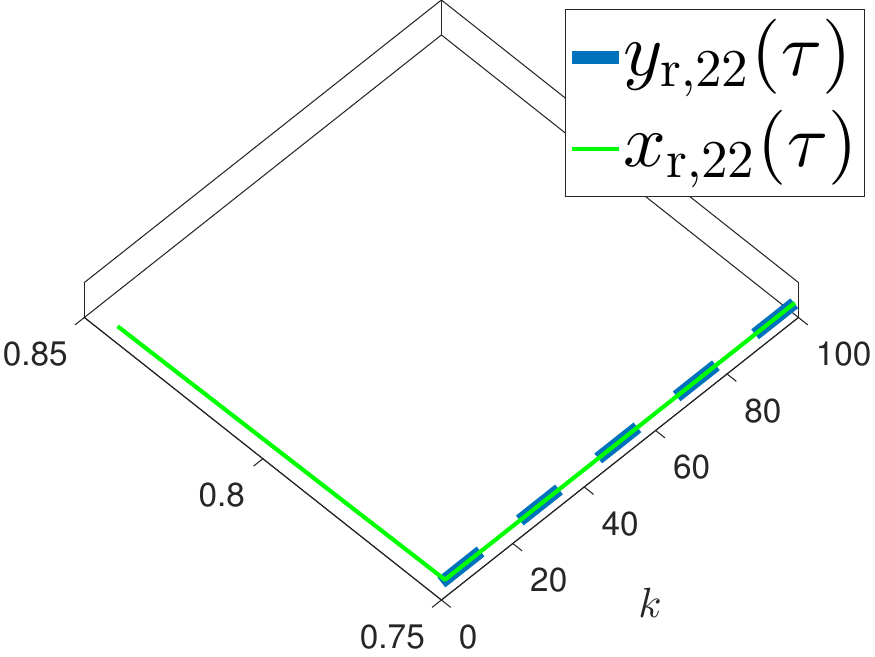}}
	\subfigure[]{\includegraphics[width=0.32\columnwidth]{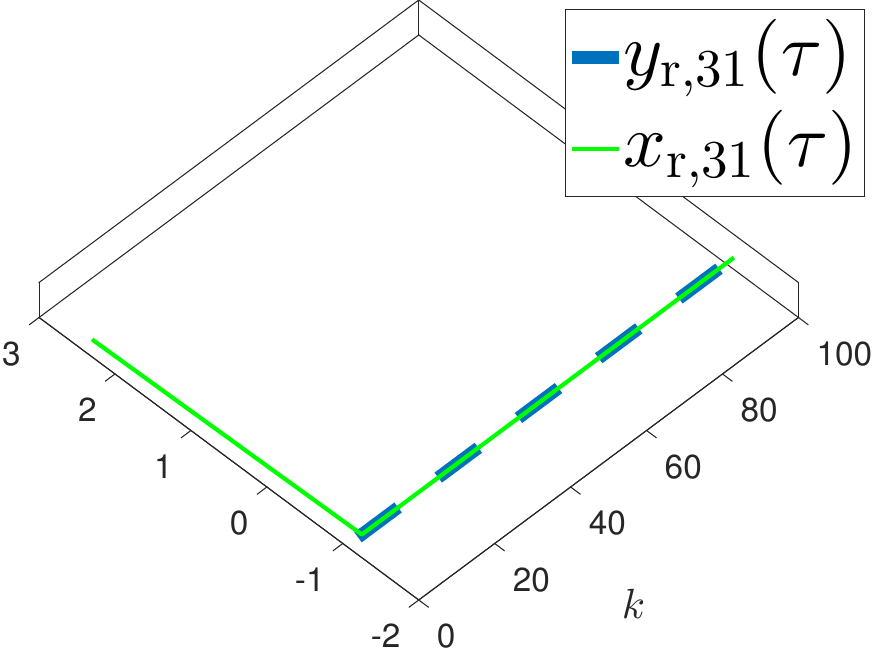}}
	\subfigure[]{\includegraphics[width=0.32\columnwidth]{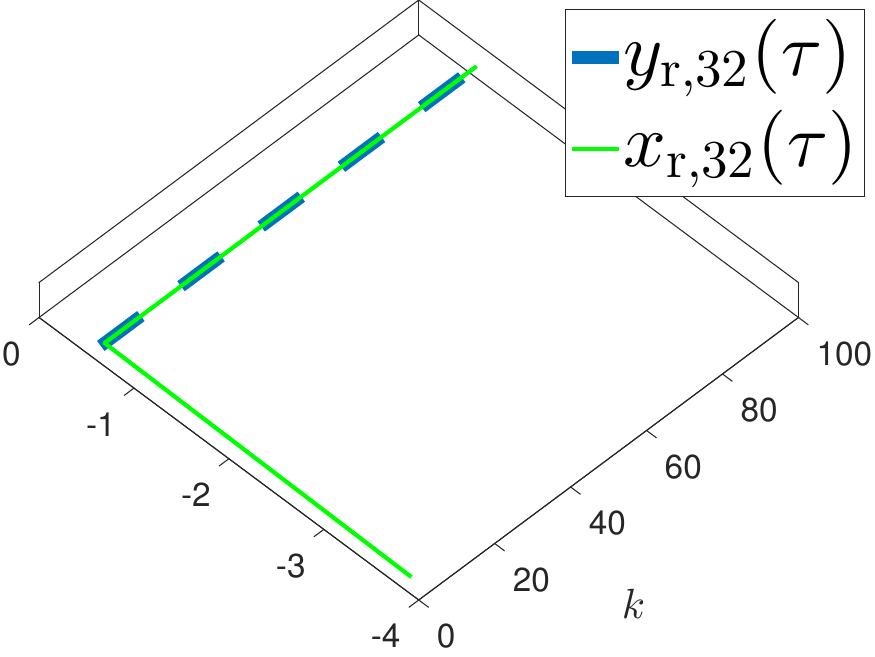}}
	\subfigure[]{\includegraphics[width=0.32\columnwidth]{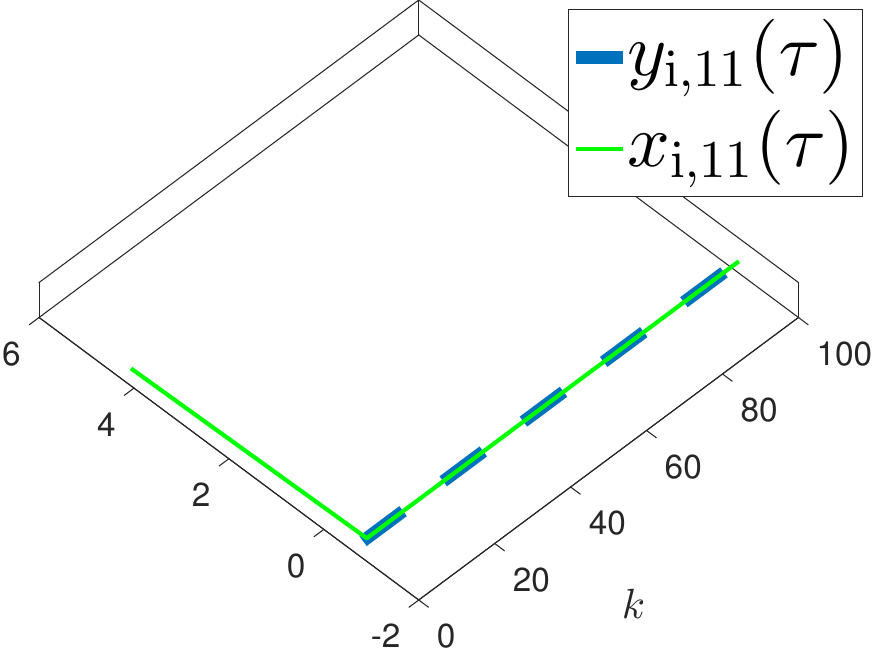}}
	\subfigure[]{\includegraphics[width=0.32\columnwidth]{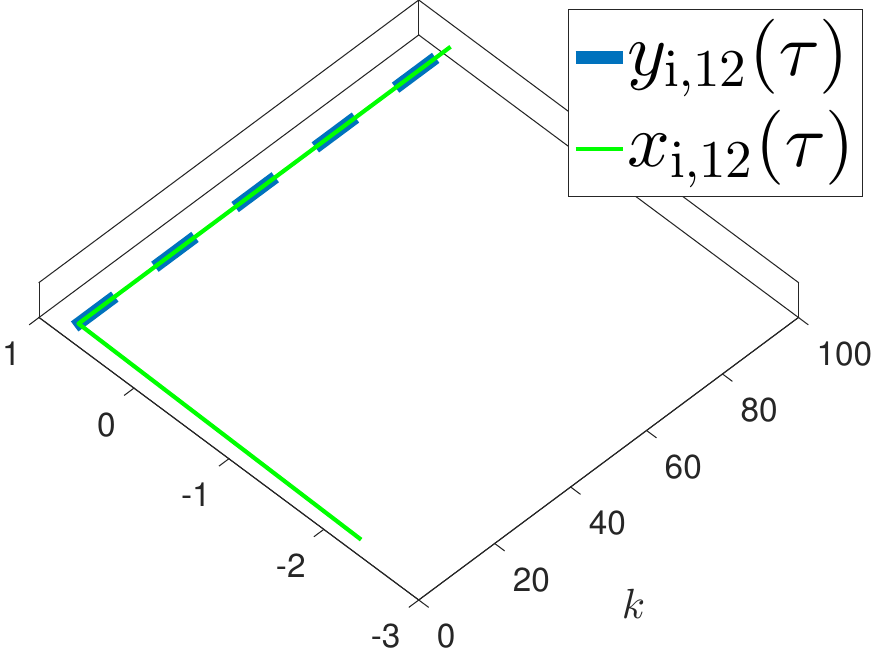}}
	\subfigure[]{\includegraphics[width=0.32\columnwidth]{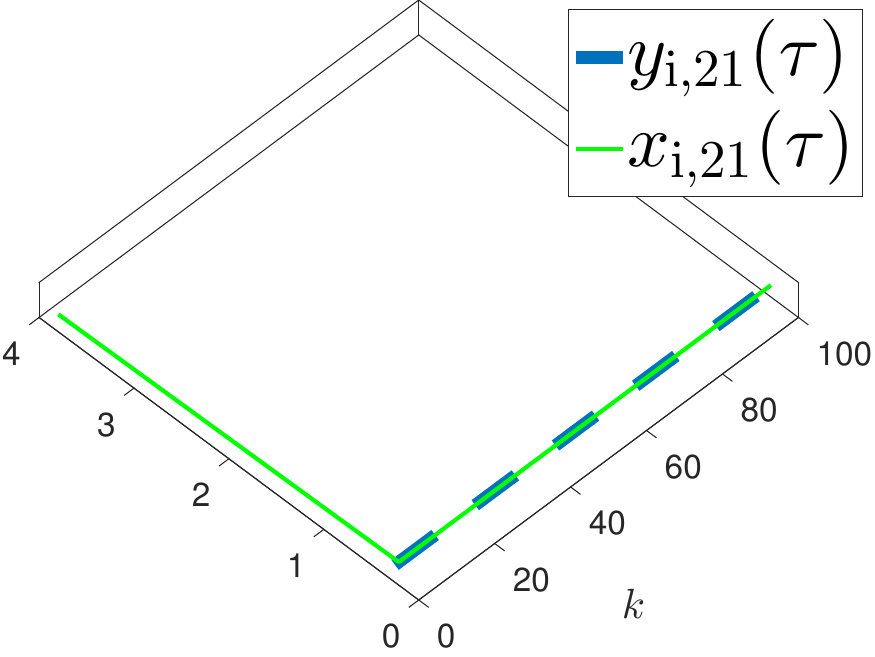}}
	\subfigure[]{\includegraphics[width=0.32\columnwidth]{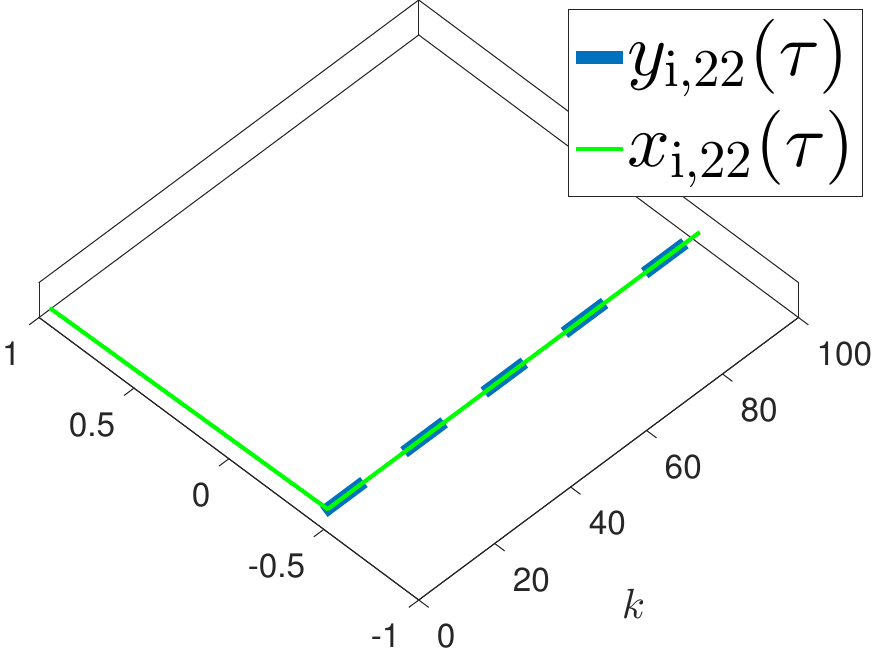}}
	\subfigure[]{\includegraphics[width=0.32\columnwidth]{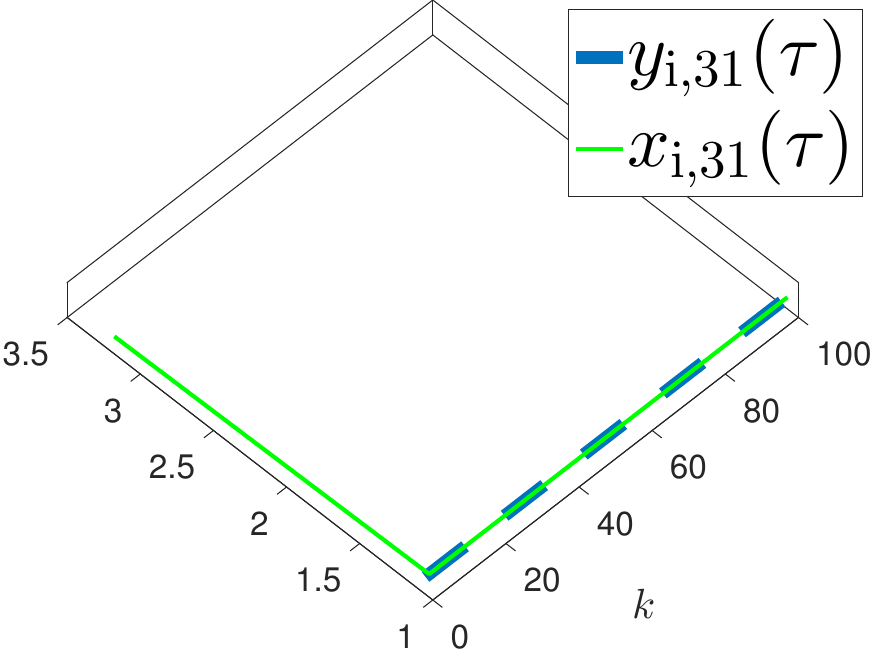}}
	\subfigure[]{\includegraphics[width=0.32\columnwidth]{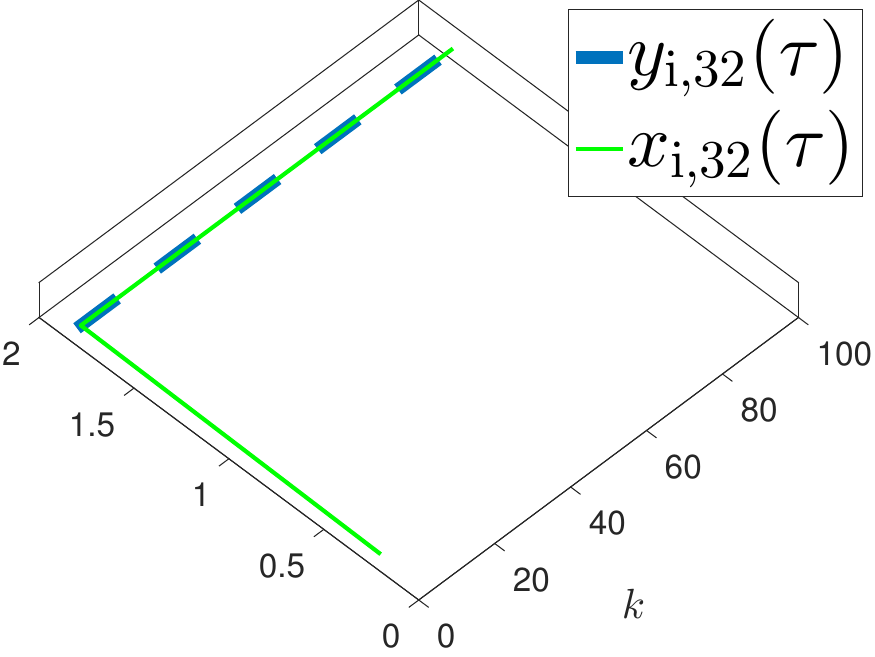}}
	\caption{Solution $X(\tau)$ computed by Con-DZND2-2i \eqref{eq.euler.forward.solve.linearerrconcznd2} model in Example \ref{example1} where $\gamma$ equals 10 and $\varepsilon$ equals 0.1.}
	\label{fig.e1.Con-DZND2-2i.solve.10.0.1}
\end{figure}

\begin{figure}[!h]\centering
	\subfigure[]{\includegraphics[width=0.32\columnwidth]{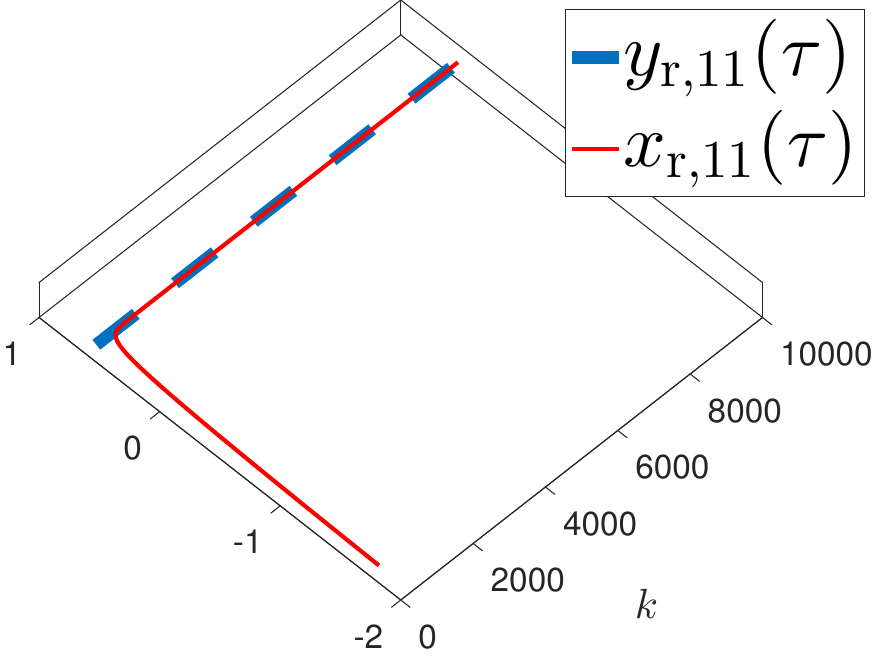}}
	\subfigure[]{\includegraphics[width=0.32\columnwidth]{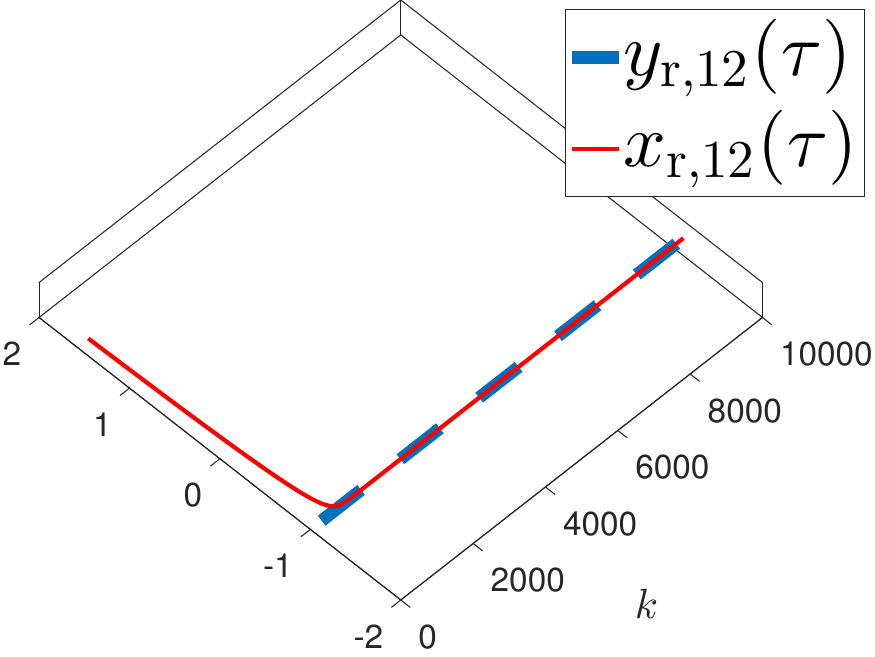}}
	\subfigure[]{\includegraphics[width=0.32\columnwidth]{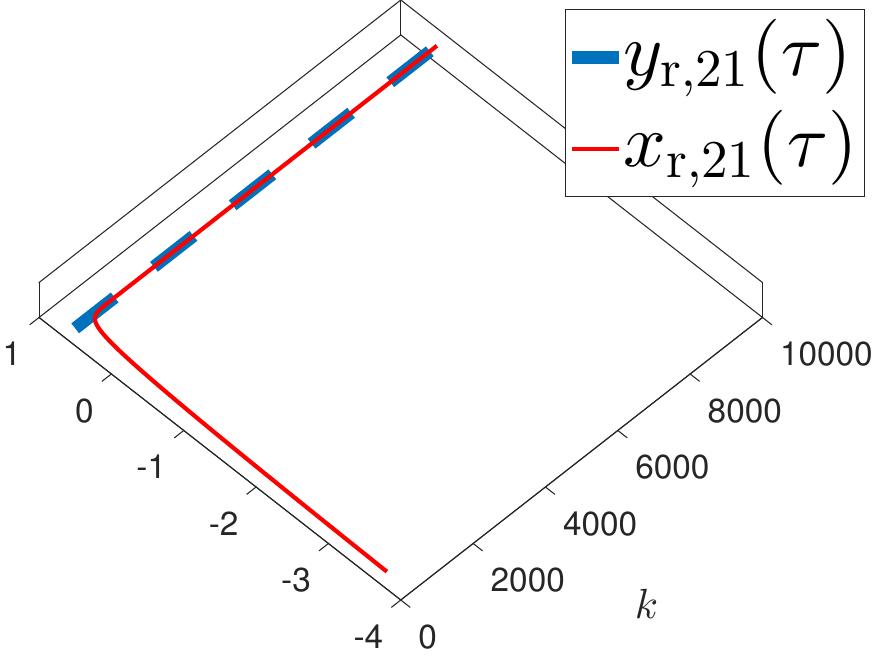}}
	\subfigure[]{\includegraphics[width=0.32\columnwidth]{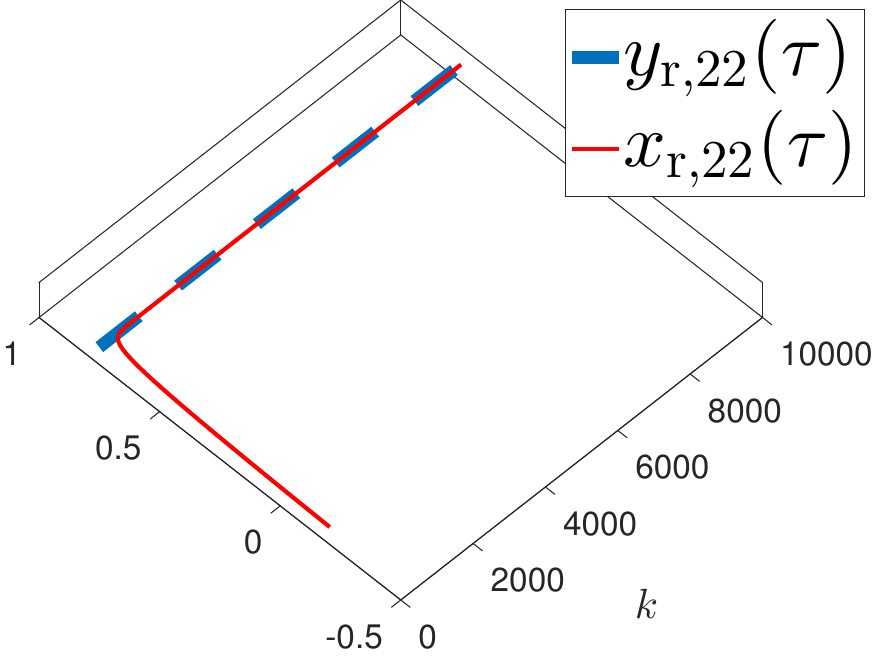}}
	\subfigure[]{\includegraphics[width=0.32\columnwidth]{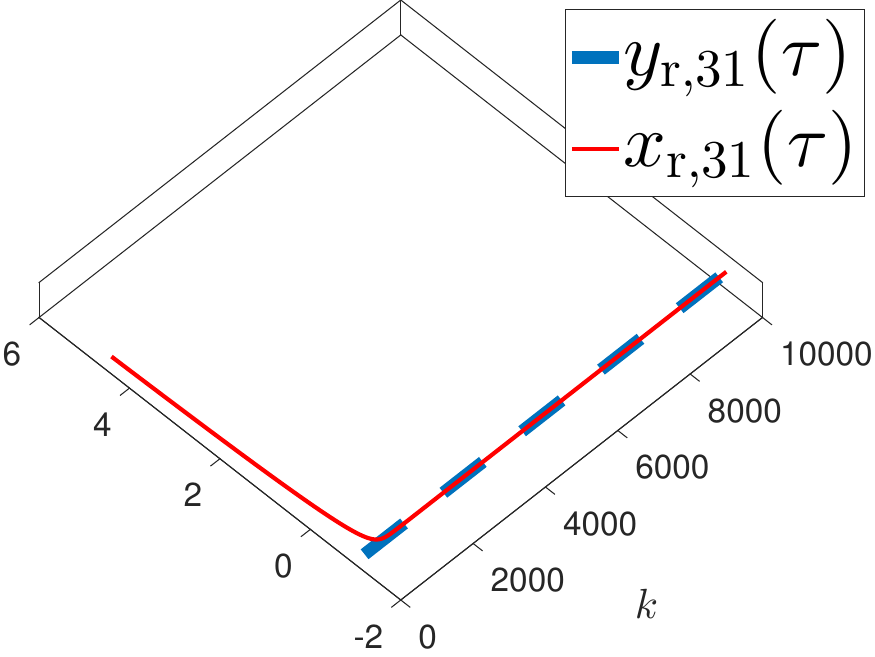}}
	\subfigure[]{\includegraphics[width=0.32\columnwidth]{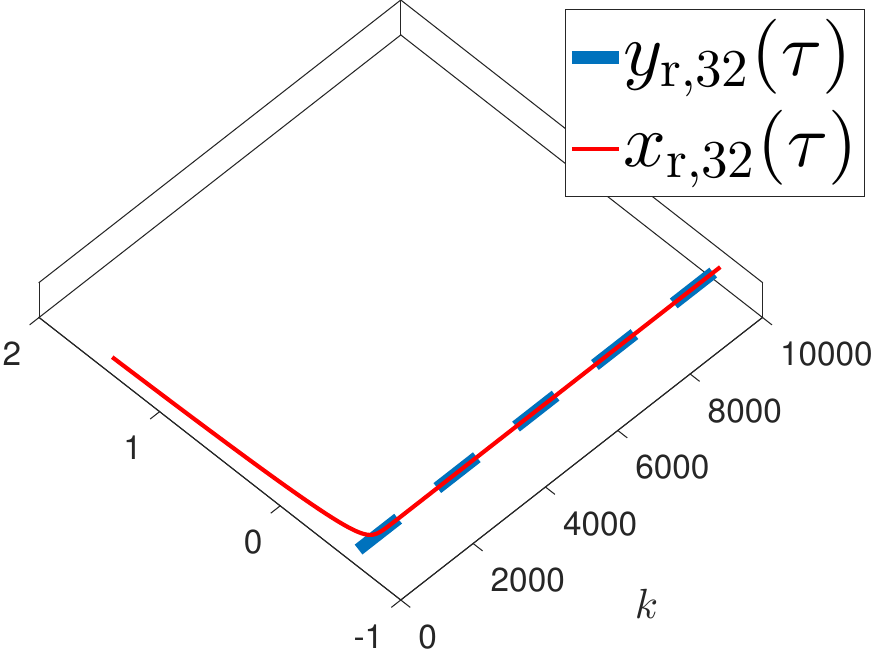}}
	\subfigure[]{\includegraphics[width=0.32\columnwidth]{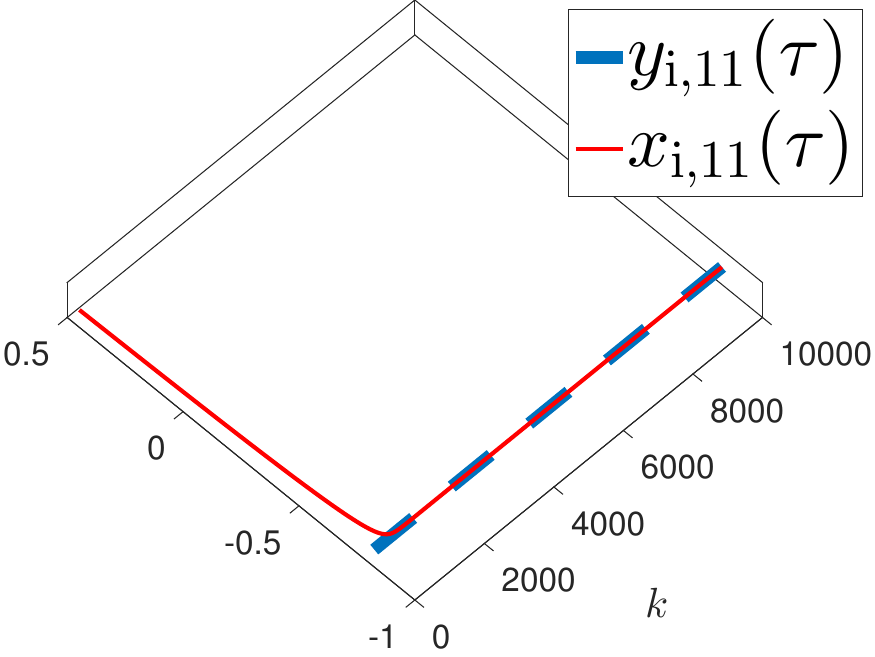}}
	\subfigure[]{\includegraphics[width=0.32\columnwidth]{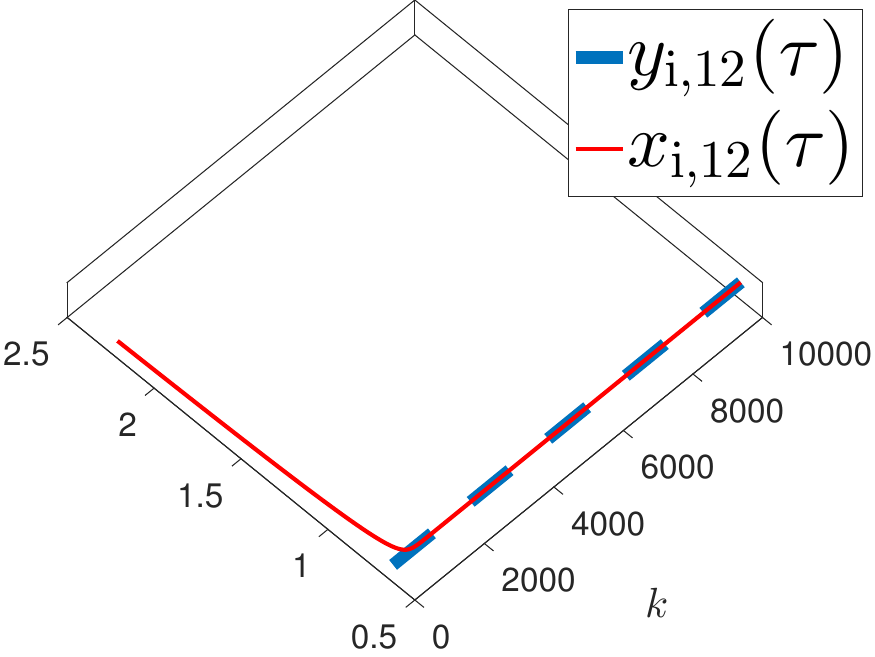}}
	\subfigure[]{\includegraphics[width=0.32\columnwidth]{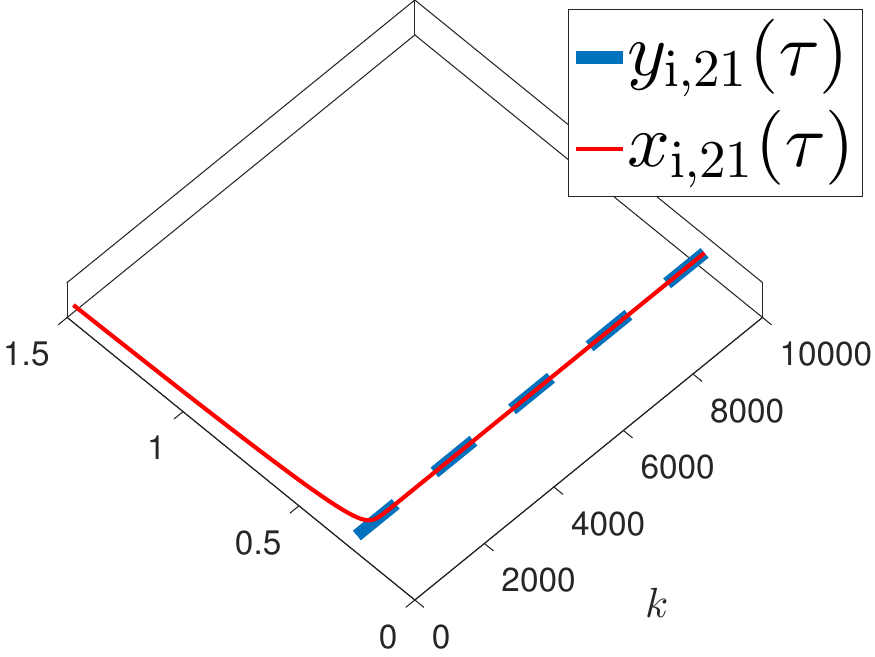}}
	\subfigure[]{\includegraphics[width=0.32\columnwidth]{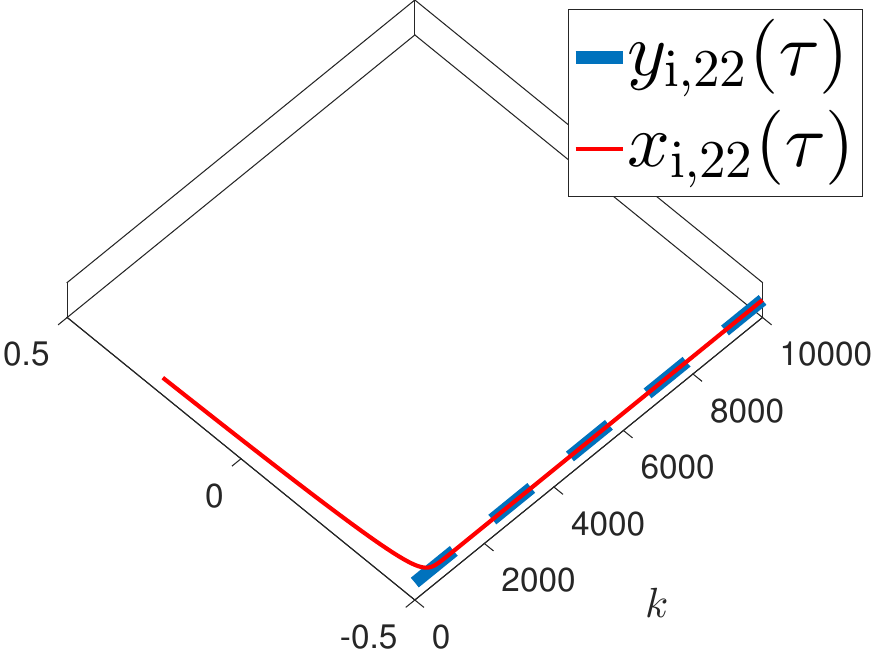}}
	\subfigure[]{\includegraphics[width=0.32\columnwidth]{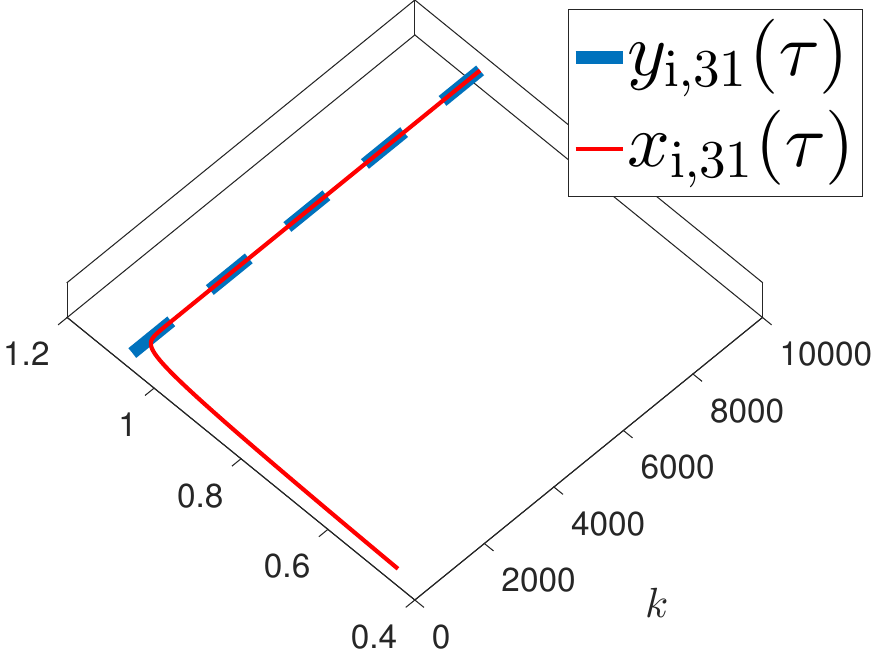}}
	\subfigure[]{\includegraphics[width=0.32\columnwidth]{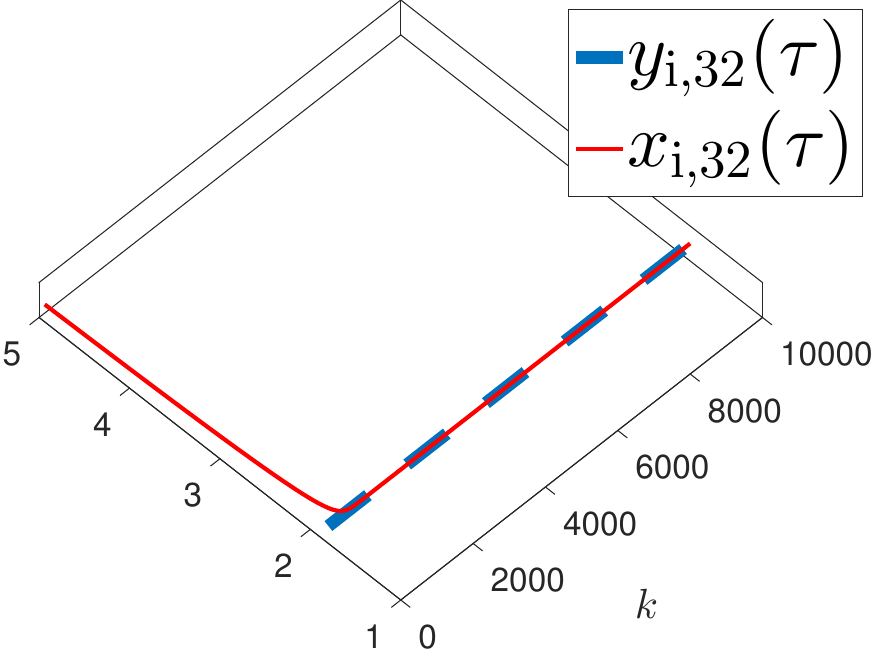}}
	\caption{Solution $X(\tau)$ computed by Con-DZND1-2i \eqref{eq.euler.forward.solve.linearerrconcznd1} model in Example \ref{example1} where $\gamma$ equals 10 and $\varepsilon$ equals 0.001.}
	\label{fig.e1.Con-DZND1-2i.solve.10.0.001}
\end{figure}

\begin{figure}[!h]\centering
	\subfigure[]{\includegraphics[width=0.32\columnwidth]{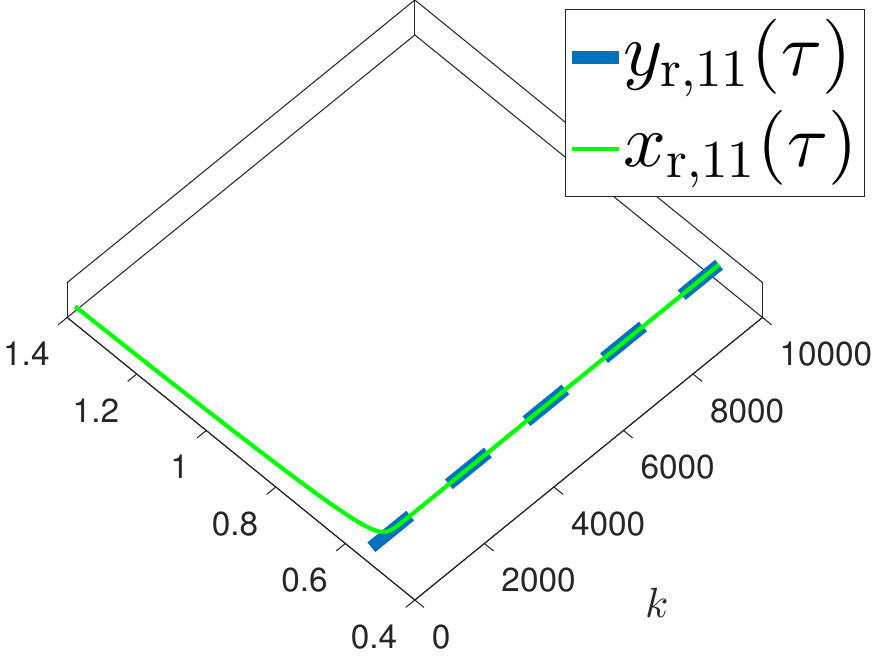}}
	\subfigure[]{\includegraphics[width=0.32\columnwidth]{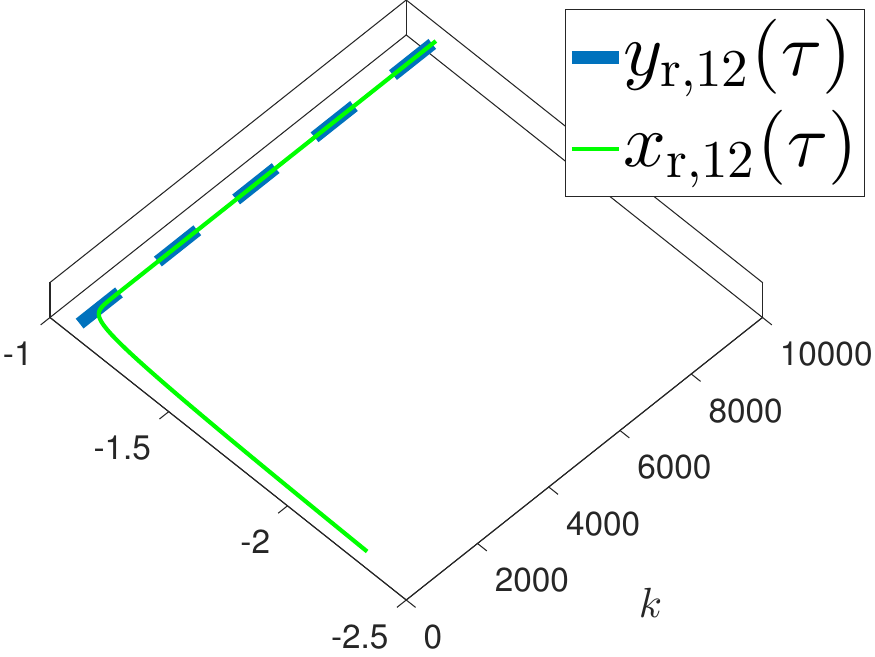}}
	\subfigure[]{\includegraphics[width=0.32\columnwidth]{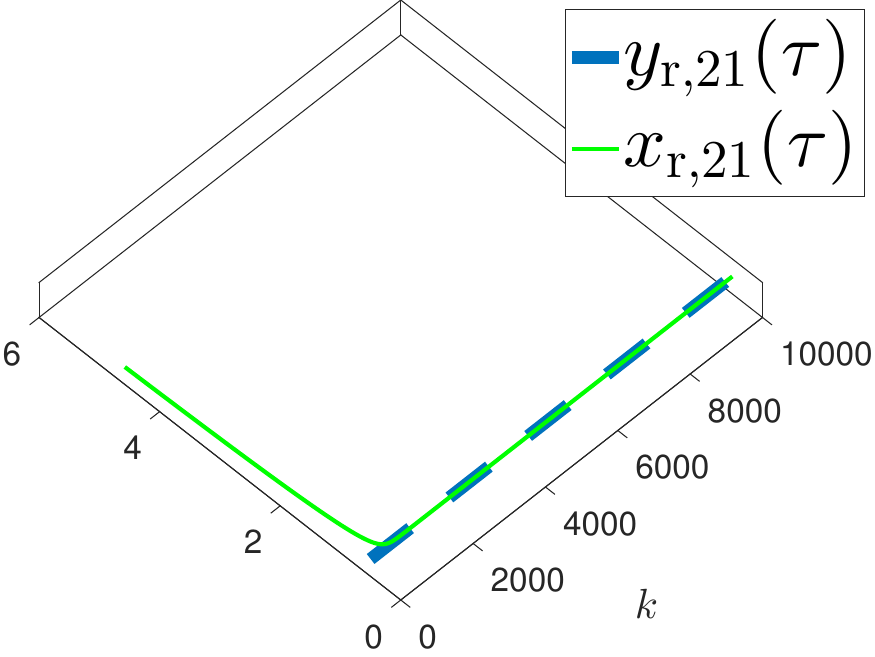}}
	\subfigure[]{\includegraphics[width=0.32\columnwidth]{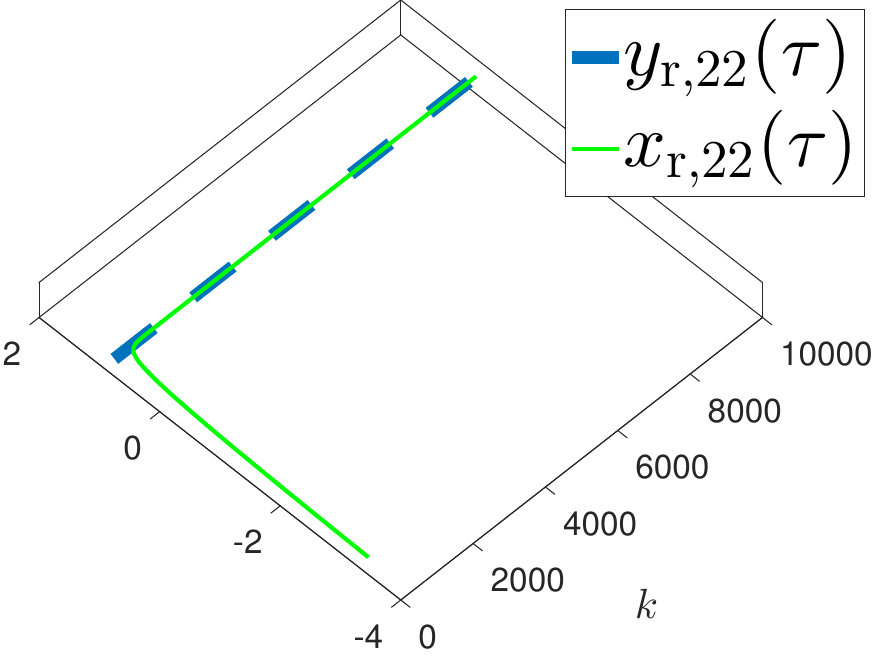}}
	\subfigure[]{\includegraphics[width=0.32\columnwidth]{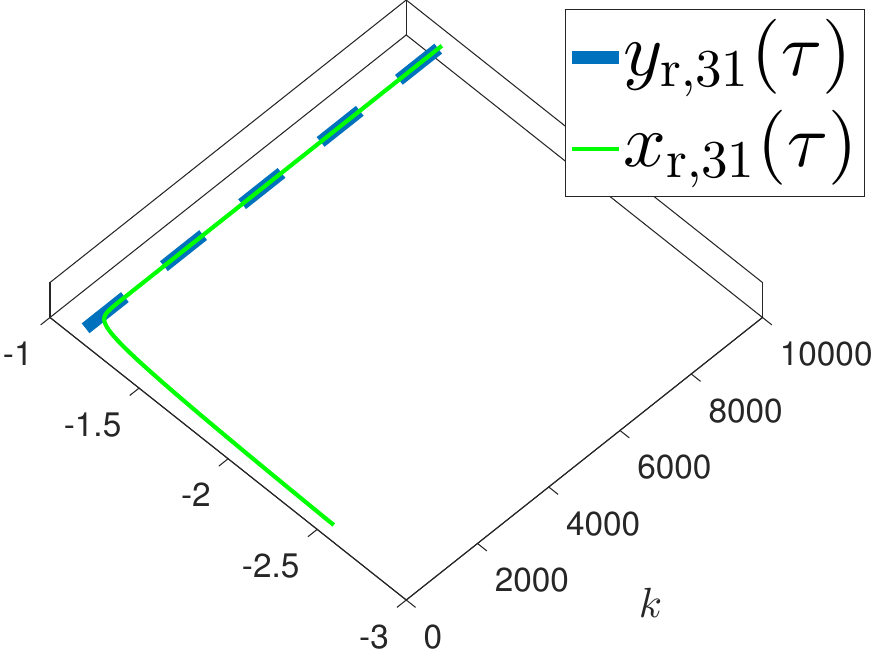}}
	\subfigure[]{\includegraphics[width=0.32\columnwidth]{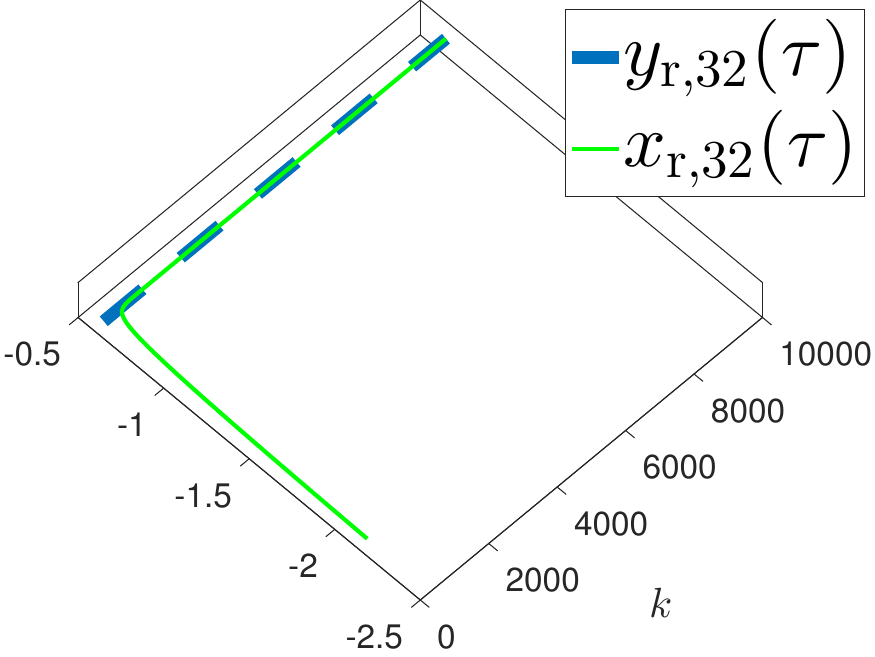}}
	\subfigure[]{\includegraphics[width=0.32\columnwidth]{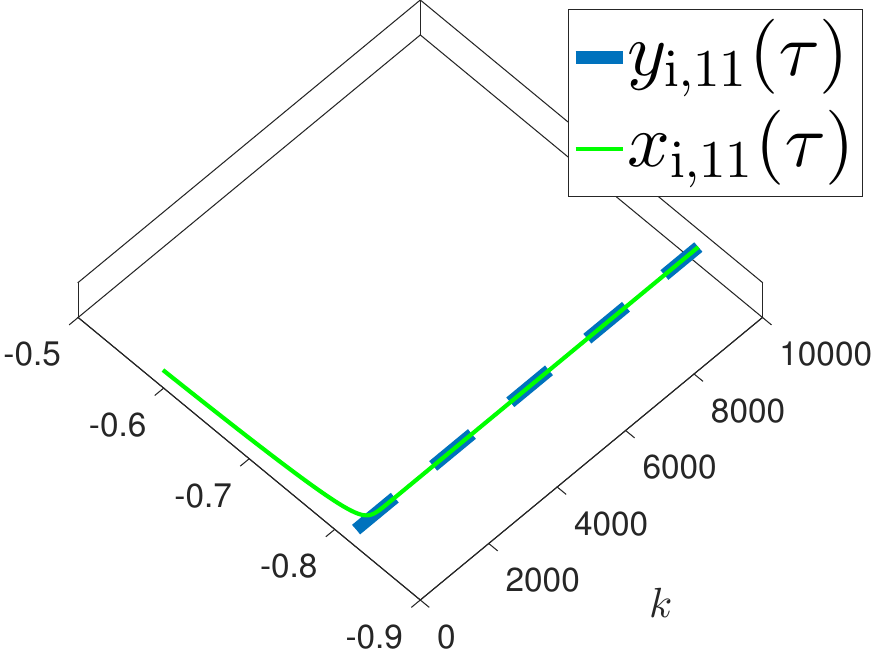}}
	\subfigure[]{\includegraphics[width=0.32\columnwidth]{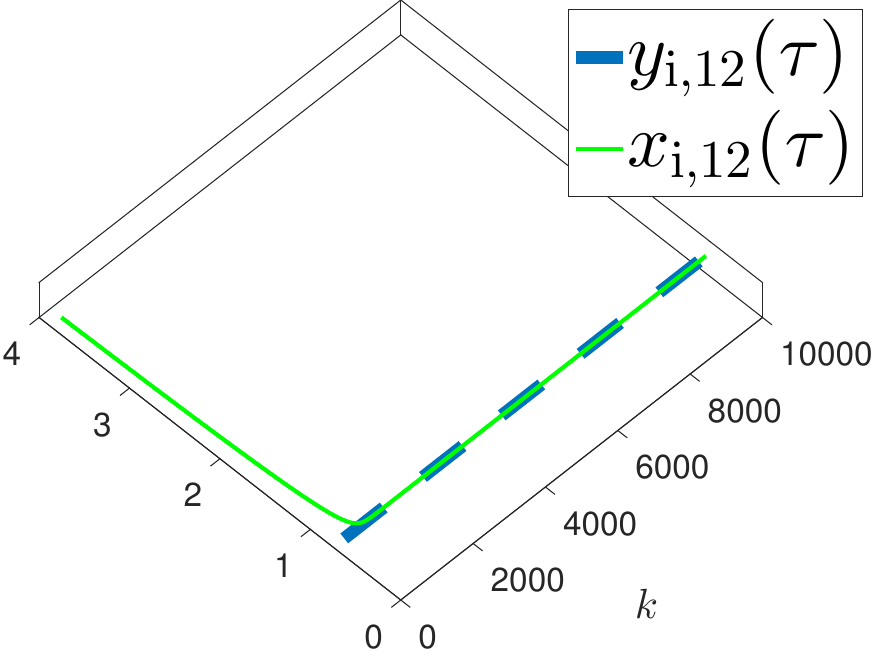}}
	\subfigure[]{\includegraphics[width=0.32\columnwidth]{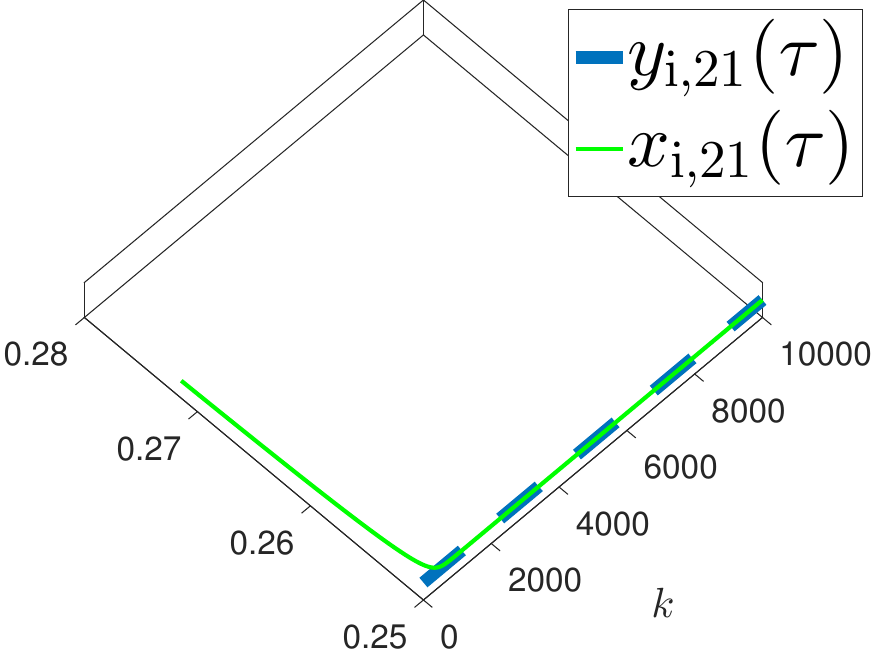}}
	\subfigure[]{\includegraphics[width=0.32\columnwidth]{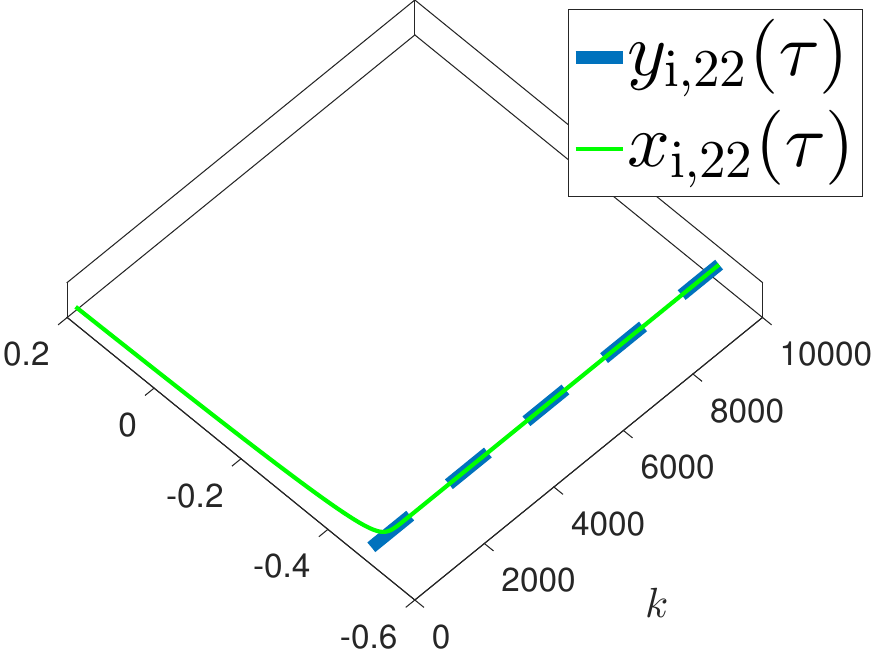}}
	\subfigure[]{\includegraphics[width=0.32\columnwidth]{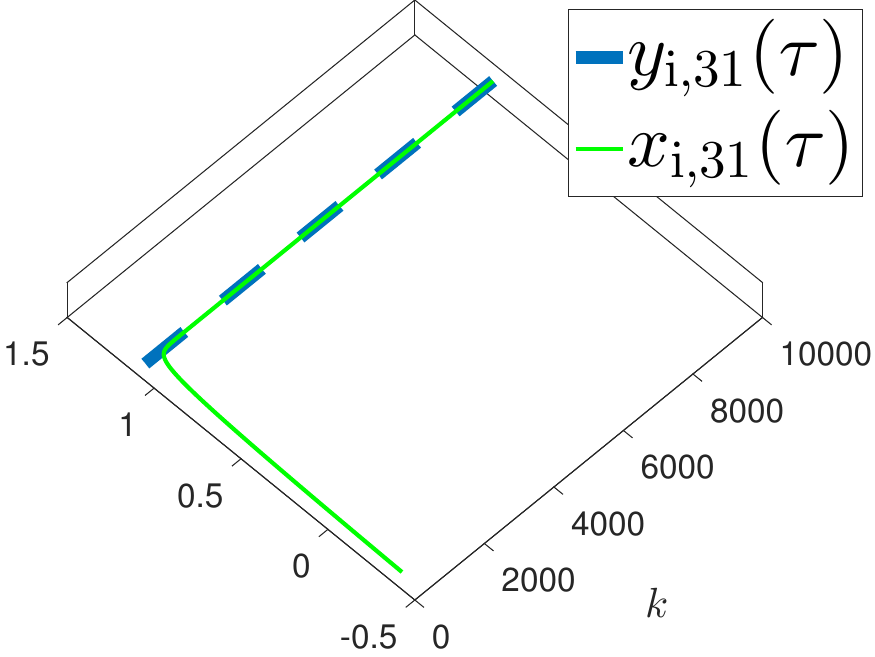}}
	\subfigure[]{\includegraphics[width=0.32\columnwidth]{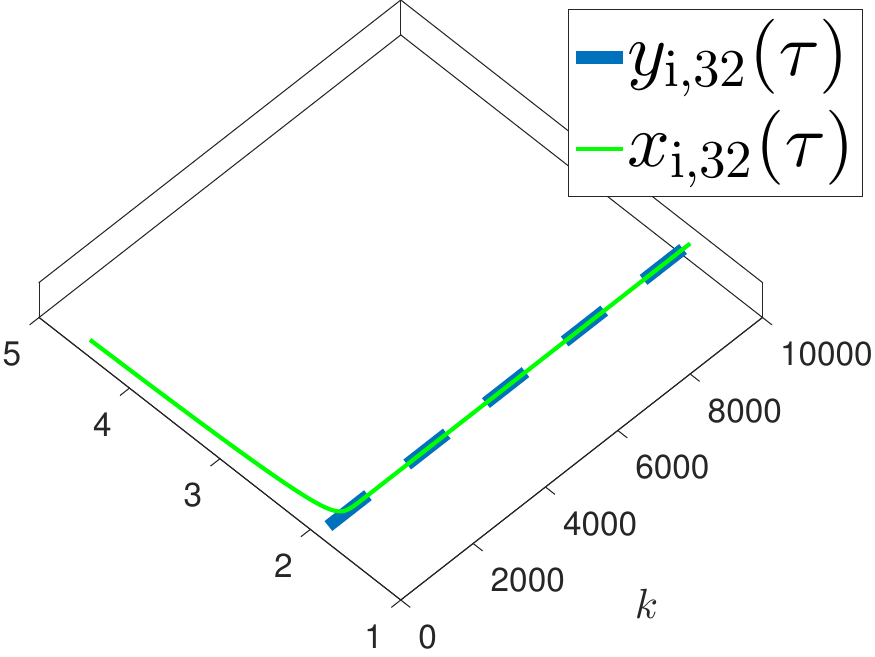}}
	\caption{Solution $X(\tau)$ computed by Con-DZND2-2i \eqref{eq.euler.forward.solve.linearerrconcznd2} model in Example \ref{example1} where $\gamma$ equals 10 and $\varepsilon$ equals 0.001.}
	\label{fig.e1.Con-DZND2-2i.solve.10.0.001}
\end{figure}

\begin{figure}[!h]\centering
	\subfigure[]{\includegraphics[width=0.70\columnwidth]{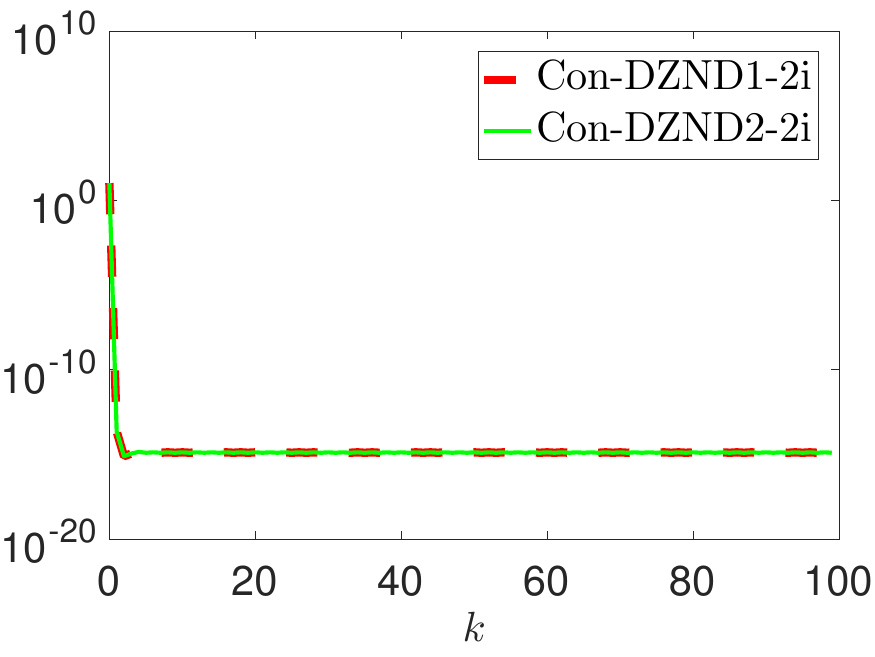}\label{fig.e1.normerror.10.0.1}}
	\subfigure[]{\includegraphics[width=0.70\columnwidth]{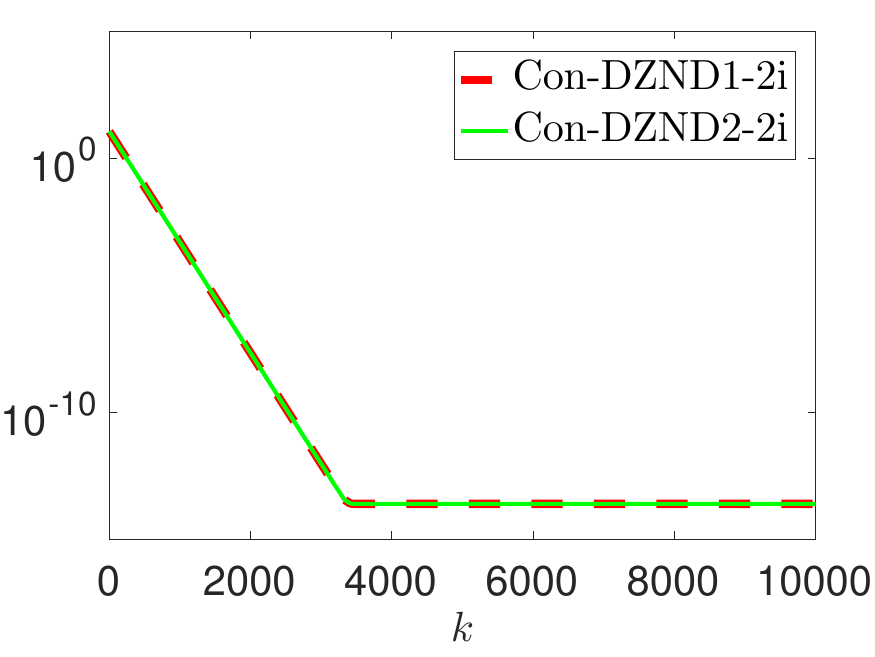}\label{fig.e1.normerror.10.0.001}}
	\caption{Logarithmic residual $\left \|X(\tau)-X^*(\tau)   \right \|_{\mathrm{F}}$ trajectories computed by Con-DZND1-2i \eqref{eq.euler.forward.solve.linearerrconcznd1} model vs. Con-DZND2-2i \eqref{eq.euler.forward.solve.linearerrconcznd2} model in Example \ref{example1} where $\gamma$ equals 10. \subref{fig.e1.normerror.10.0.1} $\varepsilon$ equals 0.1. \subref{fig.e1.normerror.10.0.001} $\varepsilon$ equals 0.001.}
	\label{fig.e1.Con-DZND1-2i.vs.Con-DZND2-2i.10}
\end{figure}

\begin{figure}[!h]\centering
	\subfigure[]{\includegraphics[width=0.375\columnwidth]{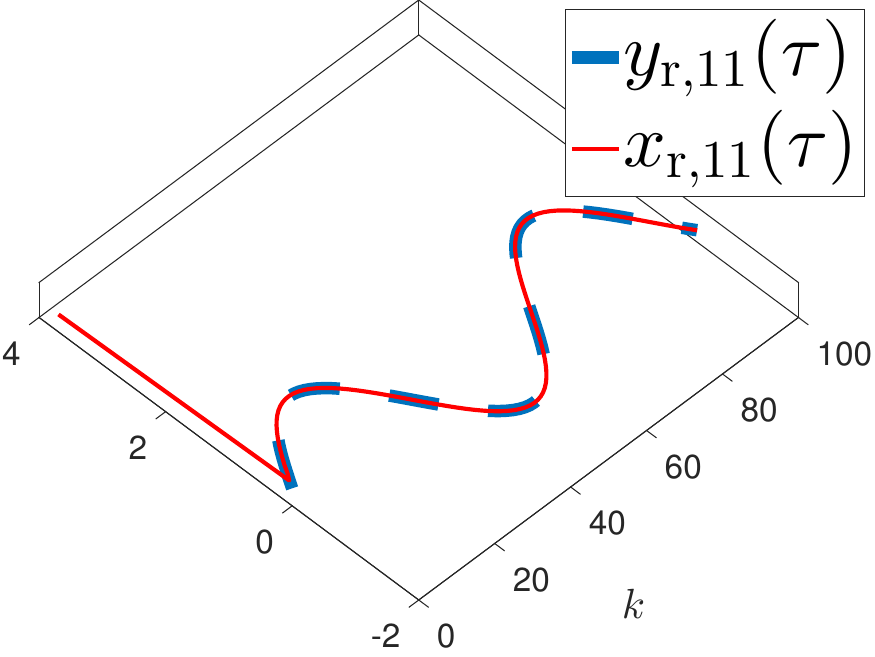}}
	\subfigure[]{\includegraphics[width=0.375\columnwidth]{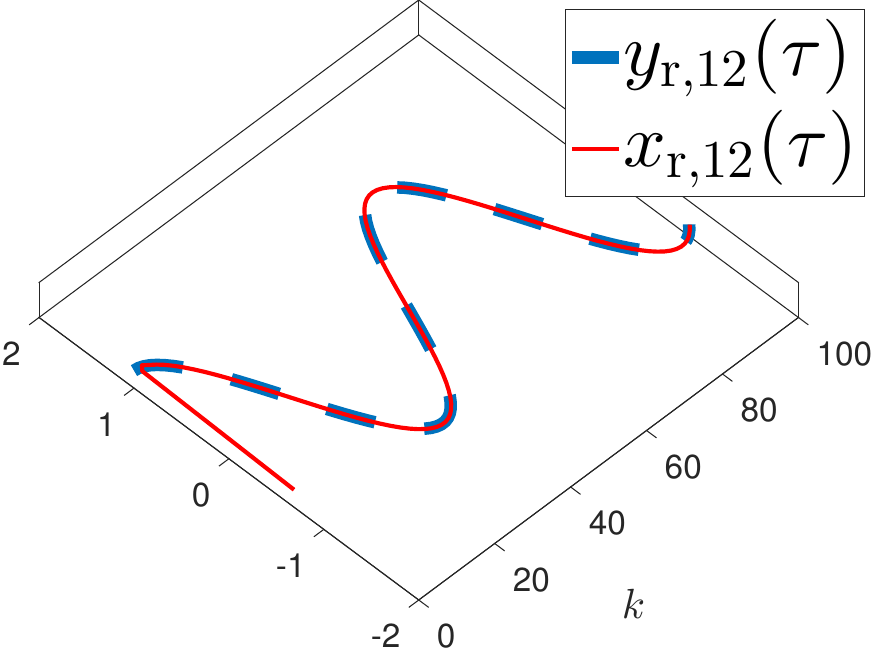}}
	\subfigure[]{\includegraphics[width=0.375\columnwidth]{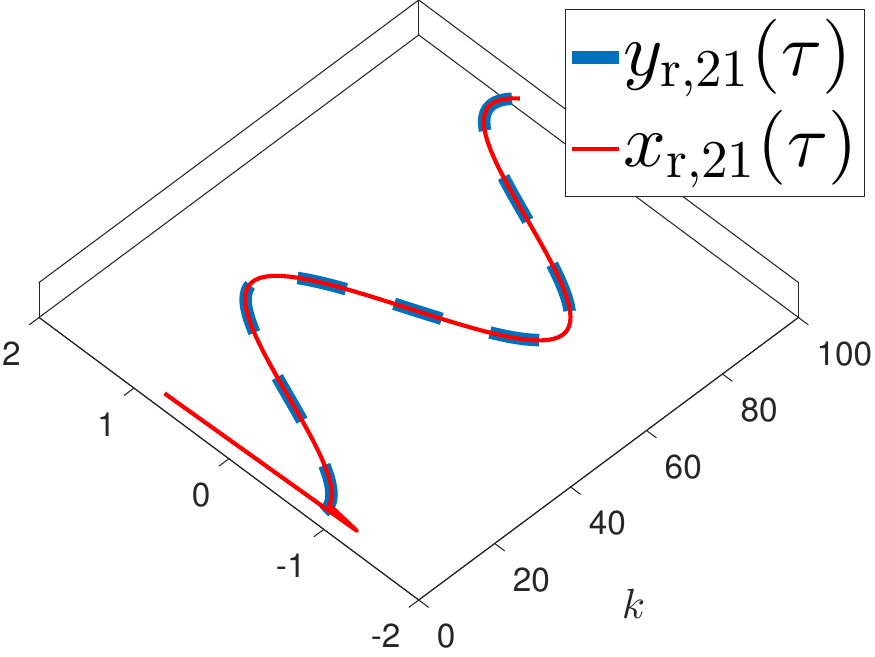}}
	\subfigure[]{\includegraphics[width=0.375\columnwidth]{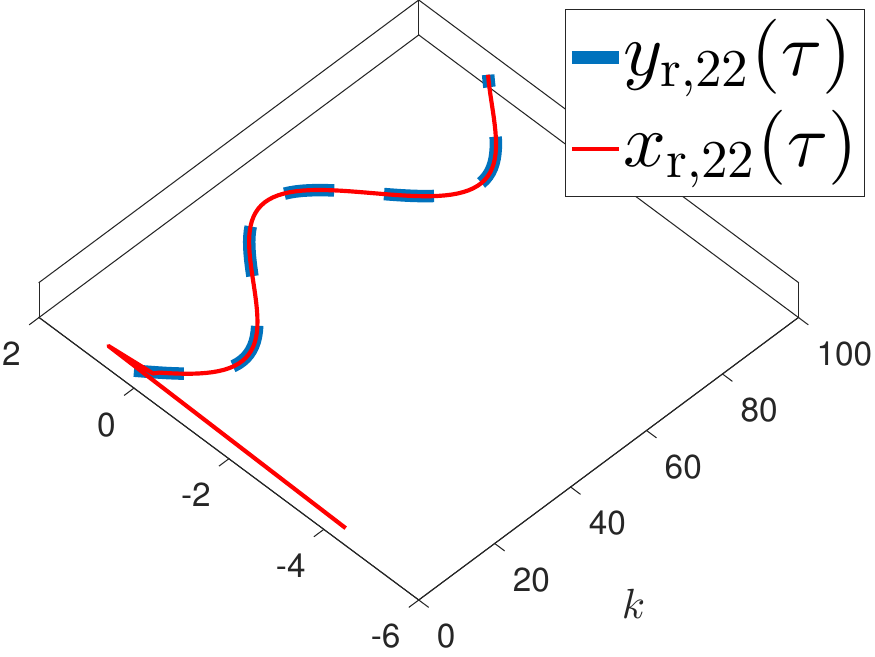}}
	\subfigure[]{\includegraphics[width=0.375\columnwidth]{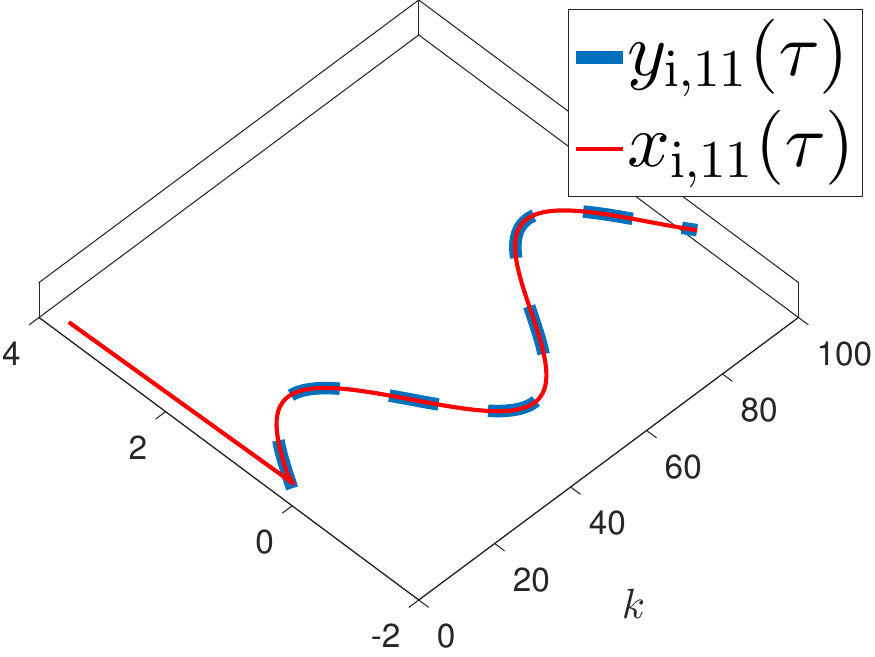}}
	\subfigure[]{\includegraphics[width=0.375\columnwidth]{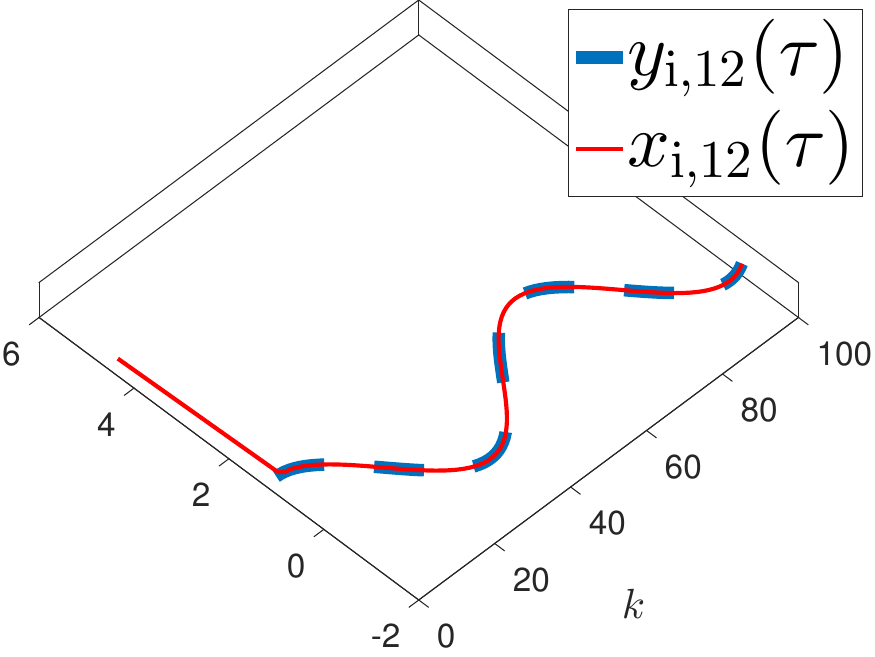}}
	\subfigure[]{\includegraphics[width=0.375\columnwidth]{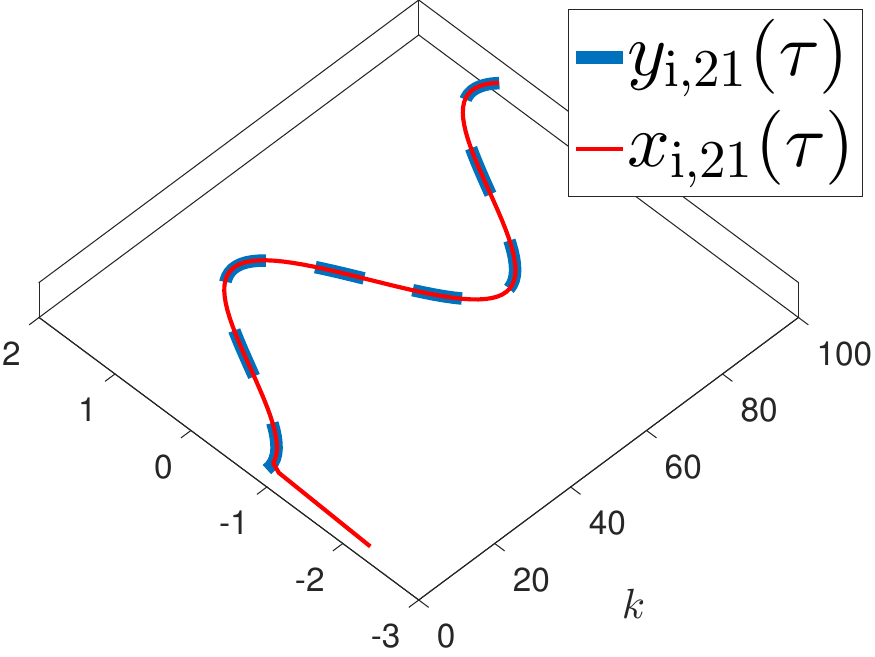}}
	\subfigure[]{\includegraphics[width=0.375\columnwidth]{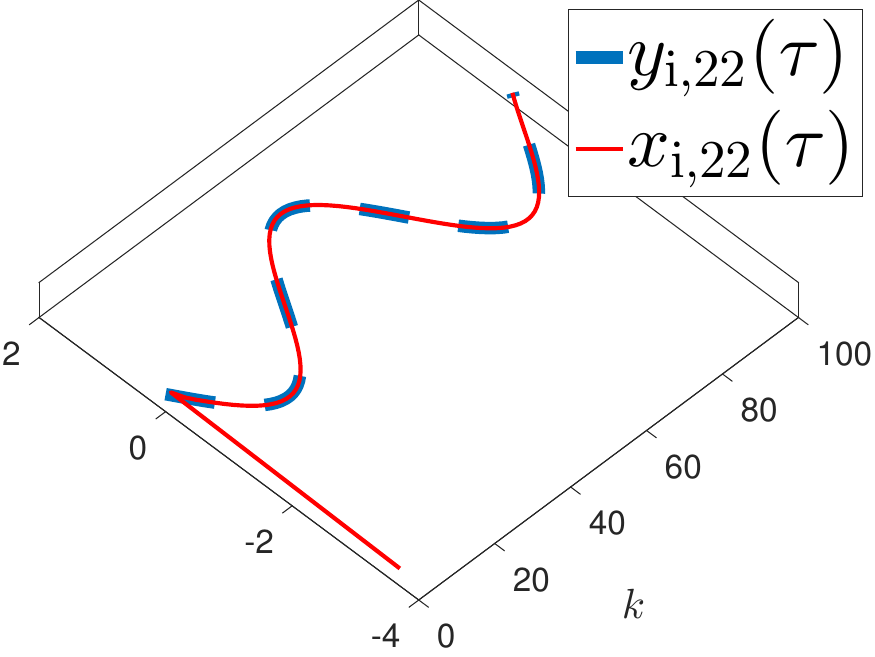}}
	\caption{Solution $X(\tau)$ computed by Con-DZND1-2i \eqref{eq.euler.forward.solve.linearerrconcznd1} model in Example \ref{example2} where $\gamma$ equals 10 and $\varepsilon$ equals 0.1.}
	\label{fig.e2.Con-DZND1-2i.solve.10.0.1}
\end{figure}

\begin{figure}[!h]\centering
\subfigure[]{\includegraphics[width=0.375\columnwidth]{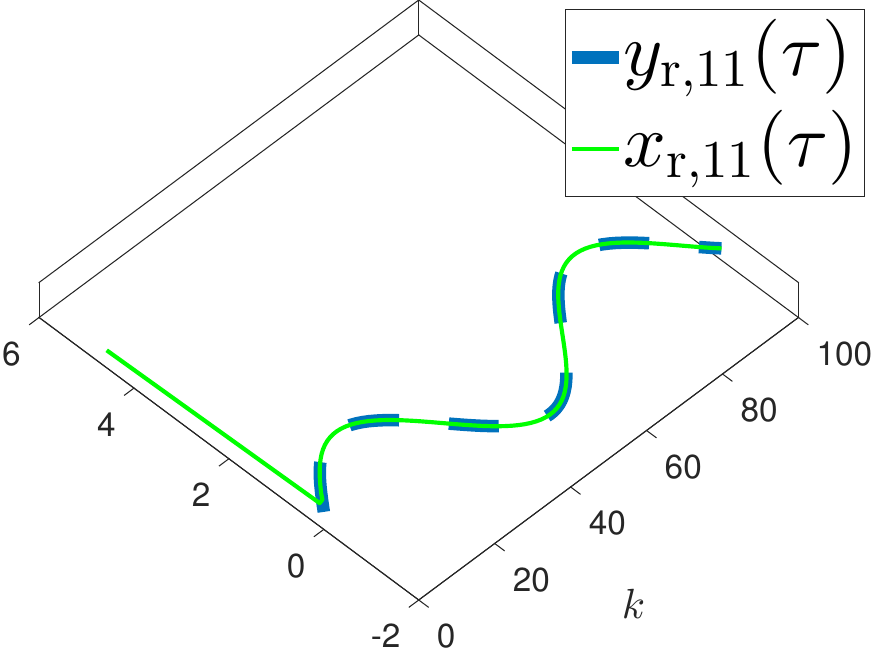}}
\subfigure[]{\includegraphics[width=0.375\columnwidth]{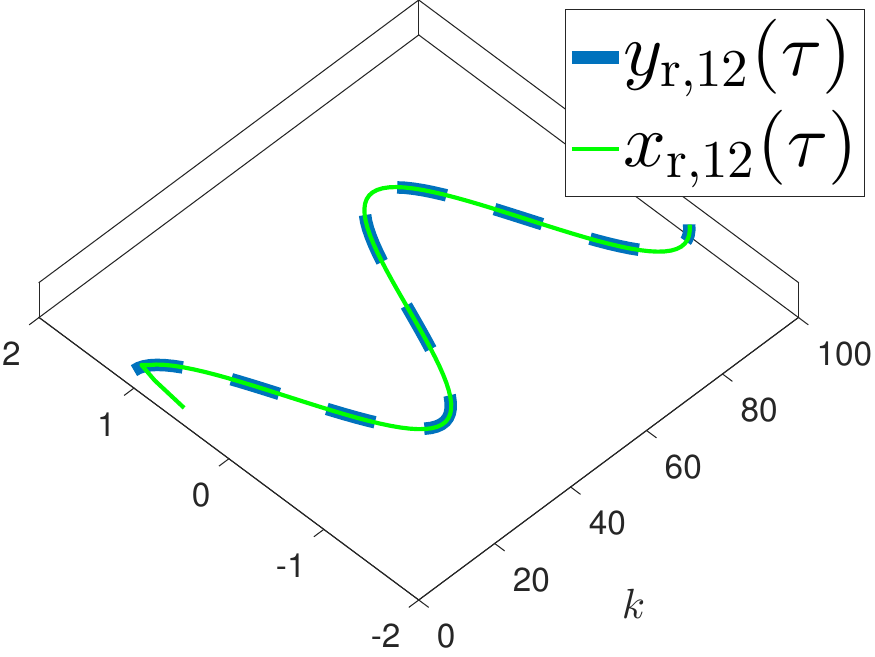}}
\subfigure[]{\includegraphics[width=0.375\columnwidth]{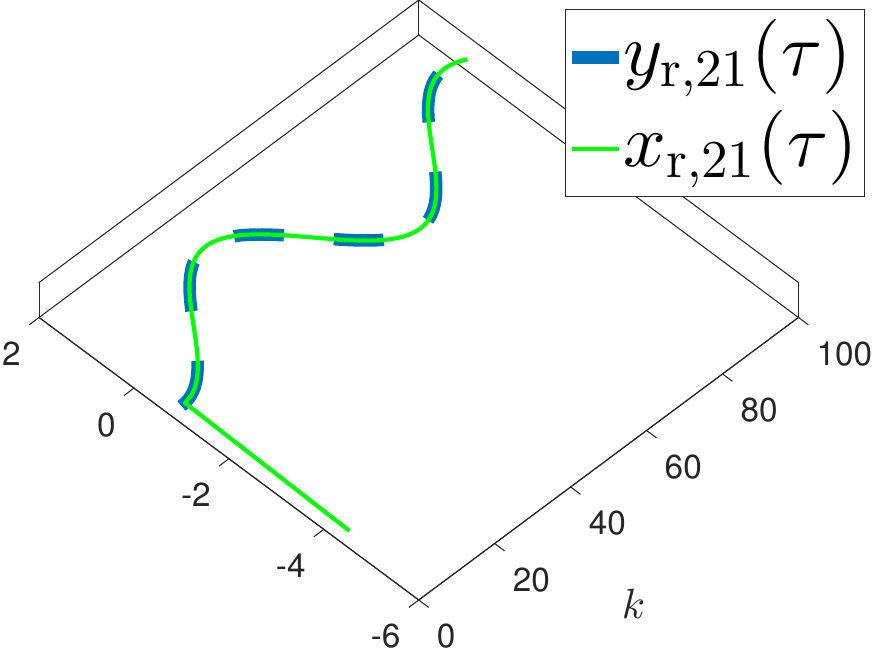}}
\subfigure[]{\includegraphics[width=0.375\columnwidth]{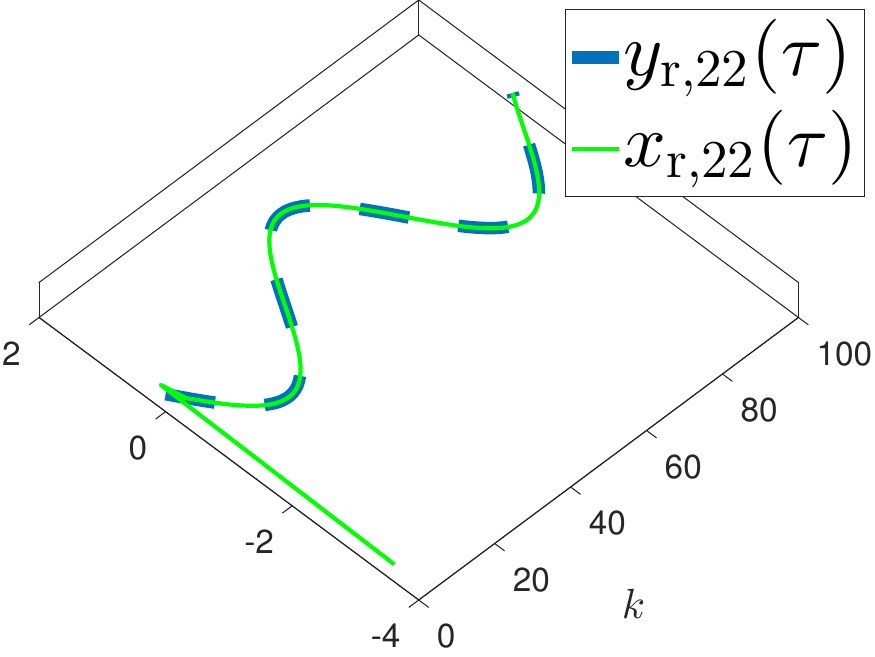}}
\subfigure[]{\includegraphics[width=0.375\columnwidth]{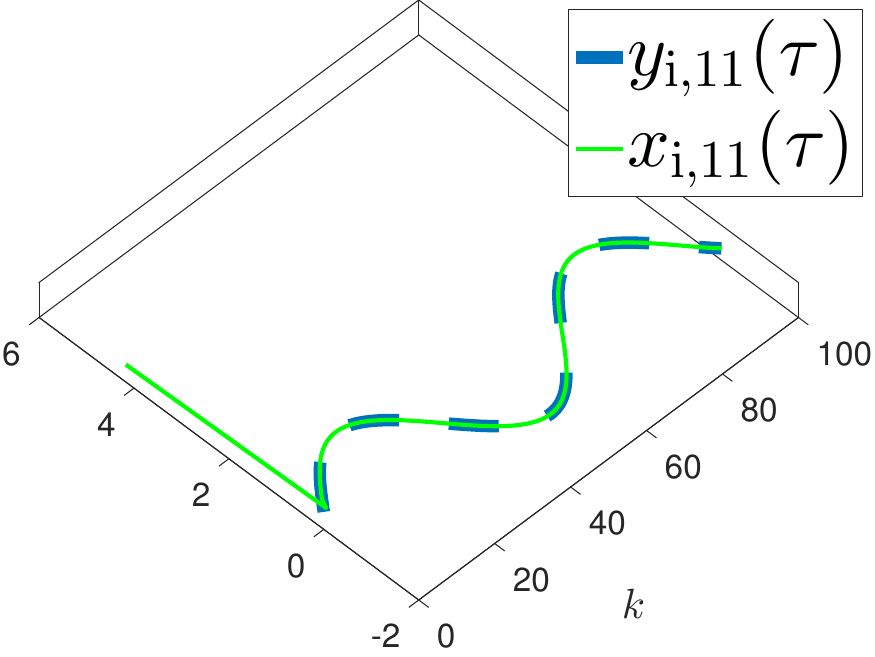}}
\subfigure[]{\includegraphics[width=0.375\columnwidth]{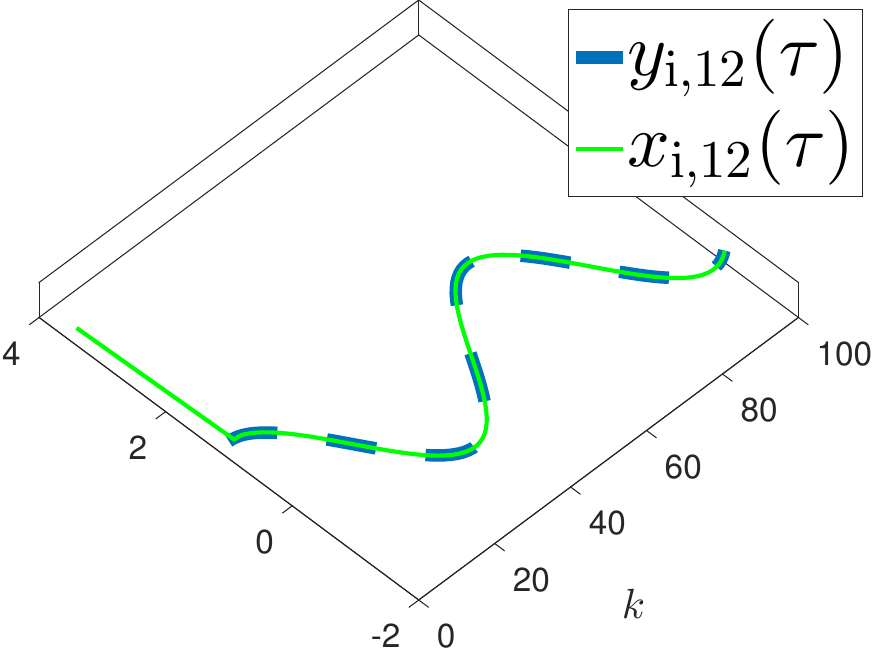}}
\subfigure[]{\includegraphics[width=0.375\columnwidth]{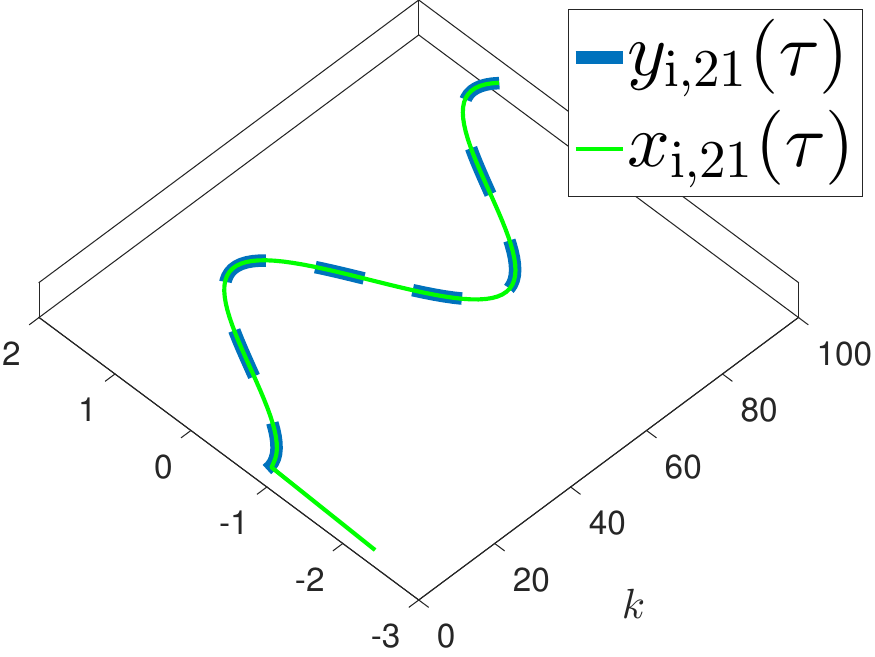}}
\subfigure[]{\includegraphics[width=0.375\columnwidth]{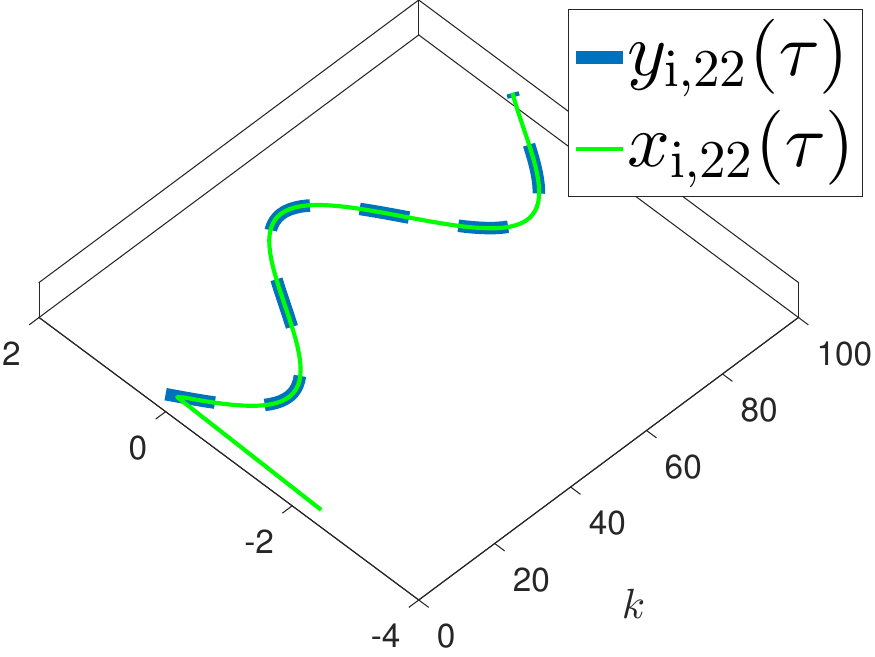}}
	\caption{Solution $X(\tau)$ computed by Con-DZND2-2i \eqref{eq.euler.forward.solve.linearerrconcznd2} model in Example \ref{example2} where $\gamma$ equals 10 and $\varepsilon$ equals 0.1.}
	\label{fig.e2.Con-DZND2-2i.solve.10.0.1}
\end{figure}

\begin{figure}[!h]\centering
\subfigure[]{\includegraphics[width=0.375\columnwidth]{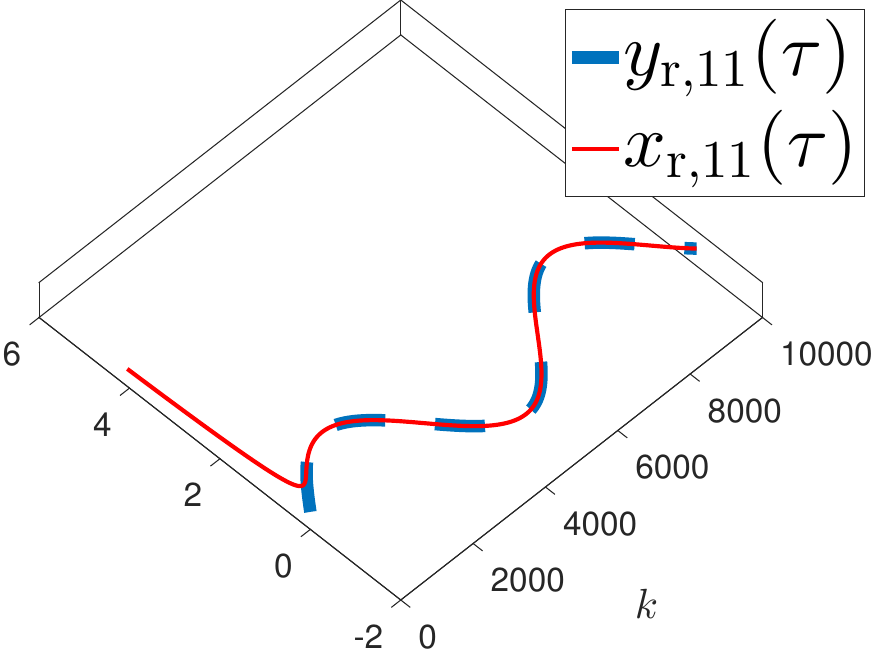}}
\subfigure[]{\includegraphics[width=0.375\columnwidth]{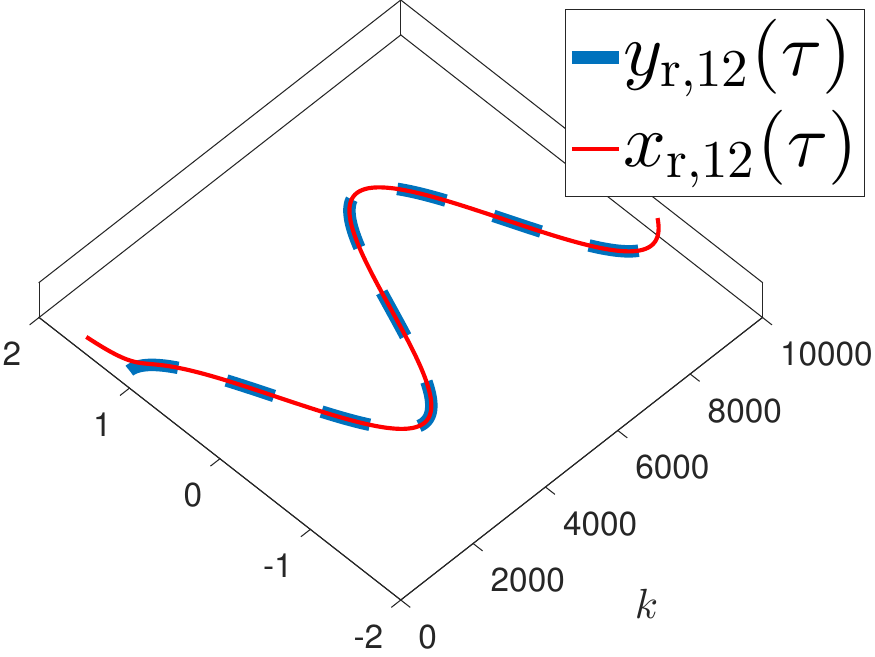}}
\subfigure[]{\includegraphics[width=0.375\columnwidth]{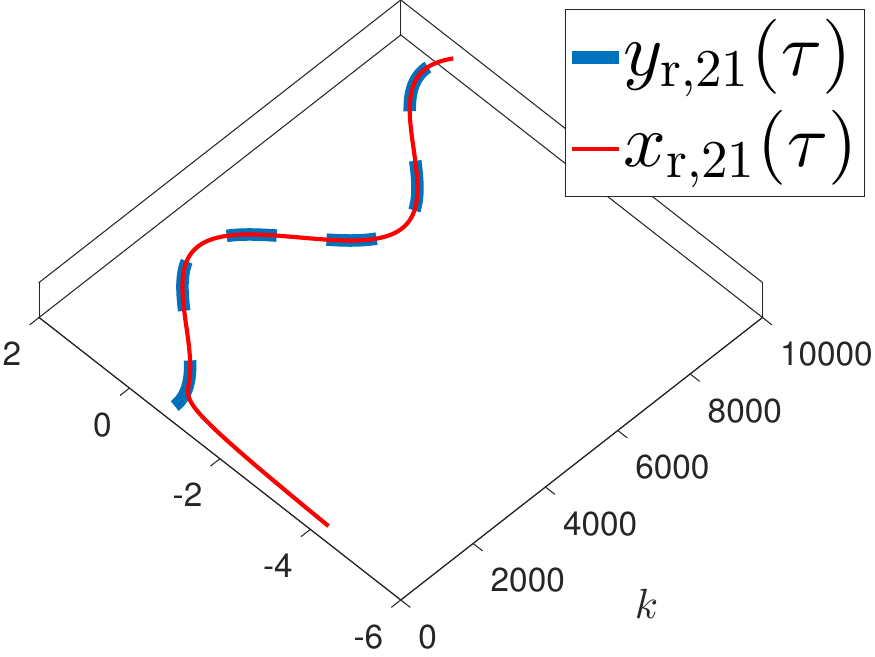}}
\subfigure[]{\includegraphics[width=0.375\columnwidth]{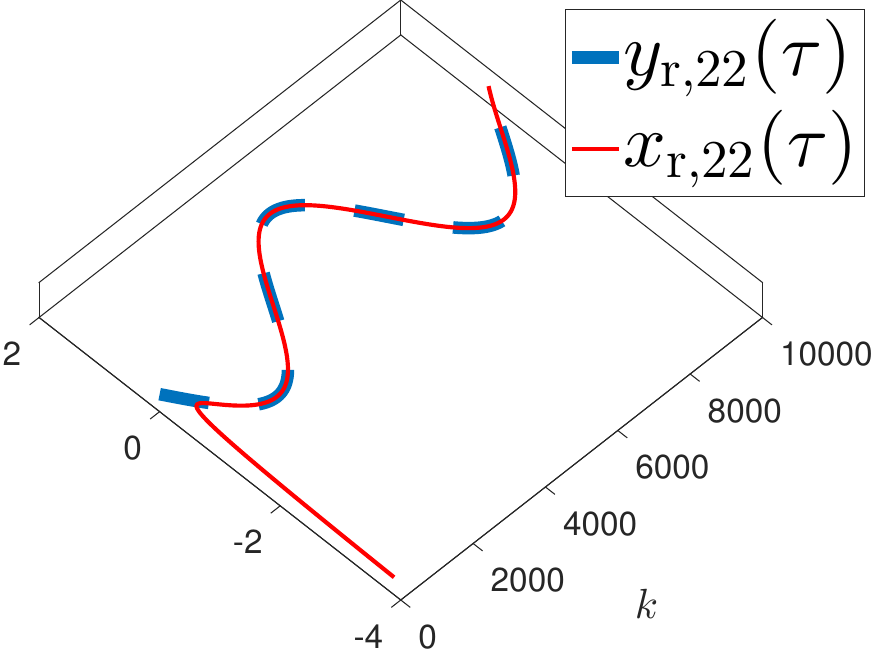}}
\subfigure[]{\includegraphics[width=0.375\columnwidth]{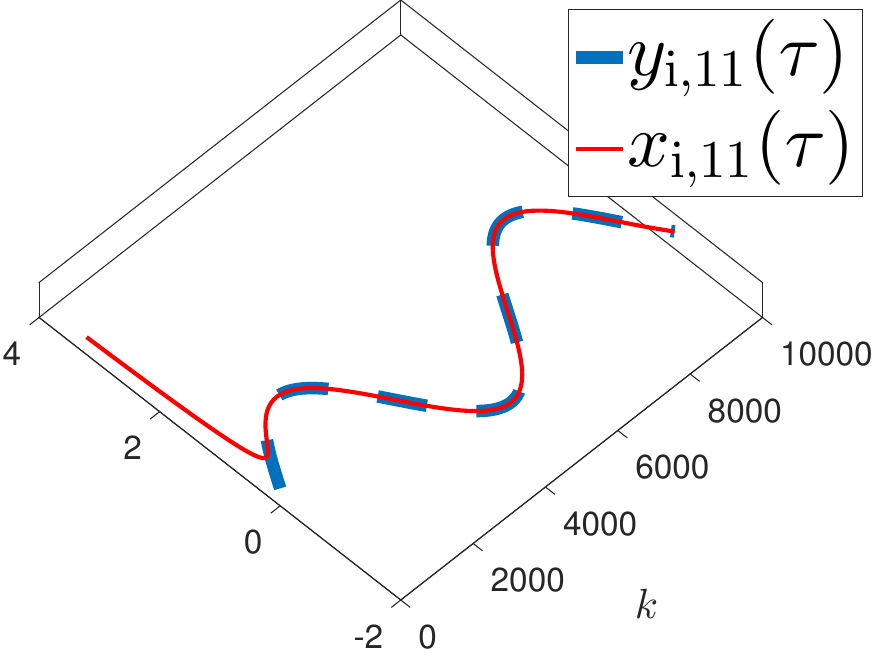}}
\subfigure[]{\includegraphics[width=0.375\columnwidth]{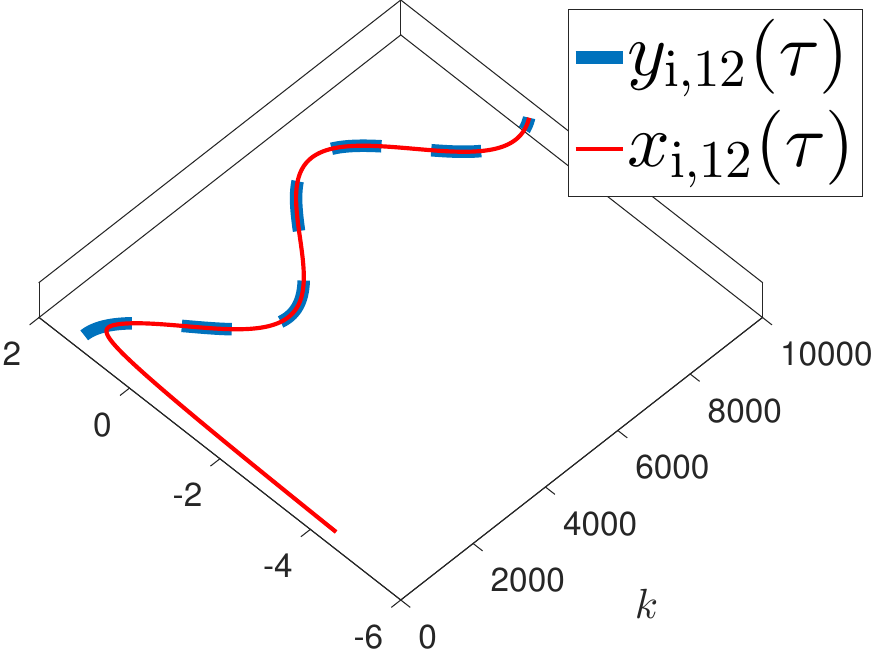}}
\subfigure[]{\includegraphics[width=0.375\columnwidth]{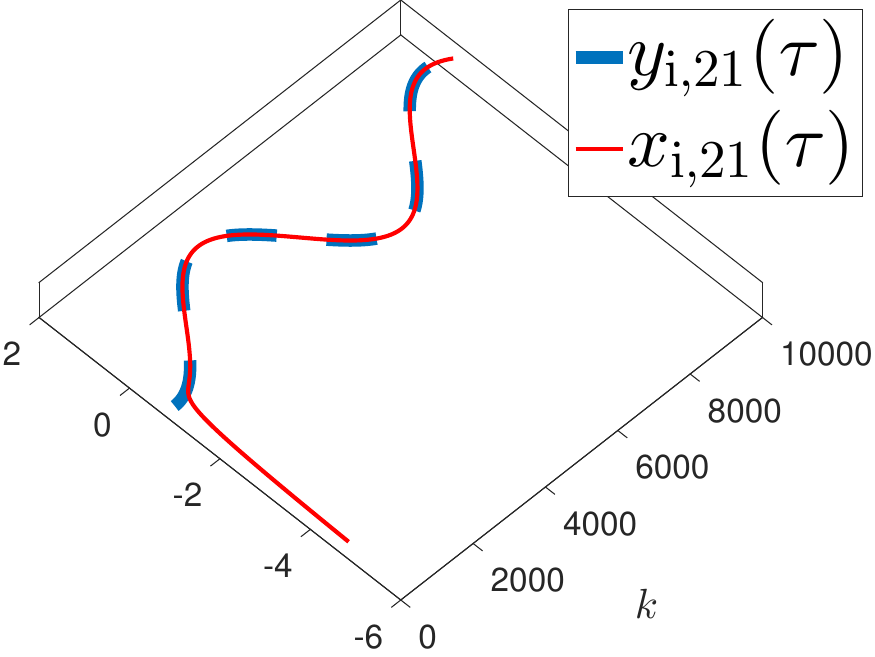}}
\subfigure[]{\includegraphics[width=0.375\columnwidth]{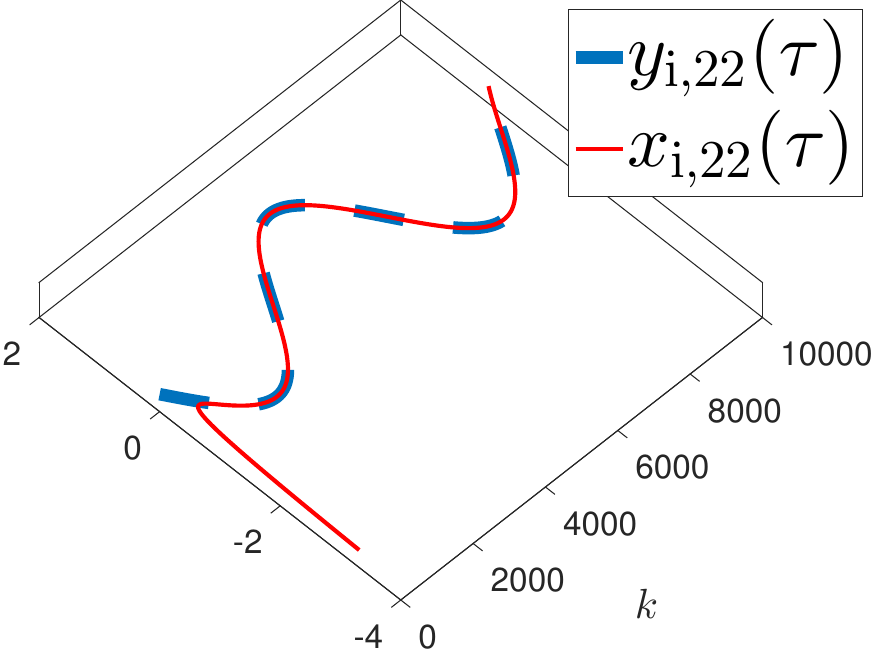}}
	\caption{Solution $X(\tau)$ computed by Con-DZND1-2i \eqref{eq.euler.forward.solve.linearerrconcznd1} model in Example \ref{example2} where $\gamma$ equals 10 and $\varepsilon$ equals 0.001.}
	\label{fig.e2.Con-DZND1-2i.solve.10.0.001}
\end{figure}

\begin{figure}[!h]\centering
\subfigure[]{\includegraphics[width=0.375\columnwidth]{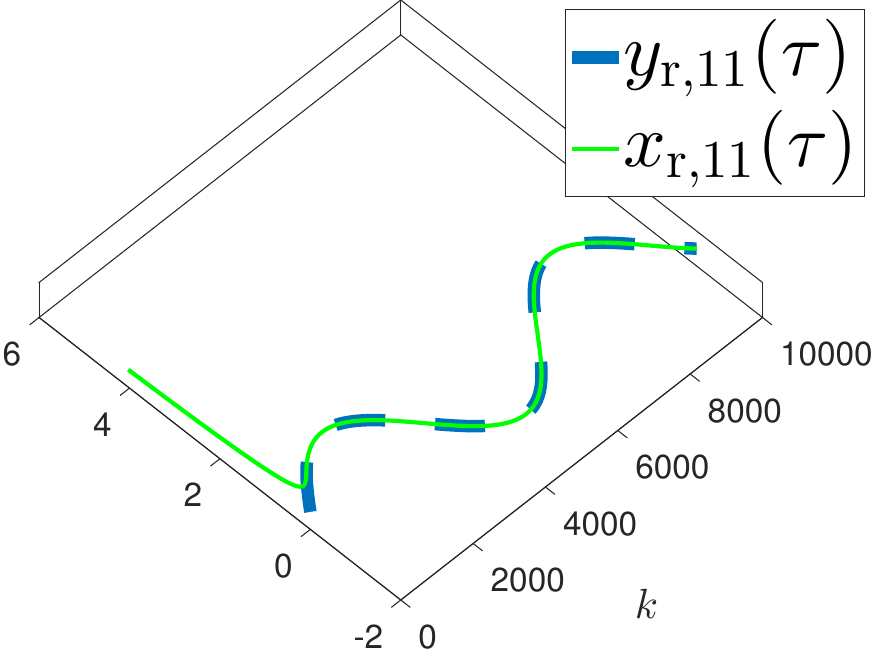}}
\subfigure[]{\includegraphics[width=0.375\columnwidth]{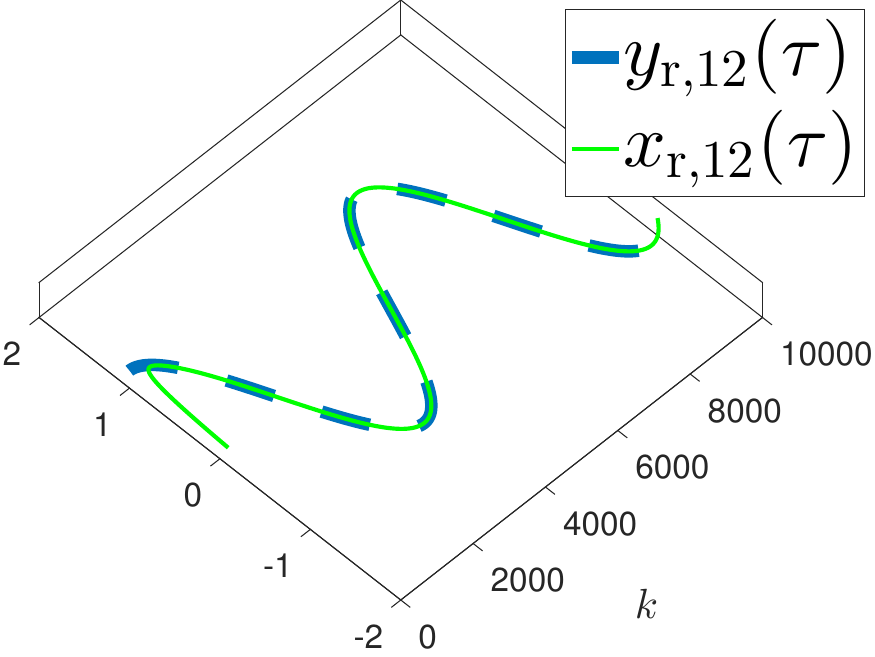}}
\subfigure[]{\includegraphics[width=0.375\columnwidth]{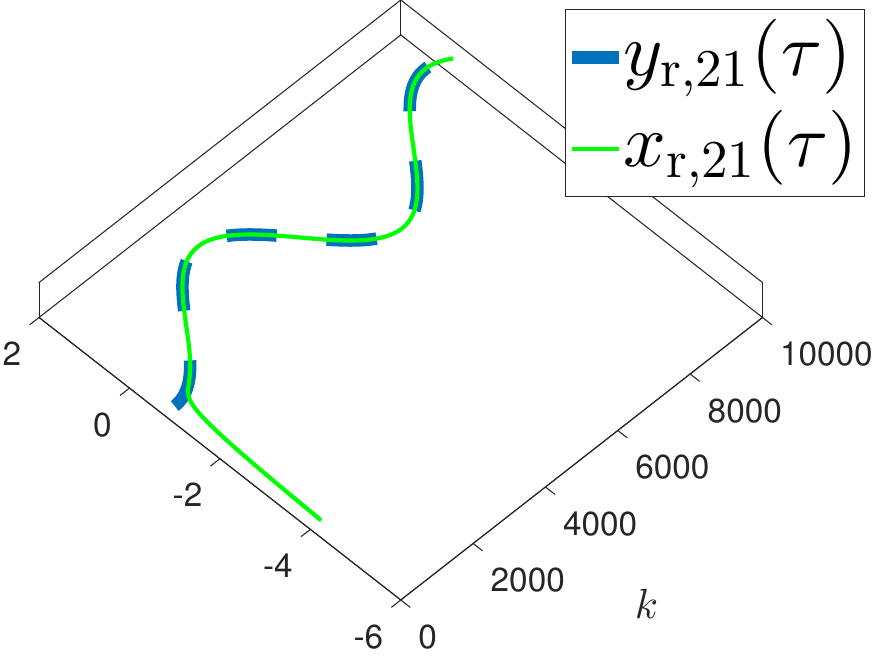}}
\subfigure[]{\includegraphics[width=0.375\columnwidth]{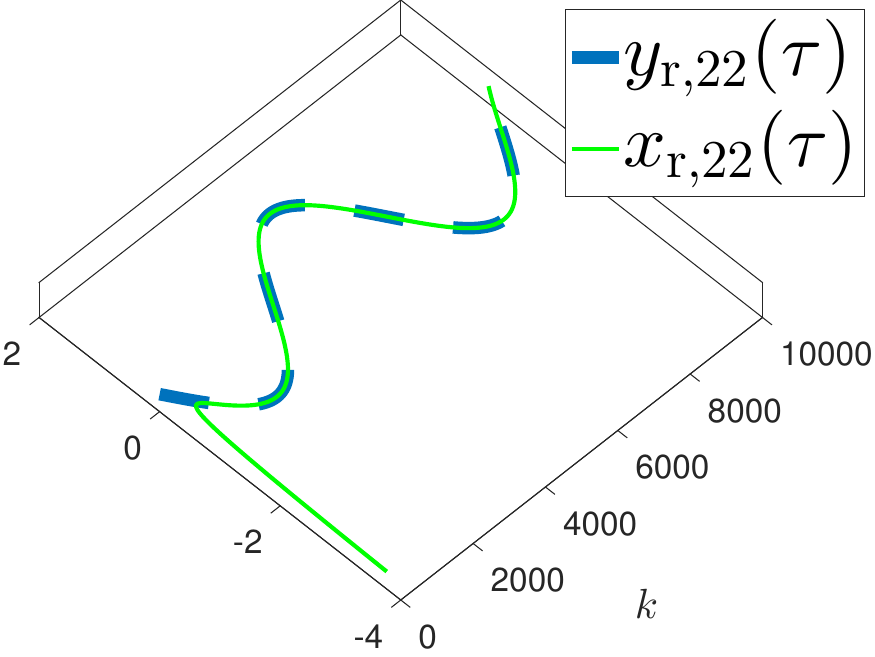}}
\subfigure[]{\includegraphics[width=0.375\columnwidth]{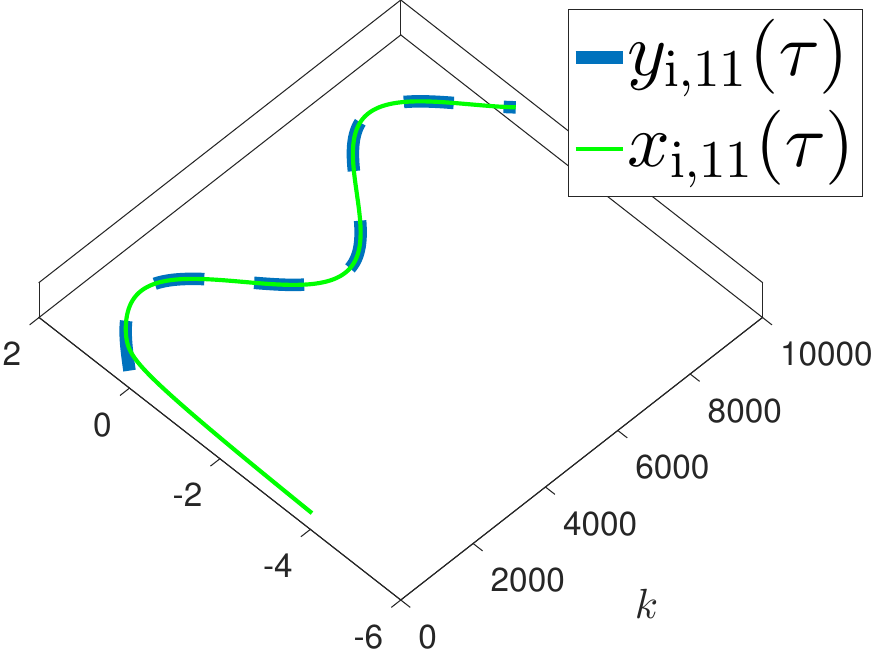}}
\subfigure[]{\includegraphics[width=0.375\columnwidth]{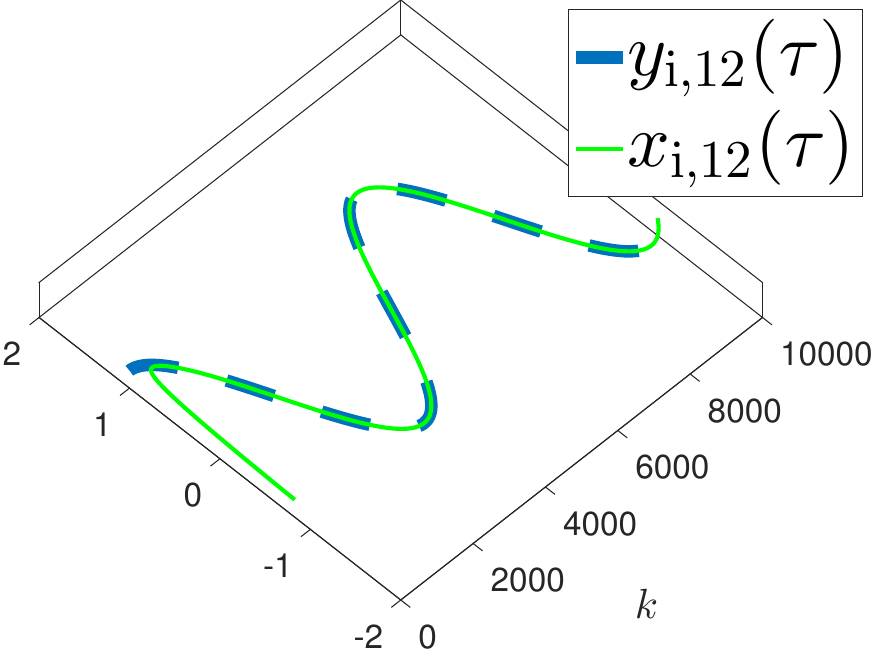}}
\subfigure[]{\includegraphics[width=0.375\columnwidth]{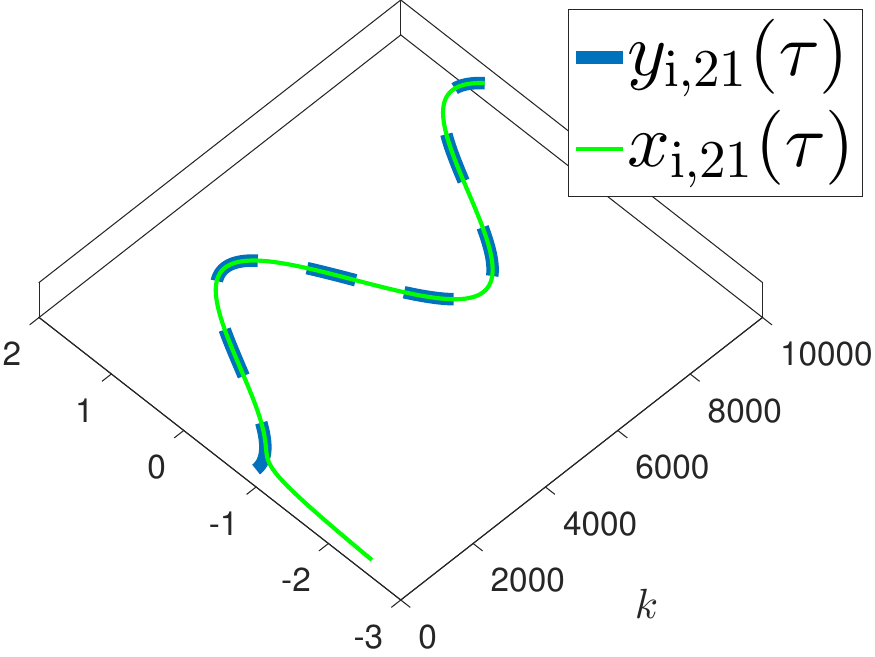}}
\subfigure[]{\includegraphics[width=0.375\columnwidth]{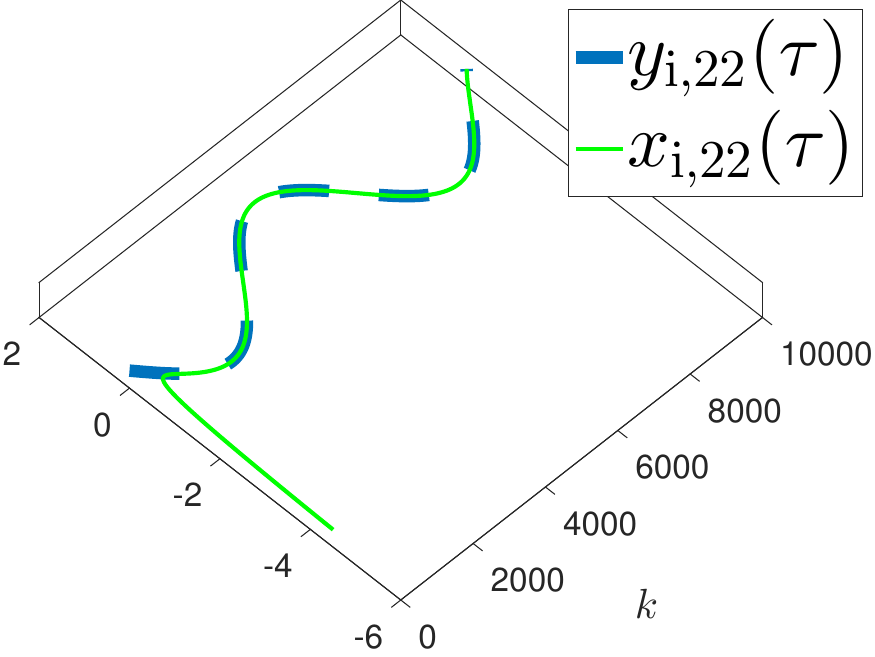}}
	\caption{Solution $X(\tau)$ computed by Con-DZND2-2i \eqref{eq.euler.forward.solve.linearerrconcznd2} model in Example \ref{example2} where $\gamma$ equals 10 and $\varepsilon$ equals 0.001.}
	\label{fig.e2.Con-DZND2-2i.solve.10.0.001}
\end{figure}

\begin{figure}[!h]\centering
	\subfigure[]{\includegraphics[width=0.70\columnwidth]{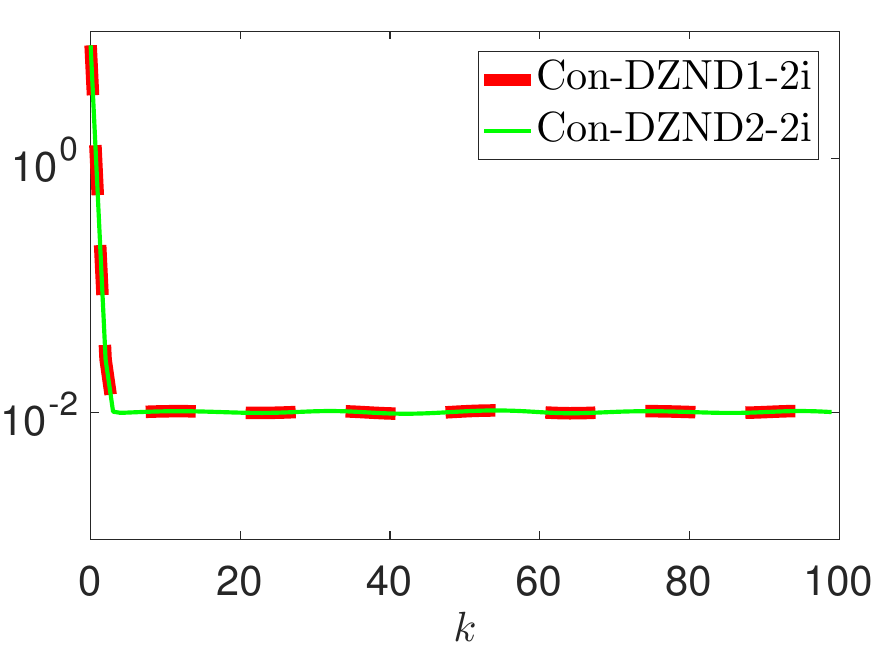}\label{fig.e2.normerror.10.0.1}}
	\subfigure[]{\includegraphics[width=0.70\columnwidth]{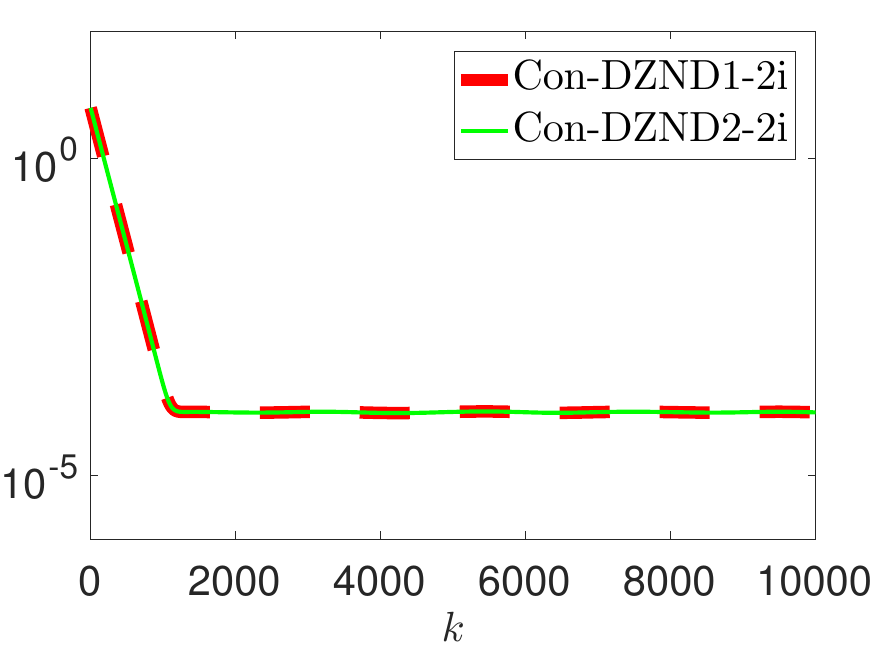}\label{fig.e2.normerror.10.0.001}}
	\caption{Logarithmic residual $\left \|X(\tau)-X^*(\tau)   \right \|_{\mathrm{F}}$ trajectories computed by Con-DZND1-2i \eqref{eq.euler.forward.solve.linearerrconcznd1} model vs. Con-DZND2-2i \eqref{eq.euler.forward.solve.linearerrconcznd2} model in Example \ref{example2} where $\gamma$ equals 10. \subref{fig.e2.normerror.10.0.1} $\varepsilon$ equals 0.1. \subref{fig.e2.normerror.10.0.001} $\varepsilon$ equals 0.001.}
	\label{fig.e2.Con-DZND1-2i.vs.Con-DZND2-2i.10}
\end{figure}

The results of Con-DZND1-2i \eqref{eq.euler.forward.solve.linearerrconcznd1} model and Con-DZND2-2i \eqref{eq.euler.forward.solve.linearerrconcznd2} model after discretion using Euler-forward formula are shown in Figs. \ref{fig.e1.Con-DZND1-2i.solve.10.0.1}
through
\ref{fig.e2.Con-DZND1-2i.vs.Con-DZND2-2i.10}.

From Figs. \ref{fig.e1.Con-DZND1-2i.solve.10.0.1}
through
\ref{fig.e1.Con-DZND1-2i.vs.Con-DZND2-2i.10}, it can be seen that in Example \ref{example1}, the high dimension error of Con-DZND2-2i \eqref{eq.euler.forward.solve.linearerrconcznd2} model is eradicated, and logarithmic residual $\left \|X(\tau)-X^*(\tau)   \right \|_{\mathrm{F}}$ trajectories are essentially consistent with those of Con-DZND1-2i \eqref{eq.euler.forward.solve.linearerrconcznd1} model, both of which are influenced by different step size $\varepsilon$.


From Figs. \ref{fig.e2.Con-DZND1-2i.solve.10.0.1}
through
\ref{fig.e2.Con-DZND1-2i.vs.Con-DZND2-2i.10}, it can be seen that in Example \ref{example2}, the high dimension error of Con-DZND2-2i \eqref{eq.euler.forward.solve.linearerrconcznd2} model is eradicated, and logarithmic residual $\left \|X(\tau)-X^*(\tau)   \right \|_{\mathrm{F}}$ trajectories are essentially consistent with those of Con-DZND1-2i \eqref{eq.euler.forward.solve.linearerrconcznd1} model, both of which are influenced by different step sizes $\varepsilon$.

However, Con-DZND1-2i \eqref{eq.euler.forward.solve.linearerrconcznd1} model essentially involves complex field error approximation. Is there space compressive approximation? Next, for Con-DZND1-2i \eqref{eq.euler.forward.solve.linearerrconcznd1} model, 
$\gamma$ will be extended to 10$+$20$\mathrm{i}$ and 10$-$20$\mathrm{i}$ to test the results of space compressive approximation and observe the convergence results of above model under different step sizes $\varepsilon$.
\subsection{The difference in Con-DZND1-2i model using step sizes $\varepsilon$ of 0.1 and 0.001 for $\gamma$ equals 10$+$20$\mathrm{i}$ and 10$-$20$\mathrm{i}$ respectively}
\begin{figure}[!h]\centering
	\subfigure[]{\includegraphics[width=0.40\columnwidth]{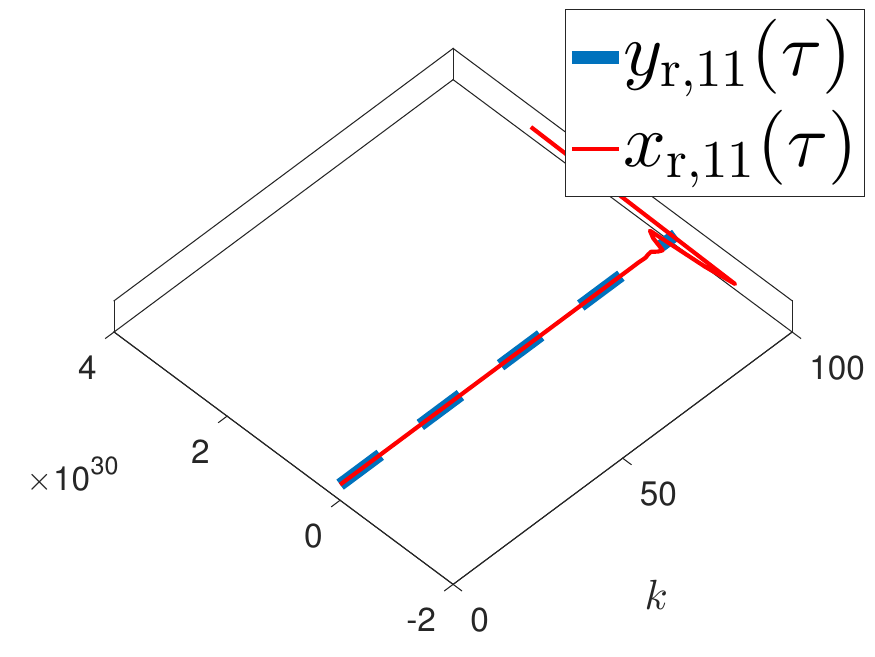}\label{fig.e1.Con-DZND1-2i.solve.10+20i.0.1.x11r.3d}}
	\subfigure[]{\includegraphics[width=0.40\columnwidth]{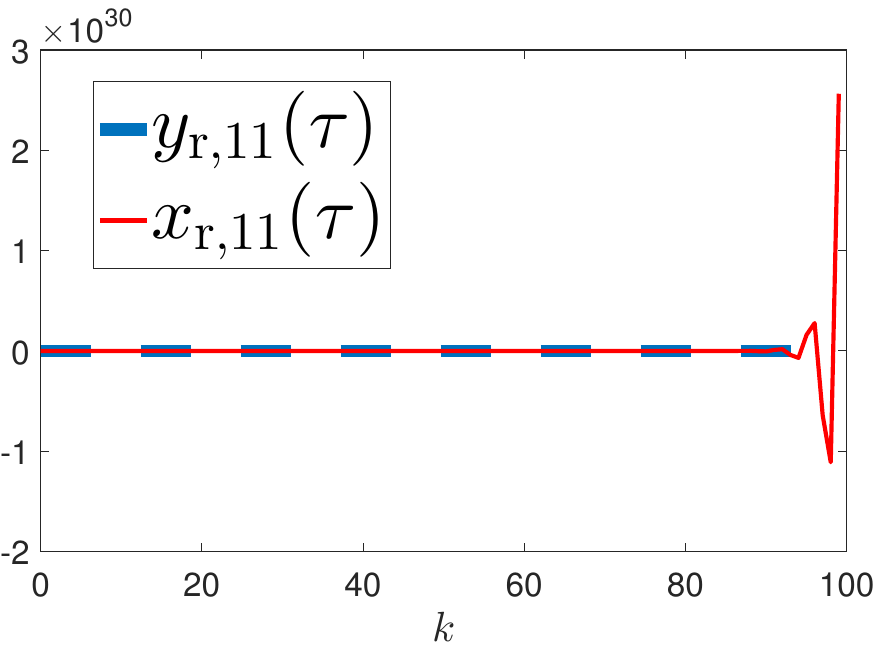}\label{fig.e1.Con-DZND1-2i.solve.10+20i.0.1.x11r.2d}}
	\subfigure[]{\includegraphics[width=0.40\columnwidth]{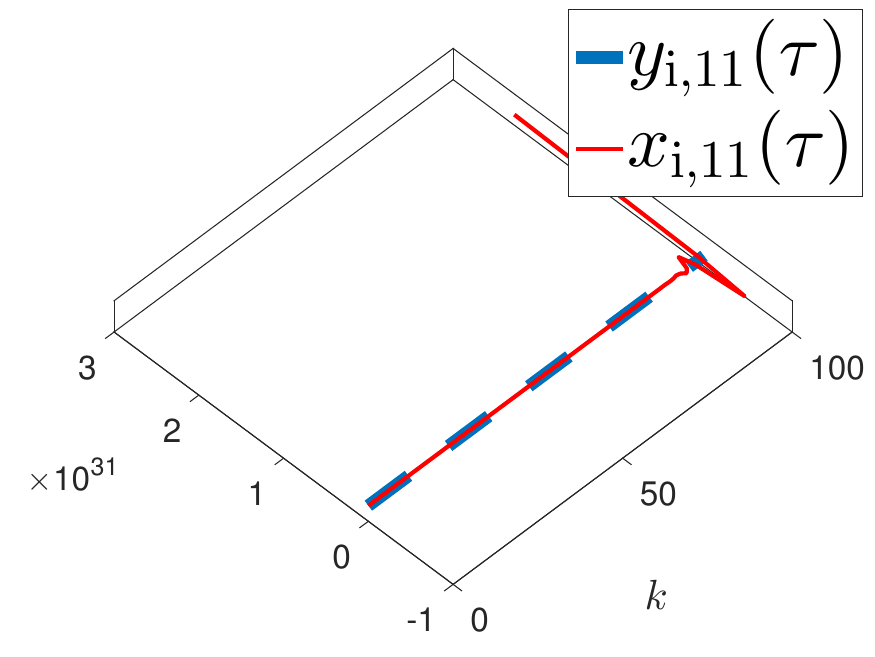}\label{fig.e1.Con-DZND1-2i.solve.10+20i.0.1.x11i.3d}}
	\subfigure[]{\includegraphics[width=0.40\columnwidth]{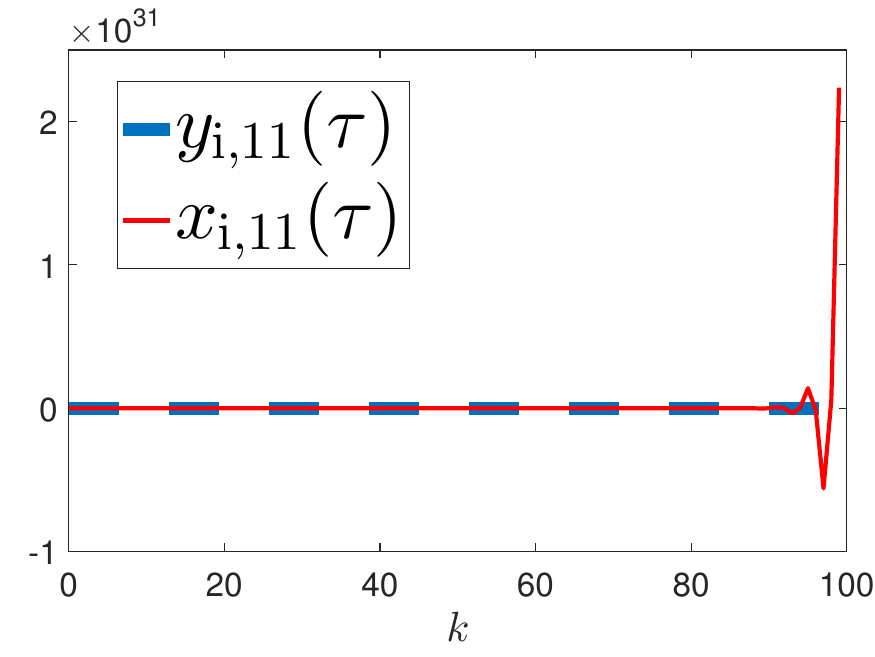}\label{fig.e1.Con-DZND1-2i.solve.10+20i.0.1.x11i.2d}}
	\subfigure[]{\includegraphics[width=0.70\columnwidth]{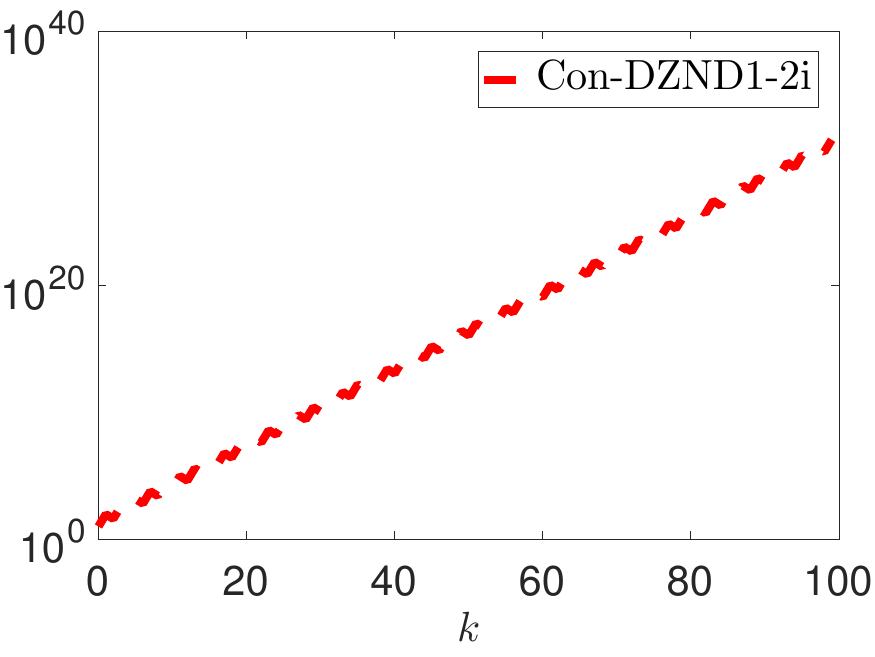}\label{fig.e1.normerror.10+20i.0.1}}
	\caption{Solution $X(\tau)$ computed by Con-DZND1-2i \eqref{eq.euler.forward.solve.linearerrconcznd1} model in Example \ref{example1} where $\gamma$ equals 10$+$20$\mathrm{i}$ and $\varepsilon$ equals 0.1. \subref{fig.e1.Con-DZND1-2i.solve.10+20i.0.1.x11r.2d} is the 2D view of \subref{fig.e1.Con-DZND1-2i.solve.10+20i.0.1.x11r.3d}, \subref{fig.e1.Con-DZND1-2i.solve.10+20i.0.1.x11i.2d} is the 2D view of \subref{fig.e1.Con-DZND1-2i.solve.10+20i.0.1.x11i.3d}, and \subref{fig.e1.normerror.10+20i.0.1} is logarithmic residual $\left \|X(\tau)-X^*(\tau)   \right \|_{\mathrm{F}}$ trajectory.}
	\label{fig.e1.Con-DZND1-2i.solve.10+20i.0.1}
\end{figure}

\begin{figure}[!h]\centering
	\subfigure[]{\includegraphics[width=0.40\columnwidth]{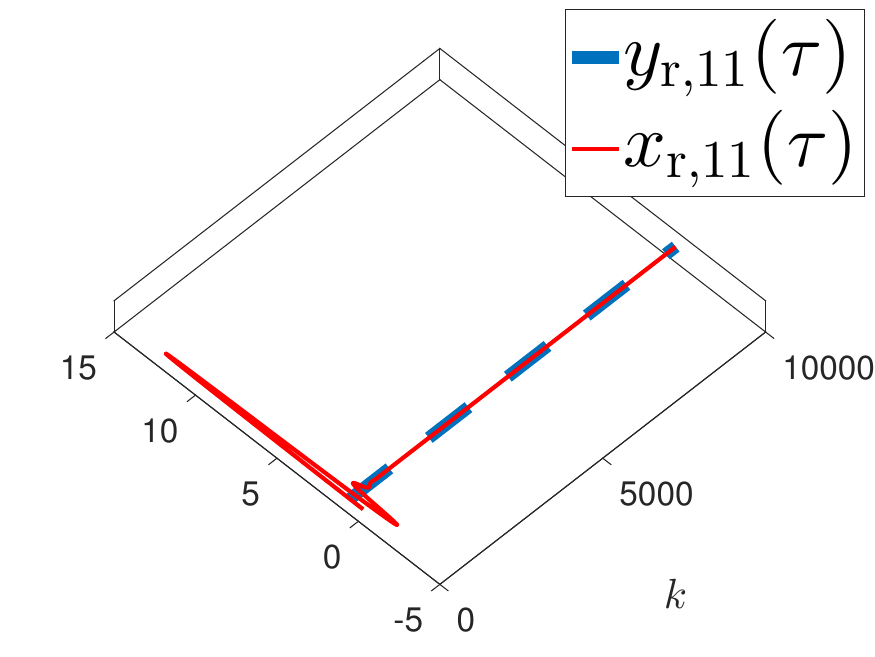}\label{fig.e1.Con-DZND1-2i.solve.10+20i.0.001.x11r.3d}}
	\subfigure[]{\includegraphics[width=0.40\columnwidth]{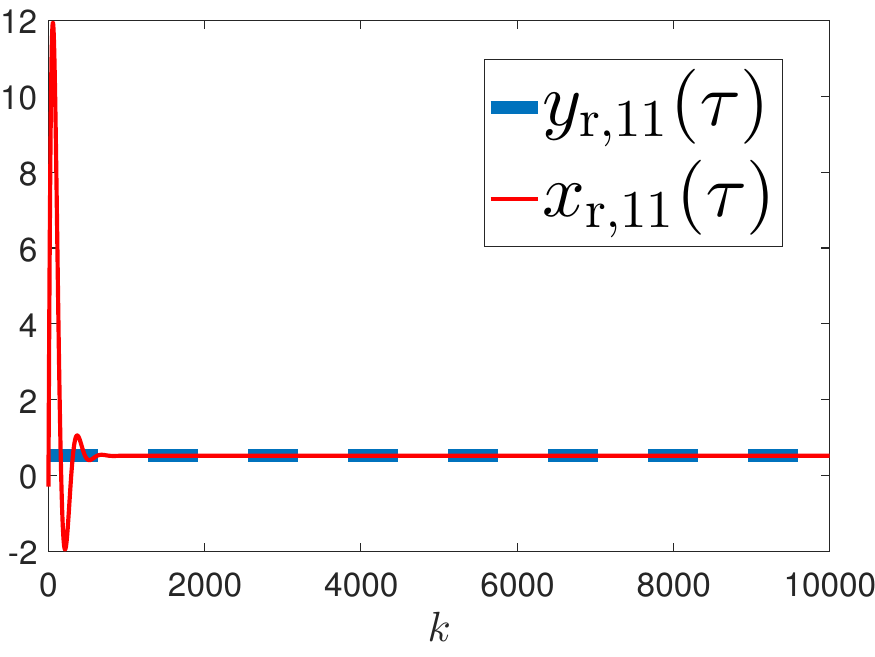}\label{fig.e1.Con-DZND1-2i.solve.10+20i.0.001.x11r.2d}}
	\subfigure[]{\includegraphics[width=0.40\columnwidth]{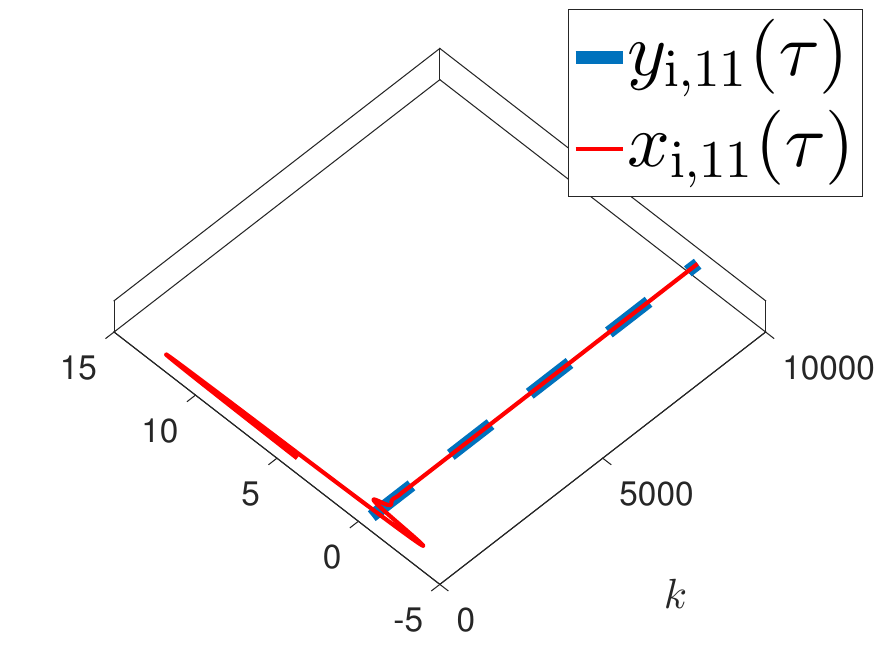}\label{fig.e1.Con-DZND1-2i.solve.10+20i.0.001.x11i.3d}}
	\subfigure[]{\includegraphics[width=0.40\columnwidth]{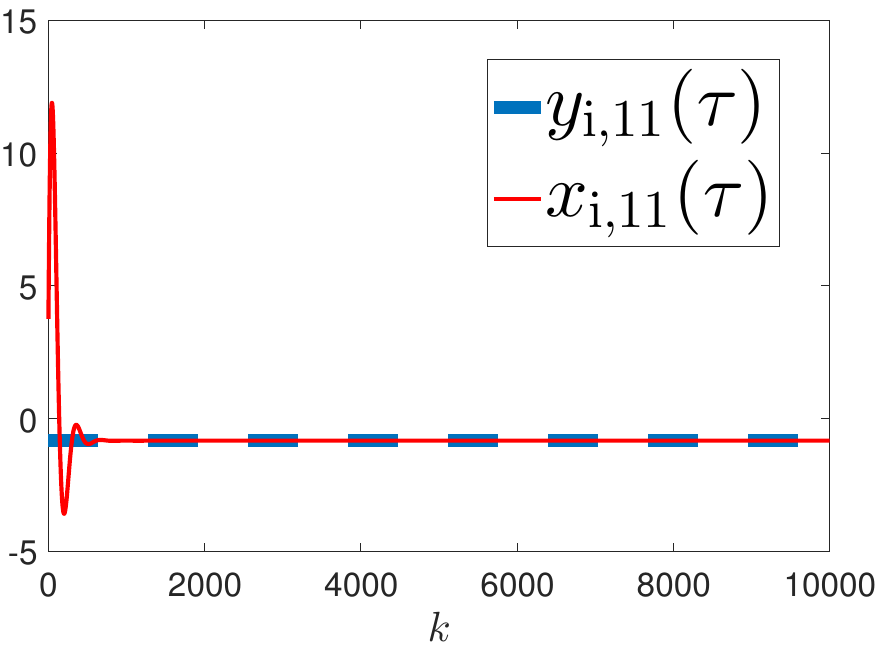}\label{fig.e1.Con-DZND1-2i.solve.10+20i.0.001.x11i.2d}}
	\subfigure[]{\includegraphics[width=0.70\columnwidth]{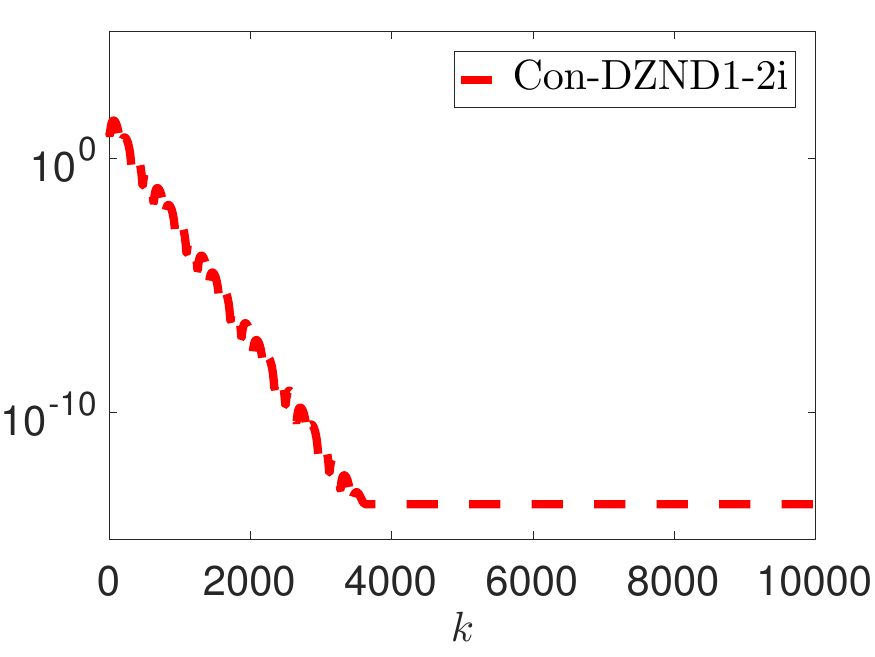}\label{fig.e1.normerror.10+20i.0.001}}
	\caption{Solution $X(\tau)$ computed by Con-DZND1-2i \eqref{eq.euler.forward.solve.linearerrconcznd1} model in Example \ref{example1} where $\gamma$ equals 10$+$20$\mathrm{i}$ and $\varepsilon$ equals 0.001.
	\subref{fig.e1.Con-DZND1-2i.solve.10+20i.0.001.x11r.2d} is the 2D view of \subref{fig.e1.Con-DZND1-2i.solve.10+20i.0.001.x11r.3d}, \subref{fig.e1.Con-DZND1-2i.solve.10+20i.0.001.x11i.2d} is the 2D view of \subref{fig.e1.Con-DZND1-2i.solve.10+20i.0.001.x11i.3d}, and \subref{fig.e1.normerror.10+20i.0.001} is logarithmic residual $\left \|X(\tau)-X^*(\tau)   \right \|_{\mathrm{F}}$ trajectory.}
	\label{fig.e1.Con-DZND1-2i.solve.10+20i.0.001}
\end{figure}

\begin{figure}[!h]\centering
	\subfigure[]{\includegraphics[width=0.40\columnwidth]{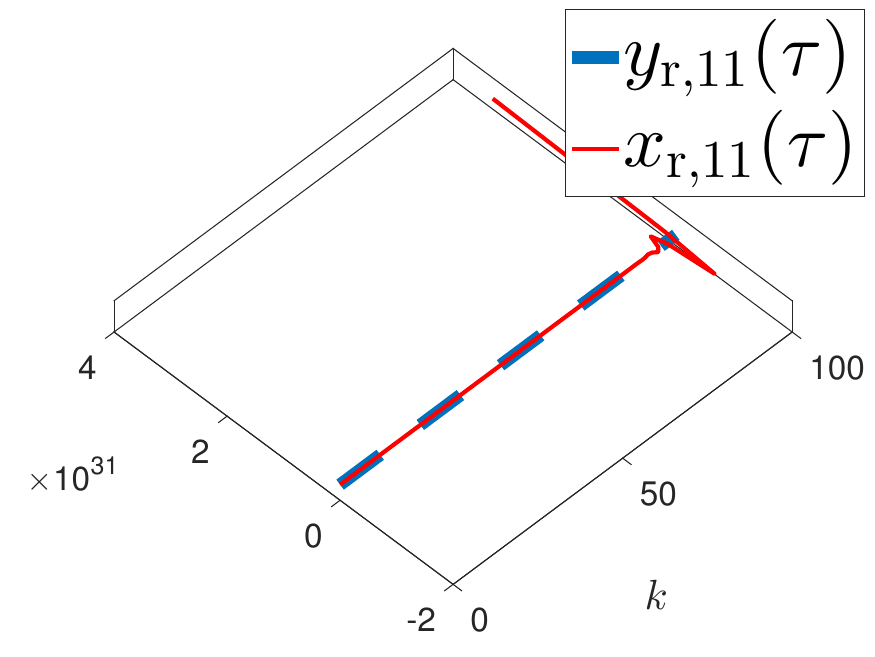}\label{fig.e1.Con-DZND1-2i.solve.10-20i.0.1.x11r.3d}}
	\subfigure[]{\includegraphics[width=0.40\columnwidth]{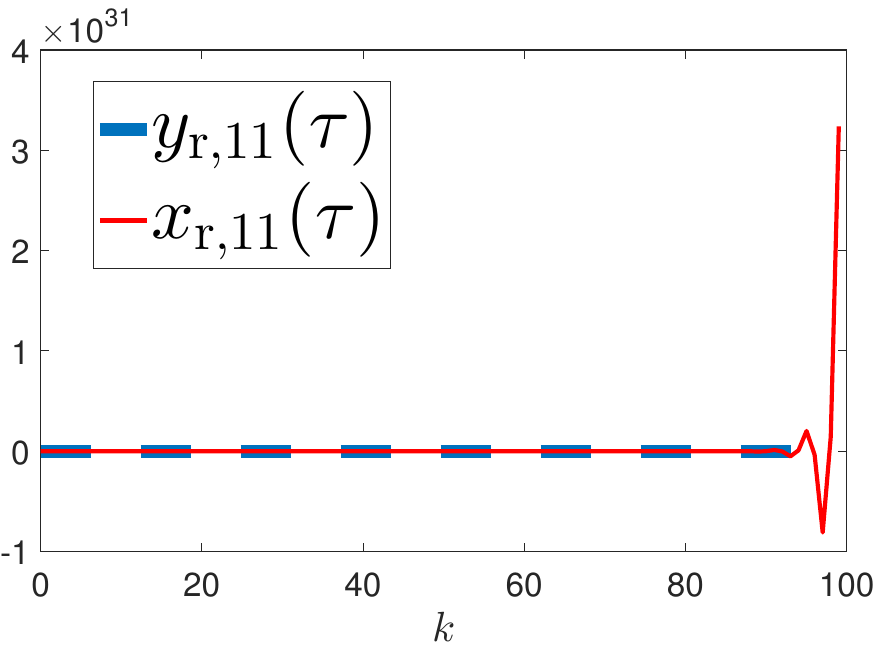}\label{fig.e1.Con-DZND1-2i.solve.10-20i.0.1.x11r.2d}}
	\subfigure[]{\includegraphics[width=0.40\columnwidth]{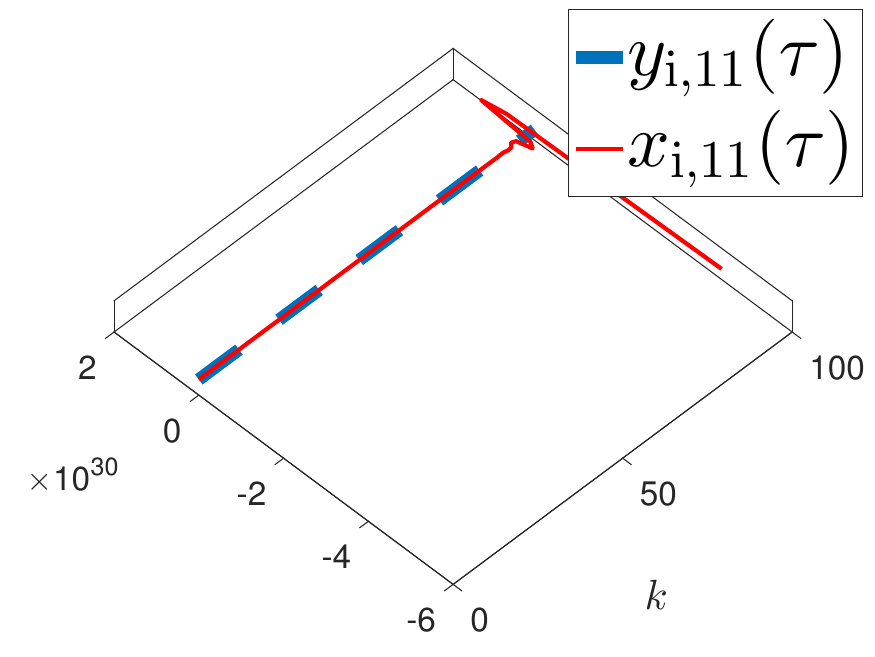}\label{fig.e1.Con-DZND1-2i.solve.10-20i.0.1.x11i.3d}}
	\subfigure[]{\includegraphics[width=0.40\columnwidth]{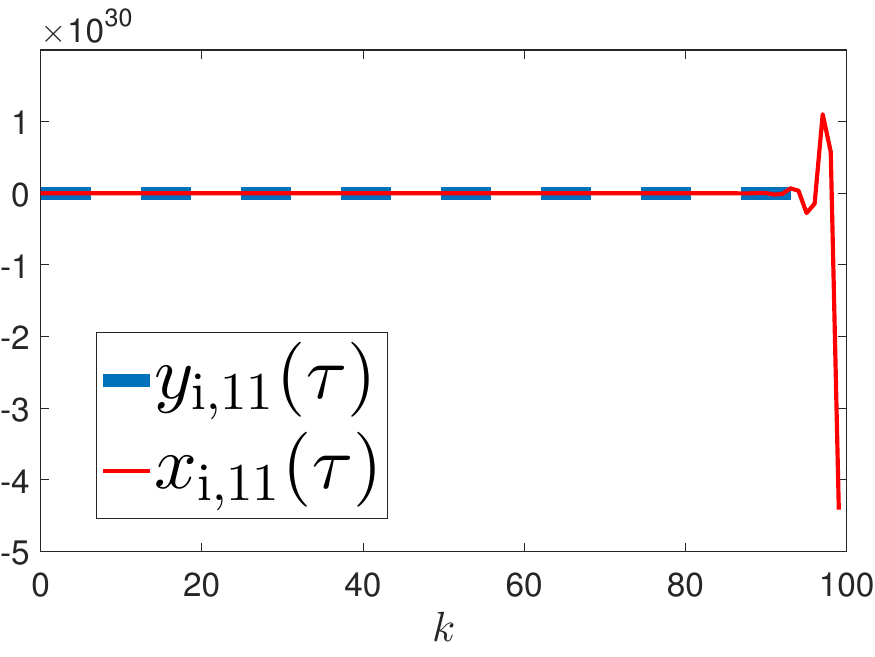}\label{fig.e1.Con-DZND1-2i.solve.10-20i.0.1.x11i.2d}}
	\subfigure[]{\includegraphics[width=0.70\columnwidth]{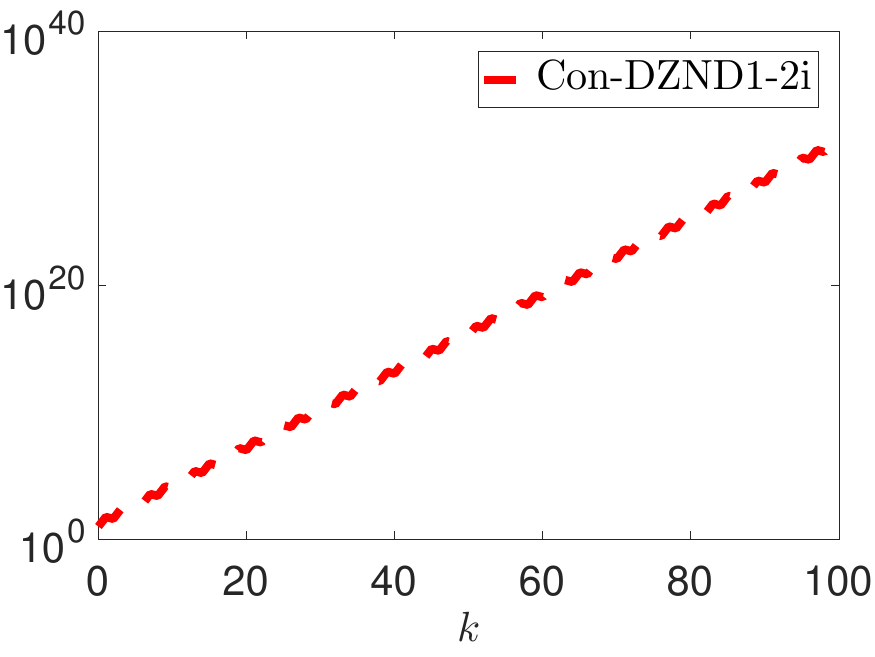}\label{fig.e1.normerror.10-20i.0.1}}
	\caption{Solution $X(\tau)$ computed by Con-DZND1-2i \eqref{eq.euler.forward.solve.linearerrconcznd1} model in Example \ref{example1} where $\gamma$ equals 10$-$20$\mathrm{i}$ and $\varepsilon$ equals 0.1.
	\subref{fig.e1.Con-DZND1-2i.solve.10-20i.0.1.x11r.2d} is the 2D view of \subref{fig.e1.Con-DZND1-2i.solve.10-20i.0.1.x11r.3d}, \subref{fig.e1.Con-DZND1-2i.solve.10-20i.0.1.x11i.2d} is the 2D view of \subref{fig.e1.Con-DZND1-2i.solve.10-20i.0.1.x11i.3d}, and \subref{fig.e1.normerror.10-20i.0.1} is logarithmic residual $\left \|X(\tau)-X^*(\tau)   \right \|_{\mathrm{F}}$ trajectory.}
	\label{fig.e1.Con-DZND1-2i.solve.10-20i.0.1}
\end{figure}

\begin{figure}[!h]\centering
	\subfigure[]{\includegraphics[width=0.40\columnwidth]{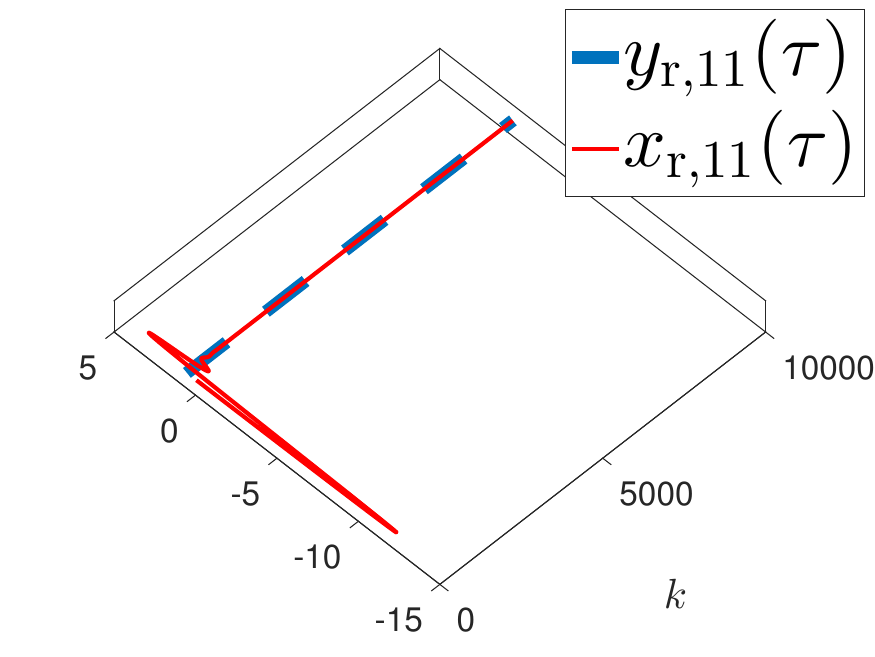}\label{fig.e1.Con-DZND1-2i.solve.10-20i.0.001.x11r.3d}}
	\subfigure[]{\includegraphics[width=0.40\columnwidth]{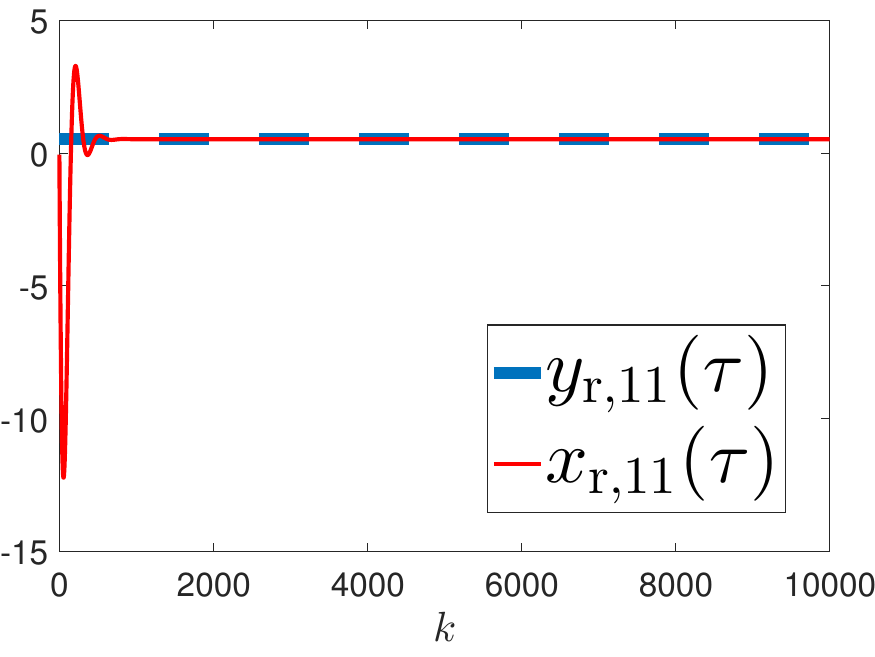}\label{fig.e1.Con-DZND1-2i.solve.10-20i.0.001.x11r.2d}}
	\subfigure[]{\includegraphics[width=0.40\columnwidth]{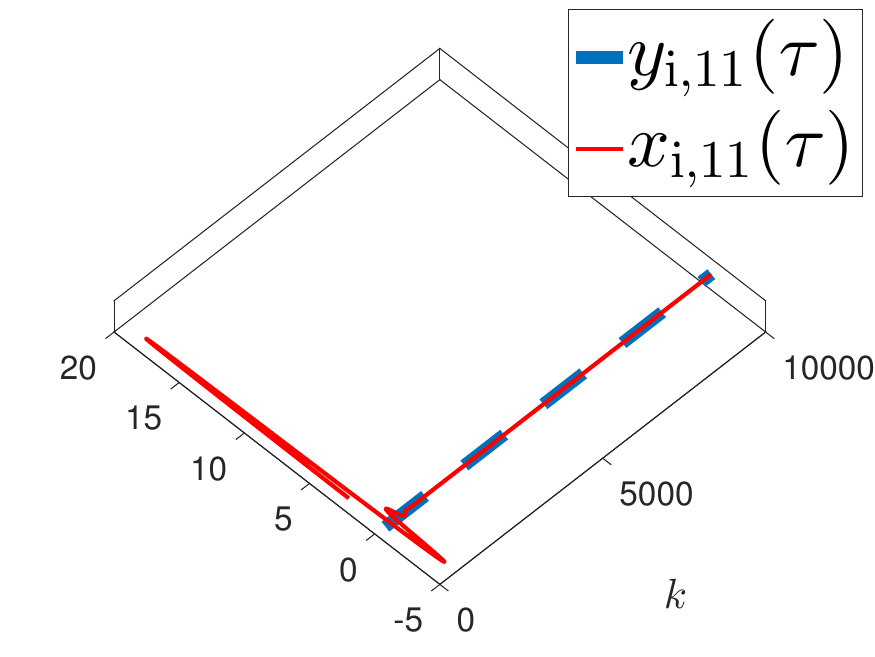}\label{fig.e1.Con-DZND1-2i.solve.10-20i.0.001.x11i.3d}}
	\subfigure[]{\includegraphics[width=0.40\columnwidth]{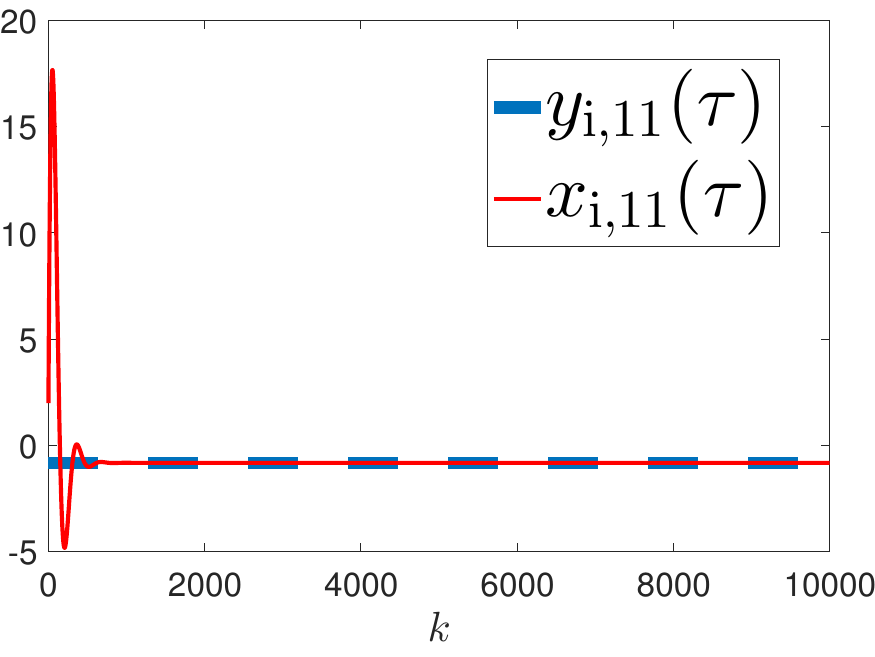}\label{fig.e1.Con-DZND1-2i.solve.10-20i.0.001.x11i.2d}}
	\subfigure[]{\includegraphics[width=0.70\columnwidth]{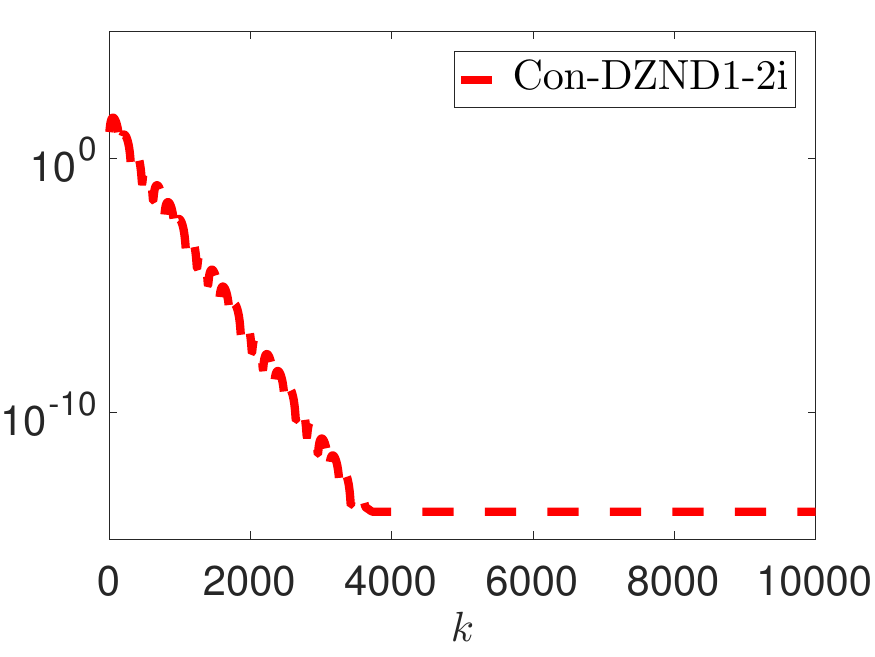}\label{fig.e1.normerror.10-20i.0.001}}
	\caption{Solution $X(\tau)$ computed by Con-DZND1-2i \eqref{eq.euler.forward.solve.linearerrconcznd1} model in Example \ref{example1} where $\gamma$ equals 10$-$20$\mathrm{i}$ and $\varepsilon$ equals 0.001.
	\subref{fig.e1.Con-DZND1-2i.solve.10-20i.0.001.x11r.2d} is the 2D view of \subref{fig.e1.Con-DZND1-2i.solve.10-20i.0.001.x11r.3d}, \subref{fig.e1.Con-DZND1-2i.solve.10-20i.0.001.x11i.2d} is the 2D view of \subref{fig.e1.Con-DZND1-2i.solve.10-20i.0.001.x11i.3d}, and \subref{fig.e1.normerror.10-20i.0.001} is logarithmic residual $\left \|X(\tau)-X^*(\tau)   \right \|_{\mathrm{F}}$ trajectory.}
	\label{fig.e1.Con-DZND1-2i.solve.10-20i.0.001}
\end{figure}

\begin{figure}[!h]\centering
	\subfigure[]{\includegraphics[width=0.40\columnwidth]{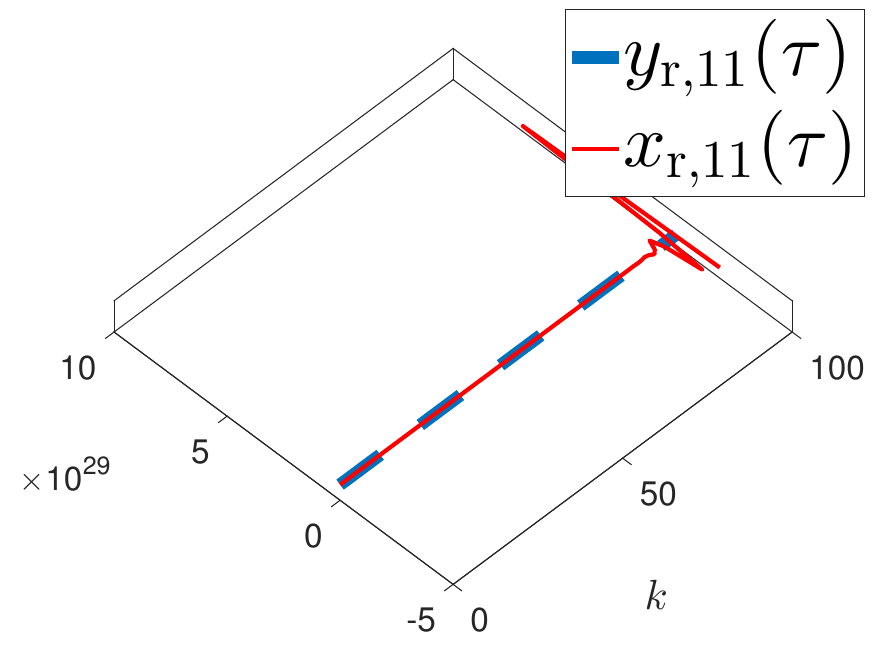}\label{fig.e2.Con-DZND1-2i.solve.10+20i.0.1.x11r.3d}}
	\subfigure[]{\includegraphics[width=0.40\columnwidth]{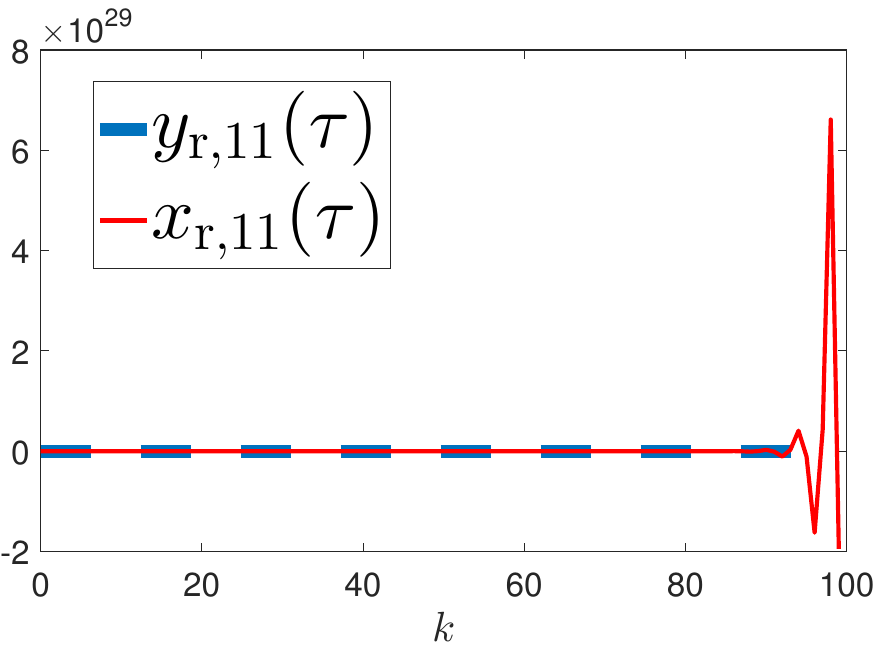}\label{fig.e2.Con-DZND1-2i.solve.10+20i.0.1.x11r.2d}}
	\subfigure[]{\includegraphics[width=0.40\columnwidth]{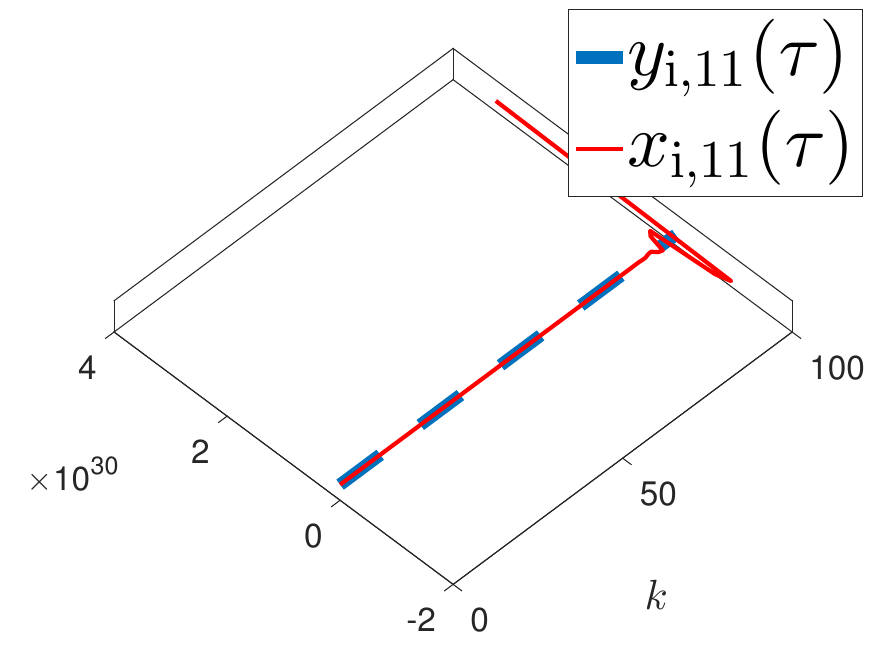}\label{fig.e2.Con-DZND1-2i.solve.10+20i.0.1.x11i.3d}}
	\subfigure[]{\includegraphics[width=0.40\columnwidth]{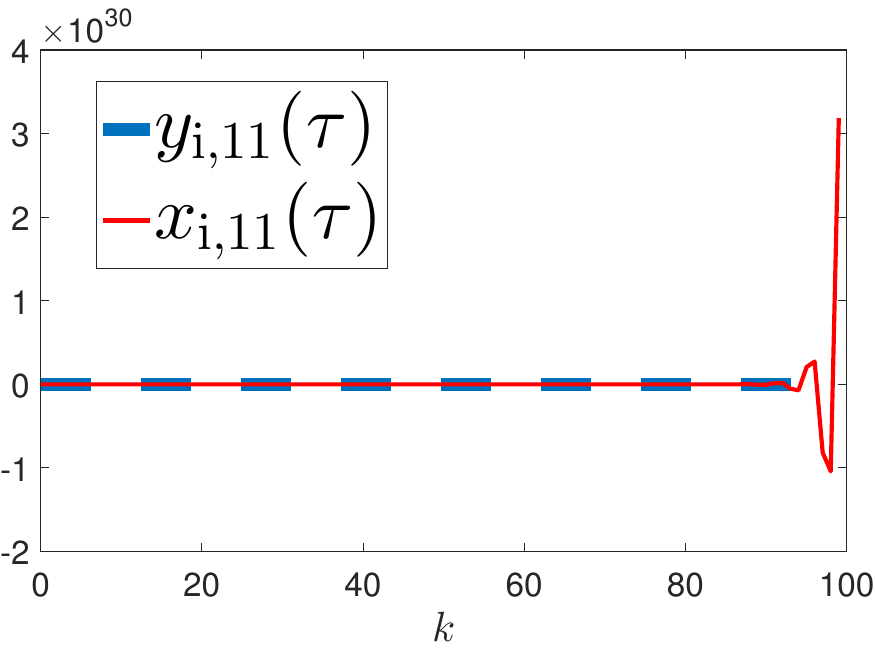}\label{fig.e2.Con-DZND1-2i.solve.10+20i.0.1.x11i.2d}}
	\subfigure[]{\includegraphics[width=0.70\columnwidth]{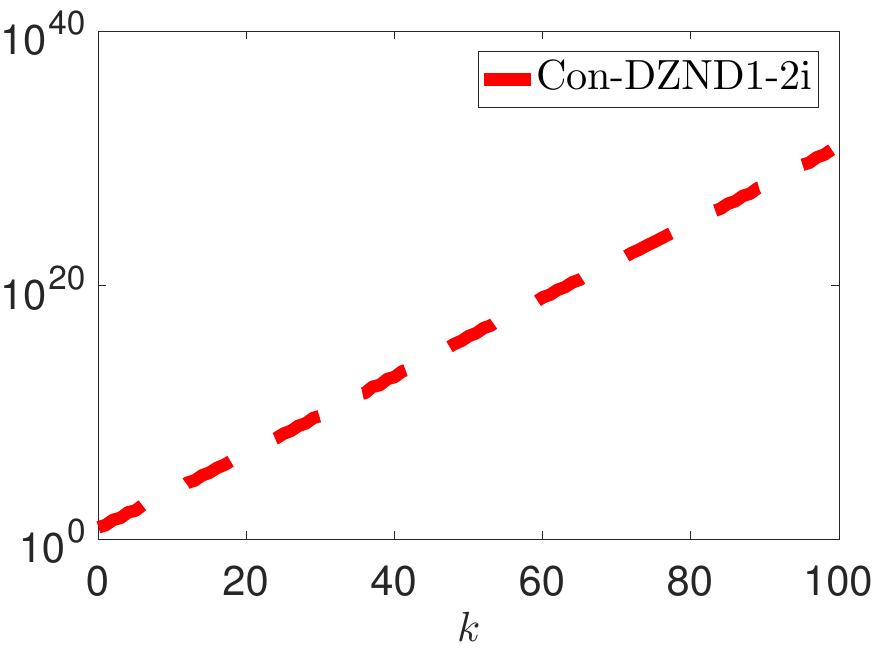}\label{fig.e2.normerror.10+20i.0.1}}
	\caption{Solution $X(\tau)$ computed by Con-DZND1-2i \eqref{eq.euler.forward.solve.linearerrconcznd1} model in Example \ref{example2} where $\gamma$ equals 10$+$20$\mathrm{i}$ and $\varepsilon$ equals 0.1.
	\subref{fig.e2.Con-DZND1-2i.solve.10+20i.0.1.x11r.2d} is the 2D view of \subref{fig.e2.Con-DZND1-2i.solve.10+20i.0.1.x11r.3d}, \subref{fig.e2.Con-DZND1-2i.solve.10+20i.0.1.x11i.2d} is the 2D view of \subref{fig.e2.Con-DZND1-2i.solve.10+20i.0.1.x11i.3d}, and \subref{fig.e2.normerror.10+20i.0.1} is logarithmic residual $\left \|X(\tau)-X^*(\tau)   \right \|_{\mathrm{F}}$ trajectory.}
	\label{fig.e2.Con-DZND1-2i.solve.10+20i.0.1}
\end{figure}

\begin{figure}[!h]\centering
	\subfigure[]{\includegraphics[width=0.40\columnwidth]{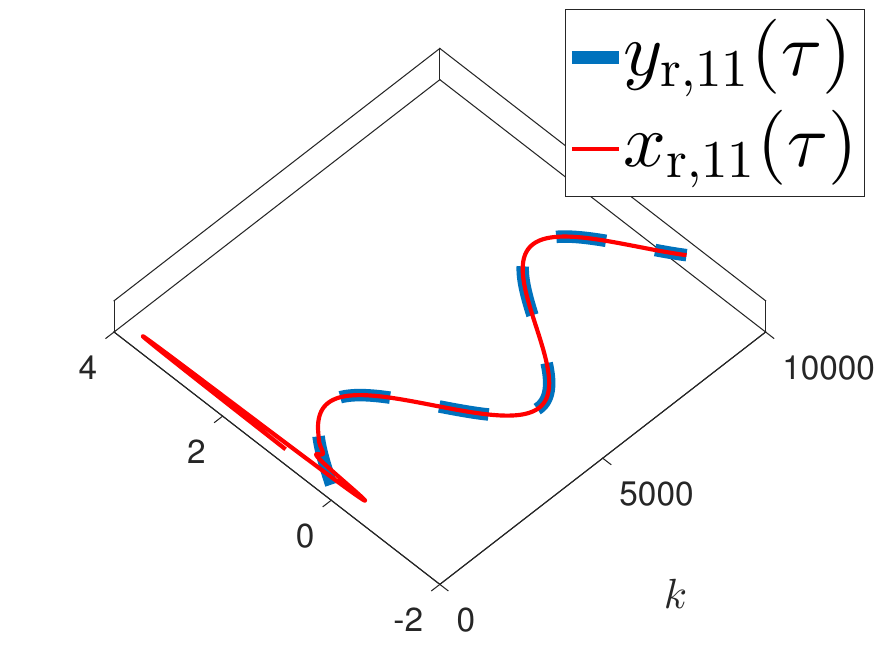}\label{fig.e2.Con-DZND1-2i.solve.10+20i.0.001.x11r.3d}}
	\subfigure[]{\includegraphics[width=0.40\columnwidth]{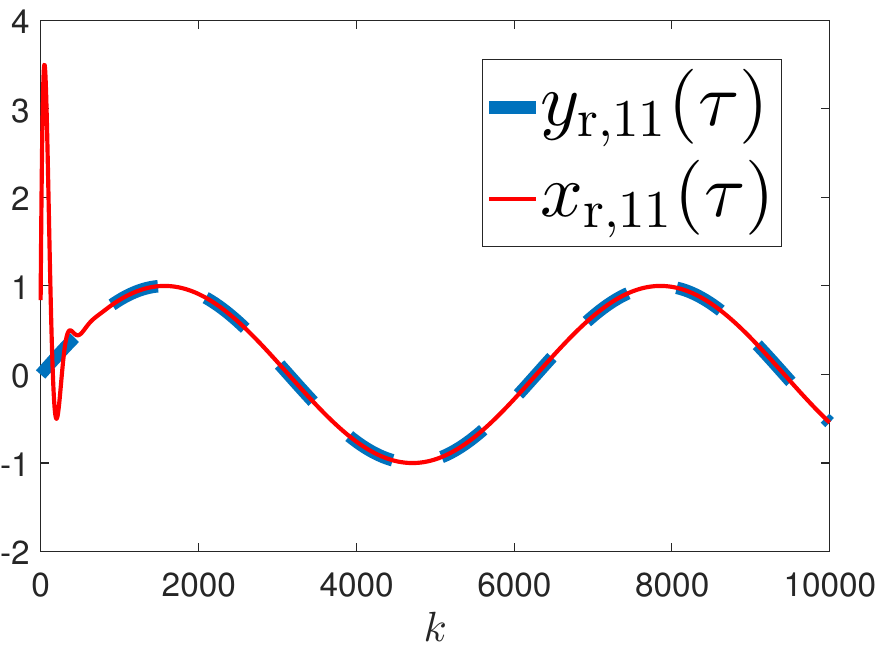}\label{fig.e2.Con-DZND1-2i.solve.10+20i.0.001.x11r.2d}}
	\subfigure[]{\includegraphics[width=0.40\columnwidth]{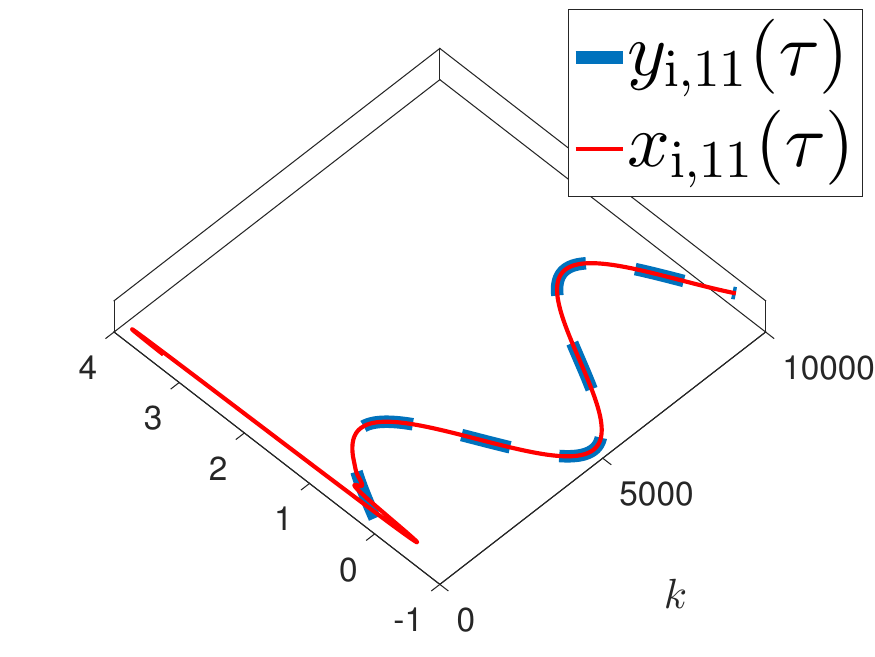}\label{fig.e2.Con-DZND1-2i.solve.10+20i.0.001.x11i.3d}}
	\subfigure[]{\includegraphics[width=0.40\columnwidth]{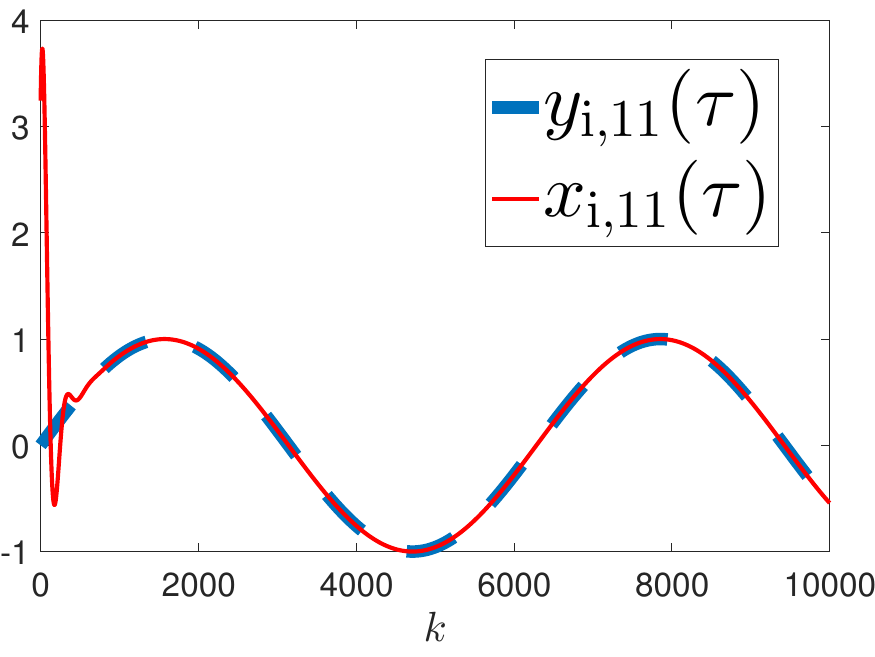}\label{fig.e2.Con-DZND1-2i.solve.10+20i.0.001.x11i.2d}}
	\subfigure[]{\includegraphics[width=0.70\columnwidth]{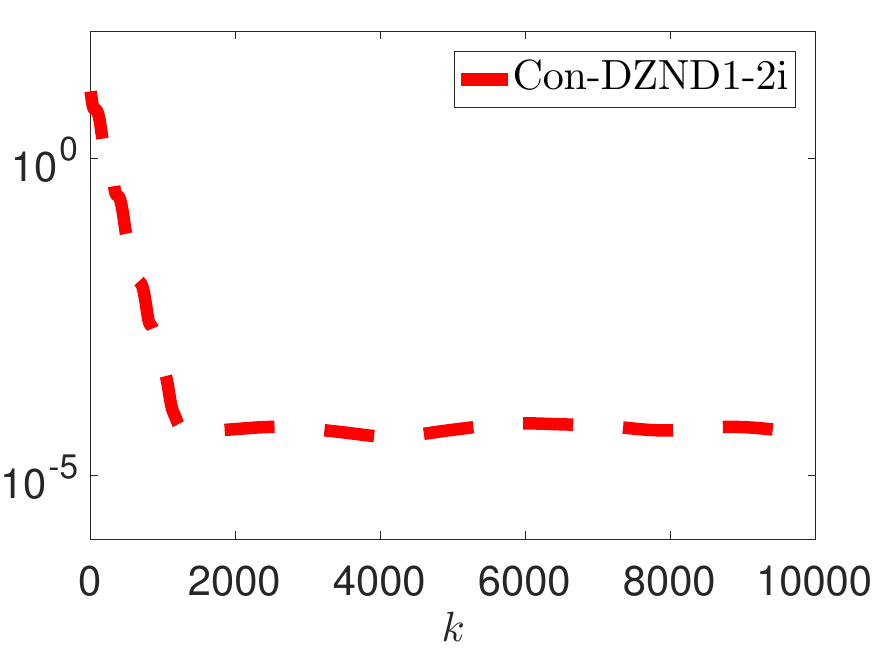}\label{fig.e2.normerror.10+20i.0.001}}
	\caption{Solution $X(\tau)$ computed by Con-DZND1-2i \eqref{eq.euler.forward.solve.linearerrconcznd1} model in Example \ref{example2} where $\gamma$ equals 10$+$20$\mathrm{i}$ and $\varepsilon$ equals 0.001.
	\subref{fig.e2.Con-DZND1-2i.solve.10+20i.0.001.x11r.2d} is the 2D view of \subref{fig.e2.Con-DZND1-2i.solve.10+20i.0.001.x11r.3d}, \subref{fig.e2.Con-DZND1-2i.solve.10+20i.0.001.x11i.2d} is the 2D view of \subref{fig.e2.Con-DZND1-2i.solve.10+20i.0.001.x11i.3d}, and \subref{fig.e2.normerror.10+20i.0.001} is logarithmic residual $\left \|X(\tau)-X^*(\tau)   \right \|_{\mathrm{F}}$ trajectory.}
	\label{fig.e2.Con-DZND1-2i.solve.10+20i.0.001}
\end{figure}

\begin{figure}[!h]\centering
	\subfigure[]{\includegraphics[width=0.40\columnwidth]{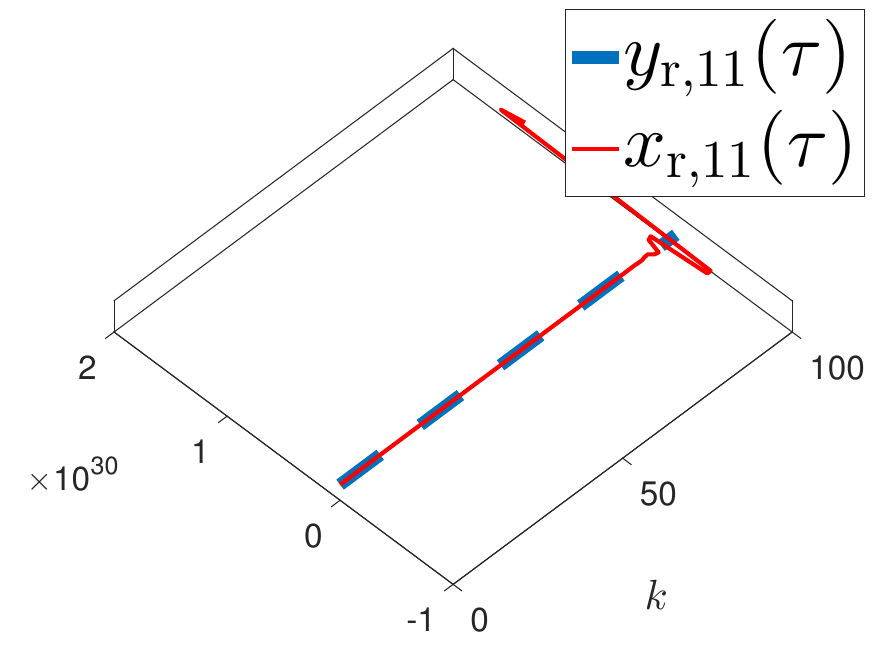}\label{fig.e2.Con-DZND1-2i.solve.10-20i.0.1.x11r.3d}}
	\subfigure[]{\includegraphics[width=0.40\columnwidth]{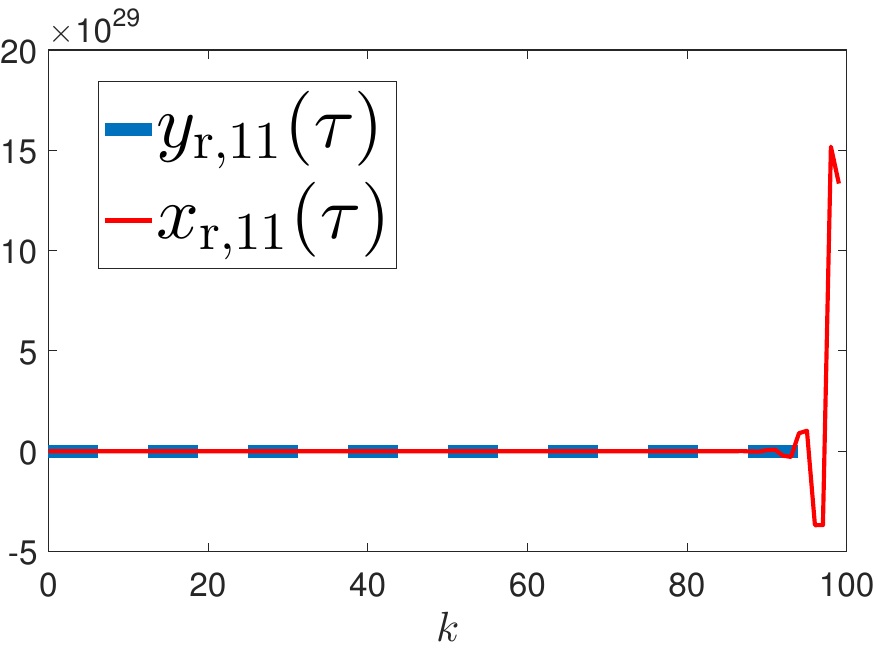}\label{fig.e2.Con-DZND1-2i.solve.10-20i.0.1.x11r.2d}}
	\subfigure[]{\includegraphics[width=0.40\columnwidth]{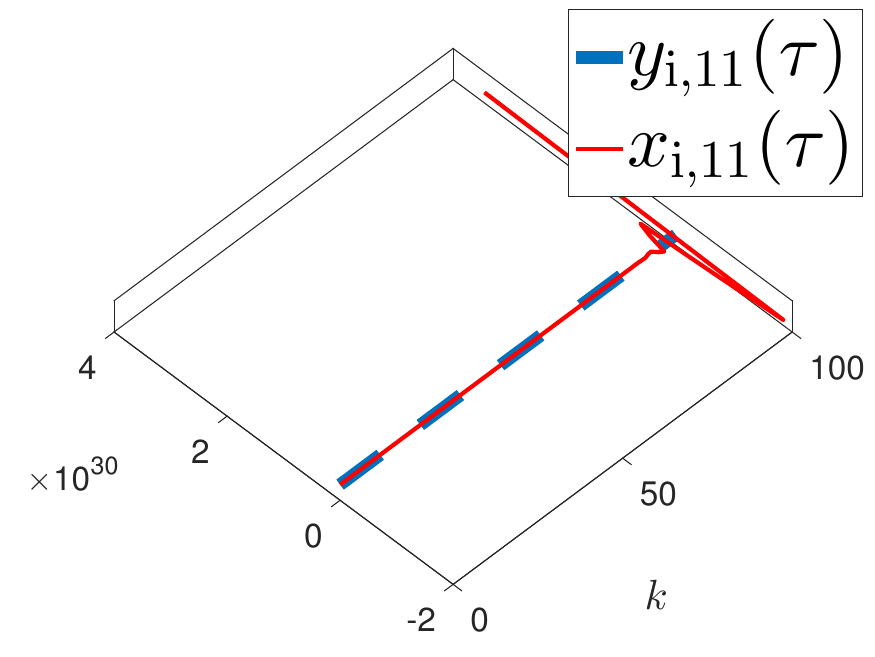}\label{fig.e2.Con-DZND1-2i.solve.10-20i.0.1.x11i.3d}}
	\subfigure[]{\includegraphics[width=0.40\columnwidth]{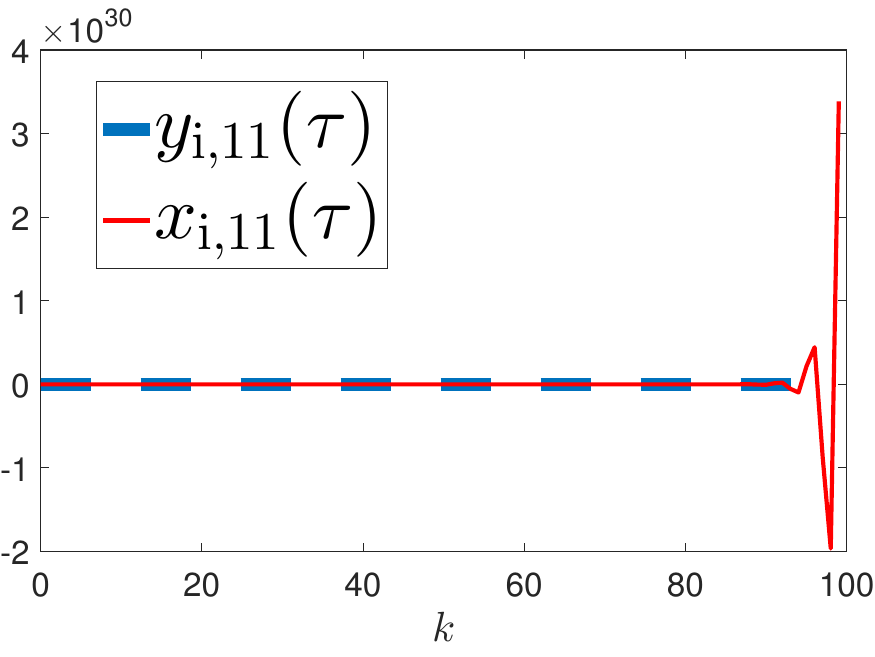}\label{fig.e2.Con-DZND1-2i.solve.10-20i.0.1.x11i.2d}}
	\subfigure[]{\includegraphics[width=0.70\columnwidth]{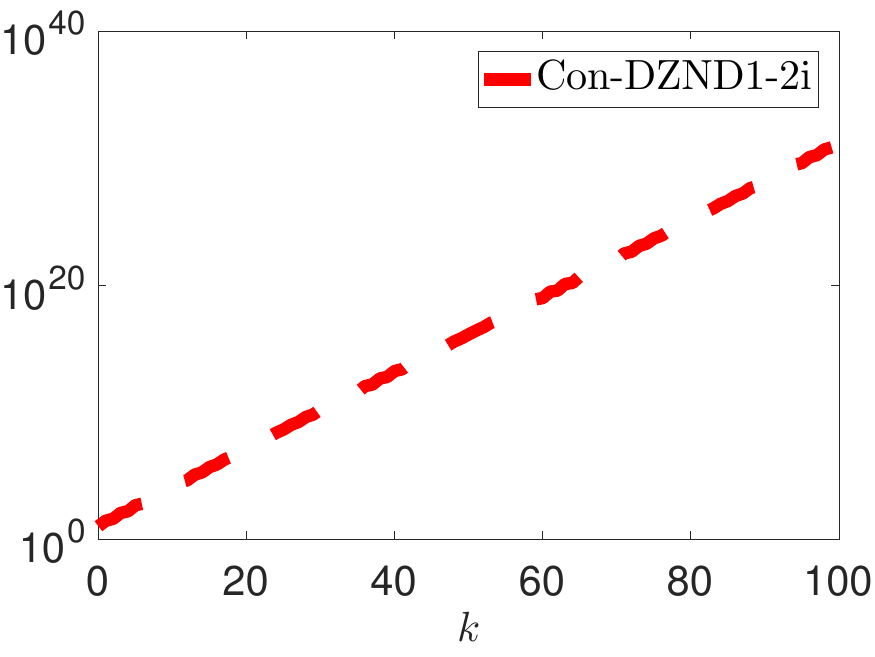}\label{fig.e2.normerror.10-20i.0.1}}
	\caption{Solution $X(\tau)$ computed by Con-DZND1-2i \eqref{eq.euler.forward.solve.linearerrconcznd1} model in Example \ref{example2} where $\gamma$ equals 10$-$20$\mathrm{i}$ and $\varepsilon$ equals 0.1.
	\subref{fig.e2.Con-DZND1-2i.solve.10-20i.0.1.x11r.2d} is the 2D view of \subref{fig.e2.Con-DZND1-2i.solve.10-20i.0.1.x11r.3d}, \subref{fig.e2.Con-DZND1-2i.solve.10-20i.0.1.x11i.2d} is the 2D view of \subref{fig.e2.Con-DZND1-2i.solve.10-20i.0.1.x11i.3d}, and \subref{fig.e2.normerror.10-20i.0.1} is logarithmic residual $\left \|X(\tau)-X^*(\tau)   \right \|_{\mathrm{F}}$ trajectory.}
	\label{fig.e2.Con-DZND1-2i.solve.10-20i.0.1}
\end{figure}

\begin{figure}[!h]\centering
	\subfigure[]{\includegraphics[width=0.40\columnwidth]{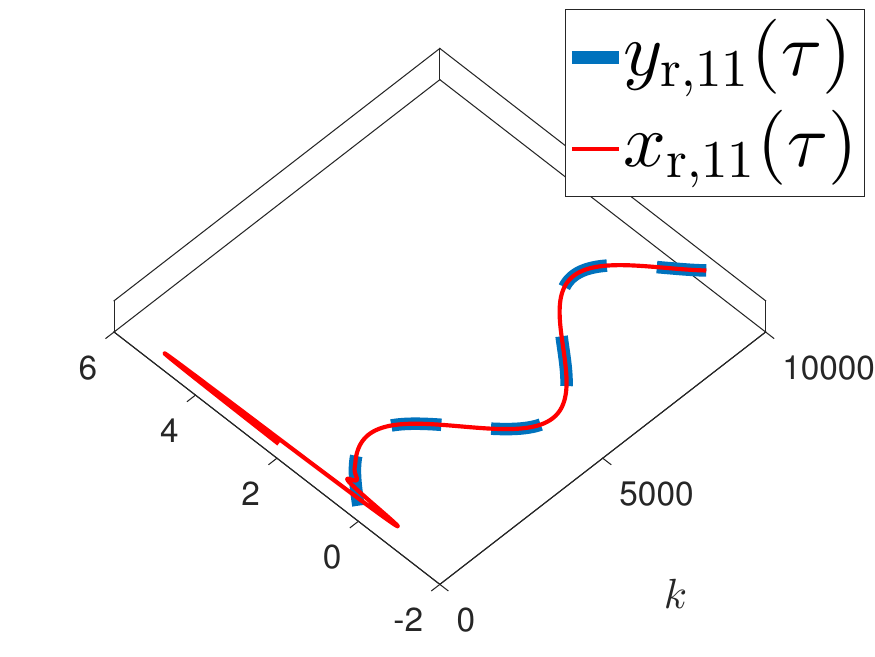}\label{fig.e2.Con-DZND1-2i.solve.10-20i.0.001.x11r.3d}}
	\subfigure[]{\includegraphics[width=0.40\columnwidth]{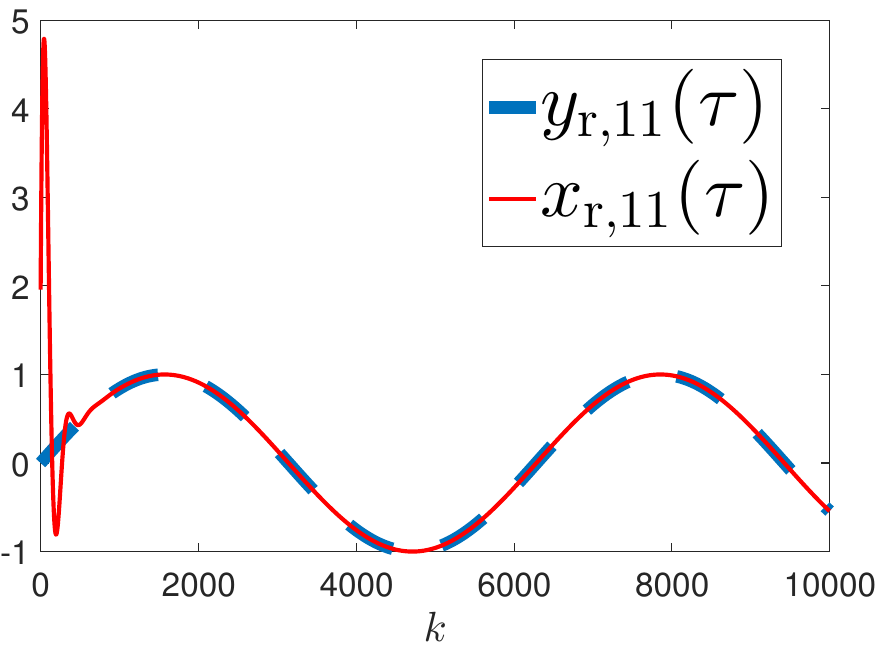}\label{fig.e2.Con-DZND1-2i.solve.10-20i.0.001.x11r.2d}}
	\subfigure[]{\includegraphics[width=0.40\columnwidth]{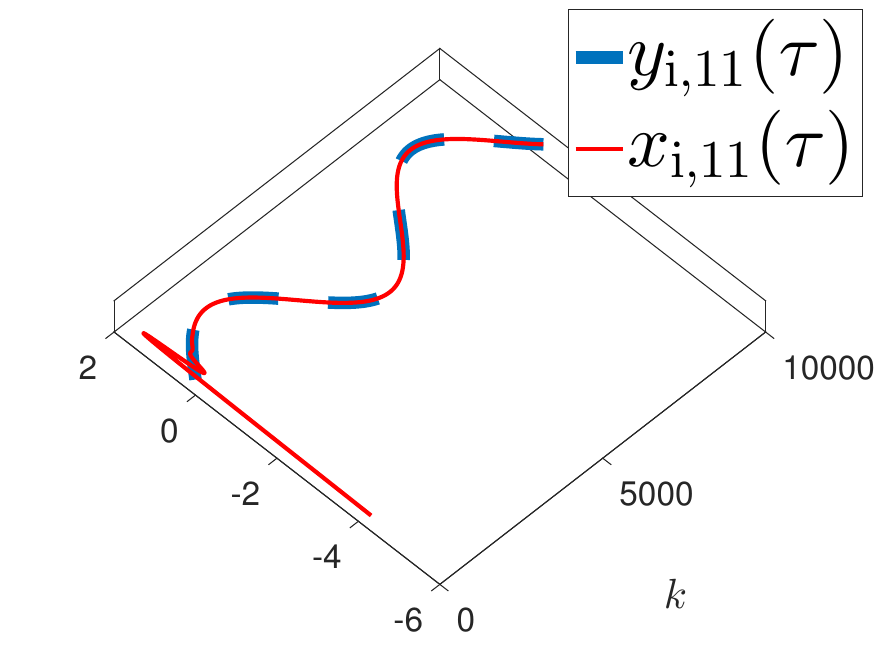}\label{fig.e2.Con-DZND1-2i.solve.10-20i.0.001.x11i.3d}}
	\subfigure[]{\includegraphics[width=0.40\columnwidth]{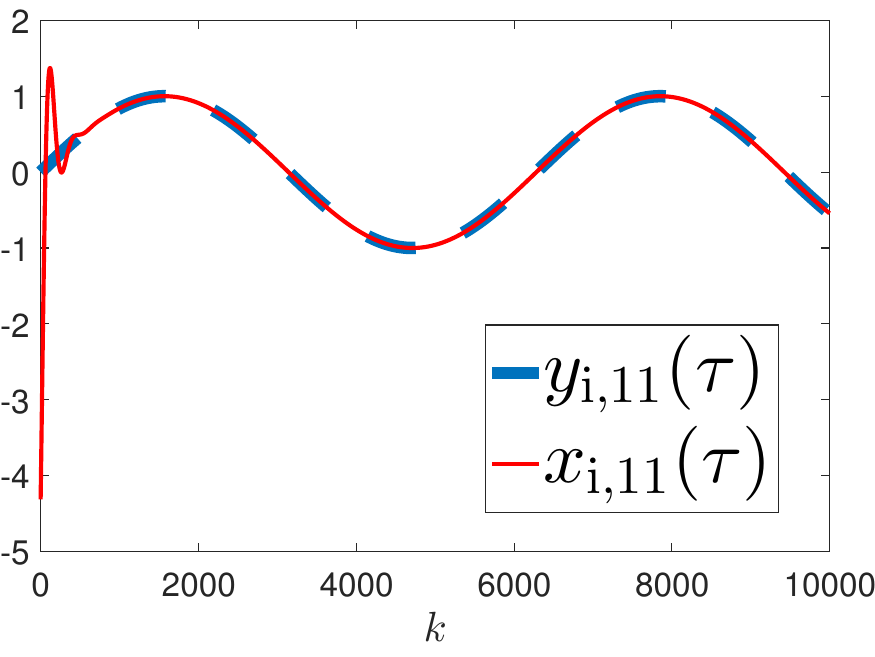}\label{fig.e2.Con-DZND1-2i.solve.10-20i.0.001.x11i.2d}}
	\subfigure[]{\includegraphics[width=0.70\columnwidth]{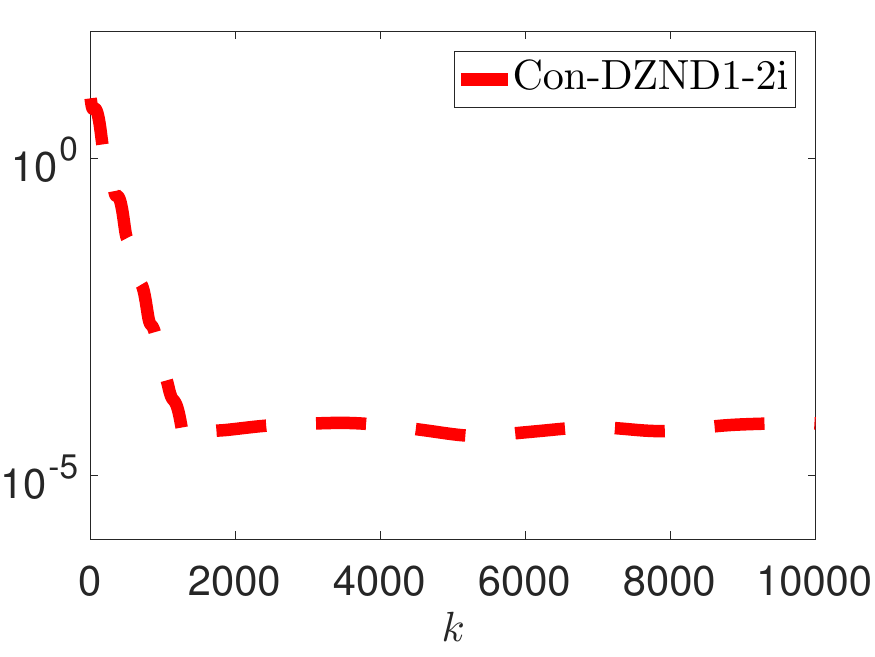}\label{fig.e2.normerror.10-20i.0.001}}
	\caption{Solution $X(\tau)$ computed by Con-DZND1-2i \eqref{eq.euler.forward.solve.linearerrconcznd1} model in Example \ref{example1} where $\gamma$ equals 10$-$20$\mathrm{i}$ and $\varepsilon$ equals 0.001.
	\subref{fig.e2.Con-DZND1-2i.solve.10-20i.0.001.x11r.2d} is the 2D view of \subref{fig.e2.Con-DZND1-2i.solve.10-20i.0.001.x11r.3d}, \subref{fig.e2.Con-DZND1-2i.solve.10-20i.0.001.x11i.2d} is the 2D view of \subref{fig.e2.Con-DZND1-2i.solve.10-20i.0.001.x11i.3d}, and \subref{fig.e2.normerror.10-20i.0.001} is logarithmic residual $\left \|X(\tau)-X^*(\tau)   \right \|_{\mathrm{F}}$ trajectory.}
	\label{fig.e2.Con-DZND1-2i.solve.10-20i.0.001}
\end{figure}

For 
$\gamma$ equals 10$+$20$\mathrm{i}$ and 10$-$20$\mathrm{i}$ respectively, to highlight the distinction between space compressive approximation errors and sampling discretion errors, this article only presents the real and imaginary parts of the solution for $x_{r,11}(\tau)$ and $x_{i,11}(\tau)$ in Examples \ref{example1} and \ref{example2} under different step sizes $\varepsilon$, along with the model's logarithmic residual $\left \|X(\tau)-X^*(\tau)   \right \|_{\mathrm{F}}$ trajectories are shown in Figs. \ref{fig.e1.Con-DZND1-2i.solve.10+20i.0.1}
through
\ref{fig.e2.Con-DZND1-2i.solve.10-20i.0.001}.


It can be seen that at step size equals $\varepsilon$ 0.1, with $\gamma$ of 
10$+$20$\mathrm{i}$ and 
10$-$20$\mathrm{i}$, Con-DZND1-2i \eqref{eq.euler.forward.solve.linearerrconcznd1} model appears to diverge (model failure), like ``nan" values. However, at step size $\varepsilon$ equals 0.001, Con-DZND1-2i \eqref{eq.euler.forward.solve.linearerrconcznd1} model remains consistent with convergence at $\gamma$ equals 10. In another way, model with 0.1 step size $\varepsilon$ crashed, while model with 0.001 step size $\varepsilon$ remains stable.
\section{Discussion}
From the above experiments, it can be seen that the performance of neural networks is influenced by many factors. Based on previous studies
\cite{kimNeuralMachineCode2023,Evanusa2024MaelstromN,Kong2024ReservoircomputingBA,blanes2024splittingmethodscomplexcoefficients,an2024fastforwardingquantumalgorithmslinear},
some hypotheses can be proposed:
\begin{enumerate}[1.]
	\item At the regulation parameter controlling the convergence rate, stable numerical algorithms can eradicate high dimension phenomenon, where the optimization of sampling discretion errors predominates. In contrast, in variable-step size numerical algorithms, optimizing space compressive approximation is evident.
	
	\item Actually, the phenomenon of space compressive approximation in neural networks still exists, but in stable sampling numerical algorithms (with constant-step size), this space compressive approximation is hidden. At this point, it is necessary to test with additional space vectors, such as imaginary vectors of the regulation parameter controlling the convergence rate, as internal noise.
	
	\item Under the discretion of Euler-forward formula, neural hypercomplex numbers space compressive approximation approach (NHNSCAA)
	\cite{He2024RevisitingTC}
	serves only as an auxiliary method. NHNSCAA aids convergence but has a limit, which is also affected by step size. Comparing to ode45
	\cite{flynn2024exploringoriginsswitchingdynamics,dabounouAdaptiveFeedforwardGradient2024}
	variable-step size solver, embedding space vector noise necessitates adjusting the step size to achieve model convergence.
	
	\item The four neural dynamic models, Con-CZND1
	\cite{He2024ZeroingND,He2024RevisitingTC}, Con-DZND1-2i \eqref{eq.euler.forward.solve.linearerrconcznd1}, Con-CZND2
	\cite{He2024ZeroingND,He2024RevisitingTC}, and Con-DZND2-2i \eqref{eq.euler.forward.solve.linearerrconcznd2}, are fundamentally different and should be considered distinctly.
	
	\item NHNSCAA should not be abused (similar to synthesis data generation), as it may be vulnerable to additional space attacks by vector noise, leading to model failure. See Appendix \ref{appendix.C} for an explanation of ``No free lunch".
	
	\item The authors believe that the main obstacles to the development of AI lie in mathematics and materials science. In mathematics, the relevant topic is ``Number theory", while in materials science, the focus is on ``High-performance computing". The error matrix is analogous to ``Integrated circuit chip".
\end{enumerate}
\section{Conclusion}
This article uses Euler-forward formula to discrete continuous ZND models for solving TVSSCME \eqref{eq.sccsme.variant}. Then Con-DZND1-2i \eqref{eq.euler.forward.solve.linearerrconcznd1} model and Con-DZND2-2i \eqref{eq.euler.forward.solve.linearerrconcznd2} model are proposed. Additionally, the convergence, stability, and errors of the above models under different step sizes and various learning rates are analyzed. Comparing differences between optimizing sampling discretion errors and space compressive approximation errors reveals some key details. These insights relate to neural network black box. After, we will focus on three main directions:
\begin{itemize}
\item Conduct a comprehensive analysis of previous studies.
\item Supplement other knowledge of CCME that are relevant to this study.
\item Explore other details of neural network black box.
\end{itemize}
\section*{Declarations}
This work is aided by the Project Supported by the Guangzhou Science
and Technology Program (with number 2023E04J1240).
\begin{appendices}
	\section{Gradient descent}\label{appendix.A}
	When the function $f(x)$ needs to find its minimum called ``$\min f(x)$",
	``Gradient descent"
	\cite{10.5555/2987994}
	essentially updates the value at iteration 
	$k+1$ as:
	\begin{equation}\label{eq.gradientdescent.method}
		x^{(k+1)} \gets x^{(k)} +\varepsilon_{k}p_{k},
	\end{equation}
	where $p_{k}$ is the search direction, taken as the negative gradient direction $p_{k}=-\nabla f(x^{(k)})$, $\nabla$ means gradient operator. $\varepsilon_{k}$ means step size, determined by searching, such that a good $\varepsilon_{k}$ ensures:
	\begin{equation}\label{eq.gradientdescent.result}
		f(x^{(k)} +\varepsilon_{k}p_{k})=\min_{\varepsilon\ge 0}f(x^{(k)} +\varepsilon p_{k}).
	\end{equation}
	\section{Linear N-step method}\label{appendix.B}	
	According to ``Linear N-step method"
	\cite{dahlquistSpecialStabilityProblem1963,Suli_Mayers_2003,griffithsNumericalMethodsOrdinary2010,WU202344},
	the following four results are highlighted:
	\begin{result}\label{result.1}
		A ``Linear N-step method" ${\textstyle \sum_{i=0}^{N}f_{i}x_{k+i}}=\varepsilon {\textstyle \sum_{i=0}^{N}h_{i}j_{k+i}}$ can be checked for 0-stability by determining the roots of its characteristic polynomial $P_{N}(\delta)={\textstyle \sum_{i=0}^{N}f_{i}\delta^{i}}$. If all roots denoted by $\delta$ of the polynomial $P_{N}(\delta)$ satisfy
		$\left | \delta \right | \le 1$ with $\left | \delta \right | = 1$ being simple, then ``Linear N-step method" is 0-stable (i.e., has 0-stability).
	\end{result}
	\begin{result}\label{result.2}
		A ``Linear N-step method" is called consistent of order `` $p$" if truncation error for theoretical solution is of order $\mathcal{O}(\varepsilon^{p})$, where $p>0$.
	\end{result}
	\begin{result}\label{result.3}
		A ``Linear N-step method" is convergent,
		i.e., $x_{\left [ \tau/\delta  \right ] }\to x^{\ast}(\tau)$, for all $\tau \in\left [0,\tau_{a}  \right ]$ as $\varepsilon \to 0$, if and only
		if the method is 0-stable and consistent. That is, 0-stability plus consistency means convergence, which is also called ``Dahlquist equivalence theorem".
	\end{result}
	\begin{result}\label{result.4}
		A linear 0-stable consistent method converges at a rate matching order of its truncation error.
	\end{result}
	\section{Principle of space attack with ``No free lunch"}\label{appendix.C}
	\begin{figure*}[!h]
		\centering  
		\subfigure[]{%
			
			\resizebox{0.4875\textwidth}{!}{%
				
				\begin{tikzpicture}
					
					\draw[thick] (0,0) -- (2,0) -- (3,1) -- (1,1) -- cycle;
					\node at (2.5,1.5) {x};
					
					\draw[thick] (1,1) -- (3,1) -- (4,2) -- (2,2) -- cycle;		
					\node at (1.5,0.5) {y};
					
					\draw[->, thick] (1.5, -0.25) -- (1.5, 1.25);
					\draw[->, thick] (2.5, 0.75) -- (2.5, 2.5);
					
					\draw[->, thick] (6, 1) -- (4, 1);
					\node at (5,1.25) {attack};
					\node at (3,-1) {$\begin{bmatrix}x\\y\end{bmatrix}$};
				\end{tikzpicture}
			}%
			
			\label{fig.attack.structure.a}
		} 
		\hfill 
		\subfigure[]{%
			\resizebox{0.4875\textwidth}{!}{%
				\def\layersep{1.5cm}
				\raisebox{0cm}{
					\begin{tikzpicture}
						\draw[->] (0, 0) -- (0, 2) node[above] {$\text{Re}$};
						\draw[->] (0, 0) -- (2, 0) node[right] {$\text{Im}$};
						
						\draw[->, thick] (0, 0) -- (2, 2) node[midway, above right] {$x +\mathrm{i}y$};
						\draw[->, thick] (6, 0) -- (4, 0);
						\node at (5,0.25) {attack};
						\node at (0, 1){$x$};
						\node at (1, 0){$\mathrm{i}y$};
						\node at (2.75, -1.25) {$\begin{bmatrix}
								x+\mathrm{i}y
							\end{bmatrix}$};
				\end{tikzpicture}}
			}%
			\label{fig.attack.structure.b}
		}  
		\caption{Principle of space attack $\begin{bmatrix}x\\y\end{bmatrix}$ and $\begin{bmatrix}
				x+\mathrm{i}y
			\end{bmatrix}$. But the authors advise viewing them from the perspectives of ``Number theory" and ``Number tables".
			\subref{fig.attack.structure.a} The horizontal vector cannot reach the vertical vector (it can only target specific points and real-world factors have minimal impact).
			\subref{fig.attack.structure.b} The horizontal vector attacks the horizontal vector (vector addition and subtraction on parallel lines is straightforward and real-world factors significantly impact the model).} 
		\label{fig.attack.structure} 
	\end{figure*}
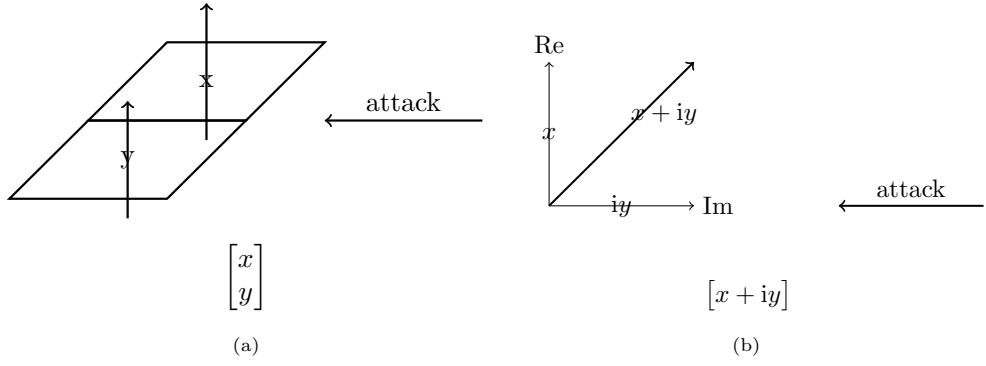
\end{appendices}
\newpage
\bibliography{bib-example}
\newpage

\begin{figure}
	\includegraphics[width=0.13\columnwidth]{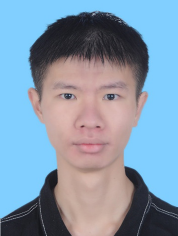}\noindent
	{\bf Jiakuang He} prefers the intersection of mathematics, chemistry, biology, medicine and computing.
	
	Currently, he focuses on precision medicine and ODE neural networks based on complex conjugate matrix equations.
\end{figure}

\begin{figure}
	\includegraphics[width=0.13\columnwidth]{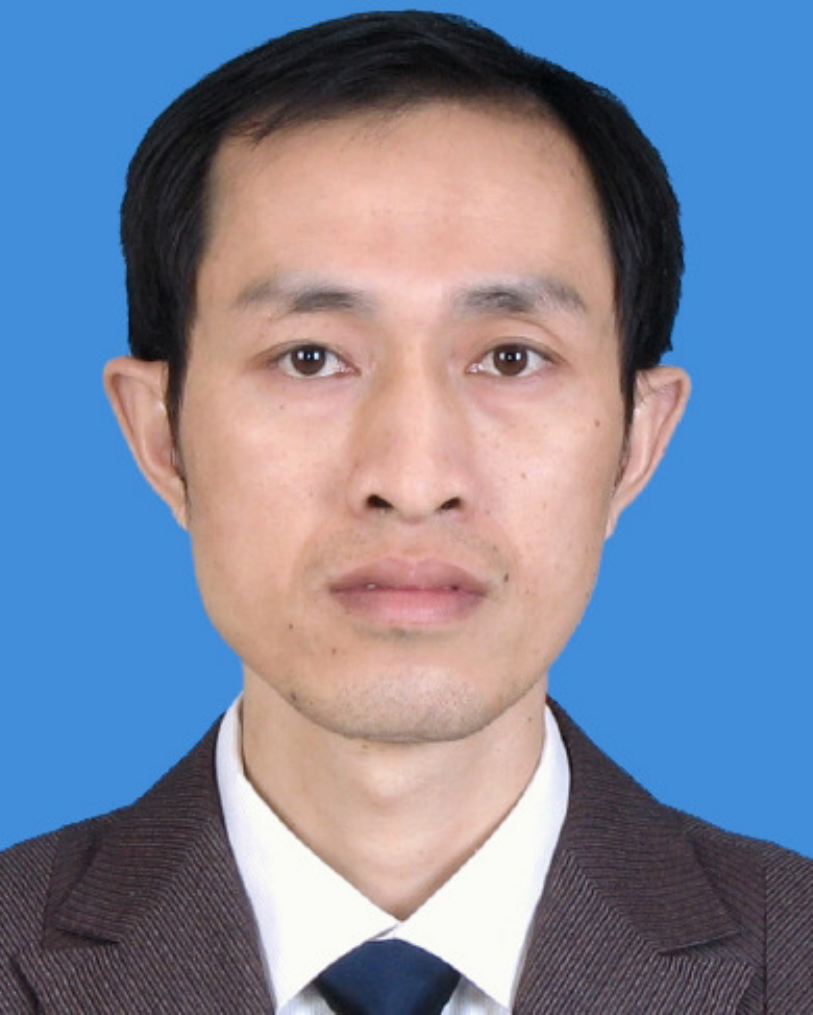}\noindent
	{\bf Dongqing Wu} received the M.S. degree in computer graphics from the Institute of Industrial Design and Graphics, South China University of Technology, Guangzhou, China, in 2005, and the Ph.D. degree in mechanical engineering from the School of Electromechanical Engineering, Guangdong University of Technology, Guangzhou, China, in 2019. He is currently a professor with the School of Mathematics and Data Science, Zhongkai University of Agriculture and Engineering, Guangzhou, China. His current research interests include neural networks, robotics, and numerical analysis.
\end{figure}	

\end{document}